\documentclass[12pt]{article}
\title{Stability for line solitary waves of Zakharov--Kuznetsov equation
\amssubj{35B32, 35B35, 35Q53.}
}
\author{YOHEI YAMAZAKI\footnote{{\it E-mail addresses:} y-youhei@math.kyoto-u.ac.jp} \\ {\footnotesize Department of Mathematics, Kyoto University} \\  {\footnotesize Kitashirakawa-Oiwakecho, Sakyo, Kyoto 606-8502, Japan,}\\ {\footnotesize Osaka City University Advanced Mathematical Institute }\\ {\footnotesize 3-3-138 Sugimoto, Sumiyoshi-ku Osaka 558-8585 Japan }}
\usepackage{amssymb}
\usepackage{mathrsfs}
\usepackage{amsthm}
\usepackage{amsmath}
\usepackage{titlefoot}
\usepackage{fancyhdr}
\def\pdfliteral #1 {}

\numberwithin{equation}{section}

\newtheorem{theorem}{Theorem}[section]
\newtheorem{corollary}[theorem]{Corollary}

\newtheorem{lemma}[theorem]{Lemma}
\newtheorem{proposition}[theorem]{Proposition}

\theoremstyle{definition}
\newtheorem{definition}[theorem]{Definition}
\newtheorem{remark}[theorem]{Remark}

\renewcommand{\eqref}[1]{(\ref{#1})}

\renewcommand{\bigskip}{\vspace{0.3cm}}

\newcommand{\N}{{\mathbb N}}
\newcommand{\R}{{\mathbb R}}
\newcommand{\C}{{\mathbb C}}

\newcommand{\Z}{{\mathbb Z}}
\newcommand{\T}{{\mathbb T}}

\newcommand{\norm}[1]{{\left \lVert #1 \right \rVert}}

\newcommand{\tbr}[1]{\langle #1 \rangle}
\newcommand{\Tbr}[1]{\left\langle #1 \right\rangle}
\newcommand{\RT}{\mathbb {R} \times \mathbb {T}}
\newcommand{\RTL}{\mathbb {R} \times \mathbb {T}_L}

\newcommand{\re}{\mbox{\rm Re }}
\newcommand{\im}{\mbox{\rm Im }}

\date{}

\begin{document}

\maketitle


\begin{abstract}
In this paper, we consider the stability for line solitary waves of the two dimensional Zakharov--Kuznetsov equation on $\RTL$ which is one of a high dimensional generalization of Korteweg--de Vries equation , where $\T_L$ is the torus with the period $2\pi L$.
The orbital and asymptotic stability of the one soliton of Korteweg--de Vries equation on the energy space has been proved by Benjamin \cite{TBB}, Pego and Weinstein \cite{P W 2} and Martel and Merle \cite{M M 1}.
We regard the one soliton of Korteweg--de Vries equation as a line solitary wave of Zakharov--Kuznetsov equation on $\RTL$.
We prove the stability and the transverse instability of the line solitary waves of Zakharov--Kuznetsov equation by applying Evans' function method and the argument of Rousset and Tzvetkov \cite{R T 1}.
Moreover, we prove the asymptotic stability for the orbitally stable line solitary wave of Zakharov--Kuznetsov equation by using  the argument Martel and Merle \cite{M M 1,M M 2,M M 3}, a Liouville type theorem and a corrected virial type estimate.
\end{abstract}


\section{Introduction}

We consider the two dimensional Zakharov--Kuznetsov equation
\begin{equation}\label{ZKeq}
 u_t + \partial_x(\Delta u + u^2)=0, \quad (t,x,y) \in \R\times \RTL,
 \end{equation}
where $\Delta=\partial_x^2 + \partial_y^2$, $u=u(t,x,y)$ is an unknown real-valued function, $\T_L=\R/2\pi L \Z$ and $L>0$.

In \cite{Z K}, Zakharov and Kuznetsov derived the Zakharov--Kuznetsov equation to describe the propagation of ionic-acoustic waves in uniformly magnetized plasma.
In \cite{L L S}, Lannes, Linares and Saut proved the rigorous derivation of the Zakharov--Kuznetsov equation from the Euler--Poisson system for uniformly magnetized plasmas.
The Cauchy problem of the Zakharov--Kuznetsov equation has been studied in the recent years.
In \cite{AVF}, Faminskii proved the global well-posedness of the Zakharov--Kuznetsov in the energy space $H^1(\R^2)$.
This result has been pushed down to $H^s(\R^2)$ for $s>\frac{3}{4}$ by Linares and Pastor \cite{L P 1}.
This result was recently improved by Gr\"unrock and Herr \cite{G H} and Molinet and Pilod \cite{M P} who proved local well-posedness in $H^s(\R^2)$ for $s> \frac{1}{2}$.
In \cite{M P}, Molinet and Pilod showed the global well-posedness of \eqref{ZKeq} in $H^1(\RTL)$.
Moreover, the well-posedness of higher dimensional Zakharov--Kuznetsov equation and the generalized Zakharov--Kuznetsov have been studied by \cite{AG1,L P 1,L P 2,L S,R V}.

The equation \eqref{ZKeq} has the following conservation laws:
\begin{align}\label{moment}
M(u) = \int_{\RTL} |u|^2 dxdy,
\end{align}
\begin{align}\label{energy}
E(u)= \int_{\RTL}\Bigl( \frac{1}{2}|\nabla u|^2 - \frac{1}{3}u^3 \Bigr) dxdy,
\end{align}
where $u \in H^1(\RTL)$.

In this paper, we show the orbital stability and the asymptotic stability of solitary waves of \eqref{ZKeq}.
By a solitary wave, we mean a non-trivial solution of \eqref{ZKeq} with form $u(t,x,y)=Q(x-ct,y)$ , where $c>0$ and $Q \in H^1(\RTL)$ is a solution of 
\begin{align}\label{Seq}
-\Delta Q + c Q - Q^2=0, \quad (x,y) \in \RTL.
\end{align}
We can write the equation \eqref{Seq} as  $S_c'(Q)=0$, where
\[S_c(u)=E(u)+cM(u)\]
and $S_c'$ is the Fr\'echet derivative of $S_c$.

The orbital stability of solitary waves is defined as follows.
\begin{definition}
We say that a solitary wave $Q(x-ct,y)$ is orbitally stable in $H^1(\RTL)$ if for any $\varepsilon >0$ there exists $\delta>0$ such that for all initial data $u_0 \in H^1(\RTL)$ with $\norm{u_0-Q}_{H^1}<\delta$, the solution $u(t)$ of \eqref{ZKeq} with $u(0)=u_0$ exists globally in time and satisfies
\[\sup_{t>0}\inf_{(x_0,y_0) \in \RTL}\norm{u(t,\cdot,\cdot)-Q(\cdot-x_0,\cdot-y_0)}_{H^1}<\varepsilon .\]
Otherwise, we say the solitary wave $Q(x-ct,y)$ is orbitally unstable in $H^1(\RTL)$.
\end{definition}
The orbital stability of positive solitary waves of the generalized Zakharov--Kuznetsov equation on $\R^N$ was showed by de Bouard \cite{AdB} under the assumption of well-posedness on the energy space.
In \cite{C M P S}, C\^ote, Mu\~noz, Pilod and Simpson have proved the asymptotic stability of positive solitary waves and multi-solitary waves of  the Zakharov--Kuznetsov equation on $\R^2$ by adapting the argument of Martel and Merle \cite{M M 1, M M 2, M M 3} to a multidimensional model.

The solution $u$ to \eqref{ZKeq} is not depend on the variable of the transverse direction $\T_L$ if and only if the solution $u$ is a solution to the Korteweg--de Vries equation
\begin{align}\label{KdV}
u_t + u_{xxx} + 2uu_x=0, \quad (t,x) \in \R \times \R.
\end{align}
The Korteweg--de Vries equation describes the propagation of ionic-acoustic waves in unmagnetized plasma.
The equation \eqref{KdV} has the soliton solution $R_c(t,x)=Q_c(x-ct)$, where $Q_c$ is the positive symmetric solution to  
\begin{align}\label{SKdV}
-\partial_x^2 Q +c Q -Q^2 = 0, \quad Q \in H^1(\R).
\end{align}
Here, $Q_c$ has the explicit form 
\[Q_c(x)=\frac{3c}{2}\cosh^{-2}\Bigl( \frac{\sqrt{c}x}{2} \Bigr).\]
The orbital stability of the soliton $R_c$ has been proved by Benjamin \cite{TBB}.
In \cite{P W 2}, Pego and Weinstein have showed the asymptotic stability of the soliton $R_c$ on the exponentially weighted space by investigating a spectral property of linearized operator around $Q_c$.
The argument of Pego and Weinstein \cite{P W 2} is useful to prove the asymptotic stability on the exponentially weighted space for nonintegrable equation.
However, the assumption of the exponential decay of initial data yields that the solution does not have a small soliton other than the main soliton.
To treat solutions with a small soliton other than the main soliton, Mizumachi \cite{TM1} has improved this result, using polynomial weighted spaces.
In \cite{M M 1, M M 2, M M 3}, Martel and Merle have proved the asymptotic stability of the soliton for initial data on $H^1(\R)$.
To prove the asymptotic stability for initial data on $H^1(\R)$, Martel and Merle have showed the Liouville type theorem for the Korteweg--de Vries equation.
The main tool to show the Liouville type theorem is the virial type estimate for solutions with some decay in space.

Then, we regard the soliton solution $R_c$ of \eqref{KdV} as a line solitary wave of \eqref{ZKeq}, namely we define the line solitary wave $\tilde{R}_c$ and the solution $\tilde{Q}_c$ of \eqref{Seq} by 
\[\tilde{R}_c(t,x,y)=\tilde{Q}_c(x-ct,y)=R_c(t,x)=Q_c(x-ct), \quad (t,x,y) \in \R \times \RTL.\]
A natural question concerning $\tilde{R}_c$ is the stability of $\tilde{R}_c$ with respect to perturbations which are periodic in the transversal direction.
The stability of the line solitary wave $\tilde{R}_c$ on Kadomtsev--Petviashvili equation have been studied by many papers.
The stability of $\tilde{R}_c$ on KP-II was confirmed the heuristic analysis by Kadomtsev and Petviashvili \cite{K P}.
In \cite{V A}, Villarroel and Ablowitz have showed the stability of line solitons $\tilde{R}_c$ of KP-II against decaying perturbations by the inverse scattering method.
In \cite{M T} Mizumachi and Tzvetkov have proved the orbital stability and the asymptotic stability of $\tilde{R}_c$ on KP-II in $L^2(\RT)$ by using the B\"ackland transformation.
The asymptotic stability for line solitons $\tilde{R}_c$ of KP-II on $\R^2$ has been proved by Mizumachi \cite{TM2}.
On $\R^2$, because of finite speed propagations of local phase shifts
along the crest of the modulating line soliton for the transverse direction, the line soliton $\tilde{R}_c$ is not orbitally stable in the usual sense.
To prove the asymptotic stability, Mizumachi have showed that the local modulations of the amplitude and the phase shift of line solitons behave like a self-similar solution of the Burgers equation. 
For KP-I equation, Rousset and Tzvetkov have proved the orbital stability and instability for line solitons $\tilde{R}_c$ of KP-I on $\RT$ in \cite{R T 1, R T 3} and on $\R^2$ in \cite{R T 0}.
For Zakharov--Kuznetsov equation, the instability for line solitons $\tilde{R}_c$ on $\R^2$ has been showed by Rousset and Tzvetkov in \cite{R T 0}.
On $\T_{L_1}\times \T_{L_2}$ with sufficiently large $L_2$, the linear instability of line periodic solitary waves of  Zakharov--Kuznetsov equation have been showed by Johnson \cite{MJ} by using Evan's function method.

The one of main results is the following:
\begin{theorem}\label{orbital-stability}
Let $c>0$.
Then, the following holds.
\begin{enumerate}
  \setlength{\parskip}{0.1cm} 
  \setlength{\itemsep}{0.05cm} 
\renewcommand{\labelenumi}{\rm (\roman{enumi})}
\item If $0<L \leq \frac{2}{\sqrt{5c}}$, then $\tilde{R}_c$ is orbitally stable.
\item If $L> \frac{2}{\sqrt{5c}}$, then $\tilde{R}_c$ is orbitally unstable.
\end{enumerate}
\end{theorem}
In Theorem \ref{orbital-stability}, the instability for line solitary waves follows a symmetry breaking bifurcation of line solitary waves in the following proposition.
\begin{proposition}\label{prop-bifurcation}
Let $c_0>0$ and $L=\frac{2}{\sqrt{5c_0}}$.
Then, there exist $\delta_0>0$ and $\varphi_{c_0} \in C^2((-\delta_0,\delta_0)^2,H^2(\RTL))$ such that for $\vec{a}=(a_1,a_2) \in (-\delta_0,\delta_0)^2$ we have $\varphi_{c_0}(\vec{a})>0$, $\varphi_{c_0}(\vec{a})(x,y)=\varphi_{c_0}(\vec{a})(-x,y)$, $\check{c}(\vec{a})=\check{c}(|\vec{a}|,0)$, 
\[-\Delta \varphi_{c_0} (\vec{a})+\check{c} (\vec{a}) \varphi_{c_0}(\vec{a}) - (\varphi_{c_0}(\vec{a}))^2=0,\]
\[\varphi_{c_0}(\vec{a})=\tilde{Q}_{c_0}+a_1 \tilde{Q}_{c_0}^{\frac{3}{2}}\cos \frac{y}{L} +a_2 \tilde{Q}_{c_0}^{\frac{3}{2}}\sin \frac{y}{L} + O(|\vec{a}|^2) \quad \mbox{ as } |\vec{a}| \to 0\]
and
\[\norm{\varphi_{c_0}(\vec{a})}_{L^2(\RTL)}^2=\norm{\tilde{Q}_{c_0}}_{L^2(\RTL)}^2+\frac{C_{2,c_0}}{2} |\vec{a}|^2+o(|\vec{a}|^2) \quad \mbox{ as } |\vec{a}| \to 0,\]
where $\check{c} (\vec{a})=c_0+\frac{\check{c} ''(0)}{2}|\vec{a}|^2+o(|\vec{a}|^2)$ as $|\vec{a}| \to 0$, $\check{c} ''(0)>0$ and 
\[C_{2,c_0}=\frac{3\check{c} ''(0)\norm{\tilde{Q}_{c_0}}_{L^2}^2}{2c_{0}} -\frac{5\norm{\tilde{Q}_{c_0}^{\frac{3}{2}}\cos \frac{y}{L}}_{L^2}^2}{2}>0\]
\end{proposition}
\begin{remark}
Proposition \ref{prop-bifurcation} follows Proposition 1 and the proof of Theorem 1.3 in \cite{YY2}.
The positivity of the constant $\check{c}''(0)$ follows the relation $L^4\check{c}''(0)=\omega''(0)$ and the positivity of $\omega''(0)$ in the proof of Theorem 1.3 in \cite{YY2}, where $\omega''(0)$ have been defined in Proposition 1 in \cite{YY2}.
The positivity of the constant $C_{2,c_0}$ have been proved from the inequality (2.25)
\begin{align*}
 R(p)\leq & \frac{4(p+1)(p^6+18p^5-11p^4-130p^3+13p^2+16p-3)}{(5-p)(p+3)^2(5p-1)(3p+1)(p-1)}\\
 &+ \frac{32p^3(p+1)^4(3p-1)}{3(7p-3)(5p-1)(3p+1)(p+3)^3(p-1)}
\end{align*}
in \cite{YY2} and the relation 
\[-\frac{6L^4C_{2,c_0}\norm{\tilde{Q}_{c_0}}_{L^2(\RTL)}^2}{5c_0}=R(2)<0,\]
where $R(p)$ is defined in the proof of Theorem 1.3 in \cite{YY2}.
\end{remark}

We define a semi-norm $\norm{\cdot}_{H^1(x>a)}$ on $H^1(\RTL)$ by
\[ \norm{u}_{H^1(x>a)}^2=\int_{x>a} (|\nabla u(x,y)|^2+|u(x,y)|^2) dxdy, \quad u \in H^1(\RTL).\]
The following theorem is a main theorem for the asymptotic stability.
\begin{theorem}\label{asymptotic-stability}
Let $c_0>0$.
\begin{enumerate}
  \setlength{\parskip}{0.1cm} 
  \setlength{\itemsep}{0.05cm} 
\renewcommand{\labelenumi}{\rm (\roman{enumi})}
\item If $0<L < \frac{2}{\sqrt{5c_0}}$, then the following holds.
For any $\beta>0$, there exists $\varepsilon _{L,\beta}>0$ such that for $u_0 \in H^1(\RTL)$ with $\norm{u_0-\tilde{Q}_{c_0}}_{H^1} < \varepsilon _{L,\beta} $, there exist  $\rho(t) \in C^1([0,\infty),\R)$ and $c_+>0$ satisfying that  
\[\norm{u(t,\cdot,\cdot)-\tilde{Q}_{c+}(\cdot-\rho(t),\cdot)}_{ H^1(x>\beta t)} \to 0 \mbox{ as } t \to \infty,\]
\[\dot{\rho}(t)-c_+ \to 0 \mbox{ as } t \to \infty\]
and $|c_0-c_+| \lesssim \norm{u_0-\tilde{Q}_{c_0}}_{H^1} $, where $u$ is the unique solution of $\eqref{ZKeq}$ with $u(0)=u_0$.
\item If $L= \frac{2}{\sqrt{5c_0}}$, then the following holds.
For any $\beta>0$, there exists $\varepsilon _{\beta}>0$ such that for $u_0 \in H^1(\RTL)$ with $\norm{u_0-\tilde{Q}_{c_0}}_{H^1} < \varepsilon _{\beta} $, there exist $\rho_1(t),\rho_2(t) \in C^1([0,\infty),\R)$, $c_+>0$ and $\vec{a}_+ \in \R^2$ satisfying that  
\[\norm{u(t,\cdot,\cdot)-\Theta (\vec{a}_+,c_+)(\cdot-\rho_1(t),\cdot-\rho_2(t))}_{H^1(x>\beta t)} \to 0  \mbox{ as } t \to \infty,\]
\[\dot{\rho}_1(t)-\hat{c}_+ \to 0, \dot{\rho}_2(t) \to 0 \mbox{ as } t \to \infty,\]
\[|c_+-c_0||\vec{a}_+|=0\]
and $|c_0-c_+|+ |\vec{a}_+|^2 \lesssim \norm{u_0-\tilde{Q}_{c_0}}_{H^1} $, where
\begin{align}\label{def-hatc}
\hat{c}_+=
\begin{cases}
c_+, & \vec{a}_+=(0,0),\\
\check{c} (\vec{a}_+), & c_+=c_0,
\end{cases}
\end{align}
\[\Theta (\vec{a},c)(x,y)=\frac{c}{c_0}\varphi_{c_0}(\vec{a})\Bigl(\sqrt{\frac{c}{c_0}}x,y\Bigr),\]
and $u$ is the unique solution of $\eqref{ZKeq}$ with $u(0)=u_0$.
\end{enumerate}
\end{theorem}

\begin{remark}
Since a neighborhood of $\tilde{Q}_{c_0} $ in $H^{1}(\RT_{L})$ contains the branch corresponding to unstable line solitary waves in the case $L=\frac{2}{\sqrt{5c_0}}$, Theorem \ref{asymptotic-stability} shows that solutions away from unstable solitary waves approach one of solitary waves in the neighborhood $\tilde{Q}_{c_0}$ as $t \to \infty$ in the sense of the norm $H^{1}(x>\beta t)$.
\end{remark}
\begin{remark}
In Theorem \ref{asymptotic-stability}, the unique solution $u$ of $\eqref{ZKeq}$ with $u(0)=u_0$ means that for $T>0$ the function $u|_{[-T, T]}$ is a unique solution of $\eqref{ZKeq}$ with $u(0)=u_0$ in $ C([-T,T],H^1(\RTL))\cap X_T^{1,\frac{1}{2}+}$ which is defined in \cite{M P}.
\end{remark}
Let us now explain the argument to prove Theorem \ref{orbital-stability}.
Since the solution $\tilde{Q}_c$ of \eqref{Seq} is not a minimizer of the functional $S_c(u)$ on $\{ u \in H^1; M(u)=M(\tilde{Q}_c)\}$ for general $c>0$, we can not apply the variational argument to prove the orbital stability.
Therefore, to prove the orbital stability of $\tilde{Q}_{c_0}$ we use the argument in \cite{G S S 1,MIW2} for $0<L<\frac{2}{\sqrt{5c_0}}$.
In the case $L=\frac{2}{\sqrt{5c_0}}$, the linearized operator of \eqref{Seq} around $\tilde{Q}_{c_0}$ has an extra eigenfunction corresponding to zero eigenvalue.
Thus, we can not show the orbital stability of $\tilde{Q}_{c_0}$ by using the standard argument in \cite{G S S 1,G S S 2,MIW2}.
Since any neighborhood of $\tilde{Q}_{c_0}$ contains the two branches which are comprised of line solitary waves $\tilde{Q}_c$ and solitary waves $\varphi_{c_0}(\vec{a})$, we can not apply the argument for the linearized operator of the evolution equation with an extra eigenfunction by Comech and Pelinovsky \cite{C P} and Maeda \cite{MM}.
Because of the degeneracy of the third order term of Lyapunov functional, we can not use the argument for the instability of a standing wave on a point of interaction of two branches of standing waves in Ohta \cite{MO}.
To prove the stability of $\tilde{Q}_{c_0}$, we apply the argument in \cite{YY2,YY3}.

To show the nonlinear instability of $\tilde{Q}_c$ from the existence of an unstable mode of the linearized operator around $\tilde{Q}_c$, we apply the argument by Grenier \cite{EG} and  Rousset and Tzvetkov \cite{R T 1}.
Since the simple criterion in \cite{R T 0, R T 2} does not seem be applicable to the linearized operator of \eqref{ZKeq} around $\tilde{Q}_c$, it is difficult to get the existence of an unstable mode of the linearized operator by the implicit function theorem.
For sufficiently large $L$, Bridges \cite{JTB} have showed the existence of an unstable mode by sophisticated arguments.
To get the existence of an unstable mode of linearized operator for all $L>\frac{2}{\sqrt{5c_0}}$, we use Evans' function method in Pego and Weinstein \cite{P W} for gKdV equation.

Next we explain the main ideas and difficulties in the proof of Theorem \ref{asymptotic-stability}.
Since the equation \eqref{ZKeq} is not complete integrable, we can not use the inverse scattering method to get the asymptotic behavior of solutions.
To prove the asymptotic stability, we apply the argument by Martel and Merle \cite{M M 1,M M 2,M M 3} and C\^ote et al. \cite{C M P S}.
This argument relies on a Liouville type theorem for decaying solutions around a solitary wave.
From the orbital stability and the monotonicity property, solutions near by a solitary wave converge to an exponentially decaying function in $H^1(x>a)$ up to subsequence of time.
Due to the Liouville type theorem, this function has to be solitary waves.
The main tool to prove Liouville type theorem is the virial type estimate.
In the case $0<L<\frac{2}{\sqrt{5c_0}}$, the linearized operator of \eqref{Seq} around $\tilde{Q}_{c_0}$ is coercive on $\{u \in H^1;M(u)=M(\tilde{Q}_{c_0})\}$ by modulating translation.
Thus, applying the estimate of \cite{M M 3} we can show the virial type estimate.
However, in the case $L=\frac{2}{\sqrt{5c_0}}$, the linearized operator of \eqref{Seq} around $\tilde{Q}_{c_0}$ is not coercive on the function space with the standard orthogonal condition.
To get the coerciveness of linearized operator, we estimate the difference between the solution and $\Theta$ instead of the difference between the solution and solitary waves, where $\Theta$ is defined in Theorem \ref{asymptotic-stability}.
However, since $\Theta$ is not solution of stationary equation, a term including $S_c'(\Theta)$ appears in the virial type estimate.
Therefore, we can not get the coerciveness of the virial type estimate by the argument in \cite{M M 3}.
To treat the term with $S_c'(\Theta)$, we investigate the virial type estimate with a correction term $S_{\hat{c}}'(\Theta)$, where $\hat{c}$ is the suitable propagation speed of $\Theta$.
To get the coerciveness of the virial type estimate with a correction, we use the precise estimate for a quadratic form and interactions among main terms.
Due to this virial type estimate with the correction, we get the Liouville type theorem around the bifurcation point $\tilde{Q}_{c_0}$.

Our plan of the present paper is as follows.
In Section 2, we show the well-posedness result on weighted space to prove the monotonicity property.
The argument of this well-posedness result follows the argument by Kato in \cite{TK1}. 
In Section 3, we prove the properties of the linearized operator of \ref{ZKeq} and the estimate of semi-group corresponding to the linearized operator.
To show the linear instability of linearized equation, we use the argument by Pego and Weinstein \cite{P W}.
In Section 4, we prove (ii) of Theorem \ref{orbital-stability} by the argument of Rousset and Tzvetkov \cite{R T 1}.
In Section 5, we show (i) of Theorem \ref{orbital-stability} by the argument of \cite{G S S 1} and \cite{YY2,YY3}.
In Section 6, we prove the coercive type estimate of a quadratic form and the Liouville property for orbitally stable solitary waves.
To get the monotonicity property, we use the Kato type local smoothing effect in Section 2.
In Section 7, we prove Theorem \ref{asymptotic-stability} by applying the Liouville property and the monotonicity property in Section 6.

\section{Preliminaries}
In this section, we show the regularity of solutions to \eqref{ZKeq} on weighted space.
To proof of smoothness of solutions to \eqref{ZKeq} on weighted space, we apply the argument on KdV in \cite{TK1}. 
For $u \in L^2(\RTL)$ we define $\hat{u}$ by the space-time Fourier transform of $u$.

From the result on well-posedness in $H^1(\RTL)$ by Molinet-Pilod \cite{M P}, for initial data $u_0 \in H^1(\RTL)$ there exists the unique solution $u(t)$ of \eqref{ZKeq} such that $u(0)=u_0$ and for $T>0$
\[ u|_{[-T,T]} \in C([-T,T],H^1(\RTL))\cap X_T^{1,\frac{1}{2}+}.\]
Moreover, for any $T>0$, there exists a neighborhood $\mathcal{U}$ of $u_0$ in $H^1(\RTL)$, such that the flow map data-solution 
\[v_0 \in \mathcal{U} \mapsto v \in C([0,T),H^1(\RTL)) \cap X_T^{1,\frac{1}{2}+}\]
is smooth.
Here,  the function space $X^{1,\frac{1}{2}+}$ is defined in \cite{M P}.
In this paper, we define $H^1$-solution by the solution in the function space $C([0,\infty),H^1(\RTL))$ satisfying the conservation laws $M(u(t))=M(u(0))$ and $E(u(t))=E(u(0))$.

Let $U_b(t) = \exp(-t(\partial_x-b)((\partial_x-b)^2+\partial_y^2))$ for $b>0$.
Then, we have for $ u \in L^2(\RTL)$ if $e^{bx}u \in L^2(\RTL)$ then
\[e^{-bx} U_b(t) e^{bx} u = U_0(t)u.\]
The following lemma decay properties of the propagator $U_b$.
\begin{lemma}\label{lem-decay-prop}
Let $b>0$, $s ,s'\in \R$, $s<s'$ and $n \in \Z_+$.
Then, there exists $C=C(n,s,b)>0$ such that for $u \in H^s(\RTL)$, $0\leq j \leq n$ and $t>0$
\begin{align}
\norm{U_b(t)u}_{H^{s'}} &\leq C t^{-\frac{s'-s}{2}}e^{b^3t}\norm{u}_{H^{s'}} \label{lem-decay-1}\\
\norm{\partial_x^j \partial_y^{n-j} U_b(t)u}_{L^2} &\leq C t^{-\frac{n}{2}} e^{b^3t} \norm{u}_{L^2} \label{lem-decay-2}\\
\norm{\partial_t U_b(t)u}_{L^2} &\leq C t^{-\frac{3}{2}} e^{b^3t} \norm{u}_{L^2} \label{lem-decay-3}
\end{align}
\end{lemma}
\proof
By the factorization we have
\[U_b(t)=\exp(tb^3)\exp(-3tb^2\partial_x)\exp(tb(3\partial_x^2+\partial_y^2))\exp(-t\partial_x\Delta).\]
Since $\exp(-3tb^2\partial_x)$ and $\exp(-t\partial_x\Delta)$ are unitary in $H^s$ and $\exp(tb(3\partial_x^2+\partial_y^2))$ is the heat semigroup, we have the estimates \eqref{lem-decay-1}--\eqref{lem-decay-3}.

\qed



\begin{proposition}\label{prop-WPW}
Let $u $ be a $H^1$-solution to $\eqref{ZKeq}$ with $e^{bx}u(0) \in L^2(\RTL)$ for some $b>0$.
Then we have $e^{bx}u \in C([0,\infty),L^2(\RTL)) \cap C^{\infty}((0,\infty),H^{\infty}(\RTL))$, with
\begin{align}
\norm{e^{bx}u(t)}_{L^2} \leq e^{Kt}\norm{e^{bx}u(0)}_{L^2}
\end{align}
where $K$ denotes various constant depending only on $b$ and $\norm{u(0)}_{L^2}$.
Moreover, for any $T>0$ and $s \geq 0$,
\begin{align}
\norm{e^{bx}u(t)}_{H^s}&\leq K' t^{-\frac{s}{2}}, \quad 0<t\leq T, \label{eq-0-2-1}\\
\norm{e^{bx}(\partial_t)^nu(t)}_{H^s}&\leq K' t^{-\frac{s+3n}{2}}, \quad 0<t\leq T , n \in \Z_+ \label{eq-0-2-2}
\end{align}
where $K'$ depend on $s,n,T,b, \norm{e^{bx}u(0)}_{L^2}$, $M(u(0))$ and $E(u(0))$.
\end{proposition}
\proof


Let 
\[ q(x)=e^{bx}(1+\varepsilon e^{2bx})^{-\frac{1}{2}}, \quad r(x)=e^{bx}(1+\varepsilon e^{2bx})^{-1}, \quad p(x)=q(x)^2.\]
Then, we have $q,r,p \in L^{\infty}(\RTL)$ and
\[\partial_xp=2br^2, \quad |\partial_x^2p|\leq 4b^2r^2, \quad |\partial_x^3p|\leq 12 b^3r^2, \quad |\partial_xr|\leq br.\]
Therefore, we have
\begin{align}
\frac{d}{dt}(pu,u)_{L^2}
\leq & -2b(3\norm{r\partial_x u}_{L^2}^2 + \norm{r\partial_yu}_{L^2}^2)+ 12b^2\norm{ru}_{L^2}^2+ \frac{8b}{3}(r^2u,u^2)_{L^2}.\label{prop-WPW-1}
\end{align}
Then,
\begin{align*}
(r^2u,u^2)_{L^2} 
\leq  \frac{1}{2} \norm{r\nabla u}_{L^2}^2+K_0\norm{ru}_{L^2}^2,
\end{align*}
where $K_0$ depend only $b$ and $\norm{u(0)}_{L^2}$.
From \eqref{prop-WPW-1} and $r<q$, we obtain 
\[\frac{d}{dt}\norm{qu}_{L^2}^2\leq -\frac{2b}{3}\norm{r\nabla u}_{L^2}^2 + K \norm{ru}_{L^2}^2 \leq -\frac{2b}{3} \norm{r\nabla u}_{L^2}^2 + K\norm{qu}_{L^2}^2.\]
It follows that $\norm{qu(t)}_{L^2}\leq e^{Kt}\norm{qu(0)}_{L^2}$.
Applying the monotone convergence theorem, we obtain that
\[\norm{e^{bx}u(t)}_{L^2} \leq e^{Kt} \norm{e^{bx}u(0)}_{L^2}, \quad t \geq 0,\]
where $K$ depend only $b $ and $\norm{u(0)}_{L^2}$.
Since $e^{bx}U_0(t)=U_b(t)e^{bx}$, by Lemma \ref{lem-decay-prop} we have for $t>0$
\begin{align*}
\norm{e^{bx}u(t)}_{L^2}\leq & \norm{U_b(t)e^{bx}u(0)}_{L^2} + \int_0^t\norm{U_b(t-\tau) e^{bx} \partial_x(u(\tau)^2)}_{L^2} d\tau\\
\leq &  C e^{b^3 t}  \norm{e^{bx}u(0)}_{L^2} + \int_0^tC(M(u),E(u)) (t-\tau)^{-3/4} \norm{e^{bx}u(\tau)}_{L^2}d\tau.
\end{align*}
Here, we use $\norm{u(t)}_{H^1} \leq C(M(u),E(u))$.
Therefore, $e^{bx}u(t) \in C([0,\infty),L^2(\RTL))$.

Next we show \eqref{eq-0-2-1}.
By the Sobolev embedding and H\"older inequality, we have
\begin{align*}
\norm{e^{bx}\partial_x(u^2)}_{H^{-\frac{1}{4}}}\leq & b \norm{e^{bx}u^2}_{H^{-\frac{1}{4}}}+\norm{e^{bx}u^2}_{H^{\frac{3}{4}}}\\
\leq & C(M(u),E(u)) \norm{e^{bx}u}_{H^1}.
\end{align*}
Thus,
\begin{align*}
\norm{e^{bx}u(t)}_{H^{\frac{5}{4}}} \leq & C t^{-\frac{5}{8}}e^{b^3T}\norm{e^{bx}u(0)}_{L^2}+ K'\int_0^{t}(t-\tau)^{-\frac{3}{4}}\norm{e^{bx}u(\tau)}_{H^1}d\tau\\
\leq &C t^{-\frac{5}{8}}e^{b^3T}\norm{e^{bx}u(0)}_{L^2}+ K'\int_0^{t}(t-\tau)^{-\frac{3}{4}}\norm{e^{bx}u(\tau)}_{H^{\frac{5}{4}}}d\tau,
\end{align*}
where $K'$ depends only on $s,n,T,b, \norm{e^{bx}u(0)}_{L^2}$, $M(u(0))$ and $E(u(0))$.
Therefore, $e^{bx}u \in C((0,\infty),H^{\frac{5}{4}}(\RTL))$ and \eqref{eq-0-2-1} holds for $s=\frac{5}{4}$.
By the interpolation, we obtain \eqref{eq-0-2-1} for $0\leq s\leq \frac{5}{4}$.
To prove $s>\frac{5}{4}$, we use the induction on $s$.
Suppose \eqref{eq-0-2-1} has been proved for $0\leq s \leq s'-\frac{1}{2}$, where $ s' \geq \frac{7}{4}$.
We shall show \eqref{eq-0-2-1} for $0 \leq s \leq s'$.
By Duhamel formula, we have
\[ t^{\frac{s}{2}}e^{bx}u(t)=\int_0^tU_b(t-\tau)\Bigl( \frac{s}{2}\tau^{\frac{s}{2}-1}e^{bx}u(\tau)-\tau^{\frac{s}{2}}e^{bx}\partial_x(u(\tau)^2) \Bigr)d\tau.\]
Since $\norm{U_b(t-\tau)}_{H^{s'-\frac{3}{2}} \to H^{s'}} \leq C (t-\tau)^{-\frac{3}{4}}$, 
\begin{align}\label{eq-0-2-3}
t^{\frac{s'}{2}}\norm{e^{bx}u(t)}_{H^{s'}} \leq C \int_0^t (t-\tau)^{-\frac{3}{4}}\Bigl( \frac{s'}{2}\tau^{\frac{s'}{2}-1} \norm{e^{bx}u(\tau)}_{H^{s'-\frac{3}{2}}}+\tau^{\frac{s'}{2}}\norm{e^{bx}\partial_x(u(\tau)^2)}_{H^{s'-\frac{3}{2}}} \Bigr)d\tau.
\end{align}
From the assumption of the induction we have
\[\tau^{\frac{s'}{2}-1}\norm{e^{bx}u(\tau)}_{H^{s'-\frac{3}{2}}}\leq K' \tau^{\frac{s'}{2}-1-\frac{s'}{2}+\frac{3}{4}}=K'\tau^{-\frac{1}{4}}.\]
On the other hand, by Appendix A in \cite{TK1} for $f, g\in H^s (\RTL) (s \geq \frac{5}{4})$
\[\norm{fg}_{H^s} \lesssim \norm{f}_{H^{\frac{3}{4}}}^{\frac{1}{2}} \norm{f}_{H^{\frac{5}{4}}}^{\frac{1}{2}}\norm{g}_{H^s}+\norm{g}_{H^{\frac{3}{4}}}^{\frac{1}{2}}\norm{g}_{H^{\frac{5}{4}}}^{\frac{1}{2}}\norm{f}_{H^s}.\]
Thus, we have
\[\norm{e^{bx}\partial_x(u(\tau)^2)}_{H^{s'-\frac{3}{4}}}\lesssim \norm{e^{bx}u(\tau)^2}_{H^{s'-\frac{1}{2}}} \lesssim \norm{e^{\frac{bx}{2}}u(\tau)}_{H^{\frac{3}{4}}}^{\frac{1}{2}}\norm{e^{\frac{bx}{2}}u(\tau)}_{H^{\frac{5}{4}}}^{\frac{1}{2}}\norm{e^{\frac{bx}{2}}u(\tau)}_{H^{s'-\frac{1}{2}}}.\]
From the assumption of the induction, we obtain 
\[ \norm{e^{bx}\partial_x(u(\tau)^2)}_{H^{s'-\frac{3}{4}}} \leq K'_{\frac{b}{2}} r^{-\frac{s'}{2}},\]
where $K'_{\frac{b}{2}}$ depends only  on $s,n,T,b, \norm{e^{\frac{bx}{2}}u(0)}_{L^2}$, $M(u(0))$ and $E(u(0))$.
Since 
\[\norm{e^{\frac{bx}{2}}u(0)}_{L^2}\leq (\norm{u(0)}_{L^2}\norm{e^{bx}u(0)}_{L^2})^{\frac{1}{2}},\]
 $K'_{\frac{b}{2}}$ depends only  on $s,n,T,b, \norm{e^{bx}u(0)}_{L^2}$, $M(u(0))$ and $E(u(0))$.
From \eqref{eq-0-2-3} we obtain
\[\norm{e^{bx}u(t)}_{H^{s'}}\leq K'' t^{-\frac{s'}{2}}\]
where $K''$ depends only on $s,n,T,b, \norm{e^{bx}u(0)}_{L^2}$, $M(u(0))$ and $E(u(0))$.
This proves \eqref{eq-0-2-1} for $0 \leq s \leq s'$, completing the induction.

Finally we prove \eqref{eq-0-2-2} by induction on $n$.
For the case $n=0$, it is known by \eqref{eq-0-2-1}.
Assuming that it has been proved for all $s\geq 0$ up to a given $n$, we prove it for $n+1$.
By the induction hypothesis, 
\begin{align}\label{eq-0-2-4}
\norm{\partial_t^n\partial_x\Delta(e^{bx}u)}_{H^{s}}\lesssim \norm{\partial_t^n(e^{bx}u)}_{H^{s+3}} \leq K' t^{-\frac{s+3+3n}{2}}.
\end{align}
On the other hand, 
\[\norm{\partial_t^ne^{bx}\partial_x(u^2)}_{H^s} \lesssim \norm{\partial_t^ne^{bx}u^2}_{H^{s+1}} \lesssim \sum_{j=0}^{n}\norm{e^{bx}(\partial_t^ju)(\partial_t^{n-j}u)}_{H^{s+1}}.\]
By Appendix A in \cite{TK1},
\begin{align*}
\norm{e^{bx}(\partial_t^ju)(\partial_t^{n-j}u)}_{H^{s+1}} \lesssim & \norm{e^{\frac{bx}{2}}\partial_t^ju}_{H^{s+1}} \norm{e^{\frac{bx}{2}}\partial_t^{n-j}u}_{H^{\frac{3}{4}}}\norm{e^{\frac{bx}{2}}\partial_t^{n-j}u}_{H^{\frac{5}{4}}}\\
&+\norm{e^{\frac{bx}{2}}\partial_t^{n-j}u}_{H^{s+1}} \norm{e^{\frac{bx}{2}}\partial_t^{j}u}_{H^{\frac{3}{4}}}\norm{e^{\frac{bx}{2}}\partial_t^{j}u}_{H^{\frac{5}{4}}}.
\end{align*}
Therefore,
\begin{align}\label{eq-0-2-5}
\norm{\partial_t^ne^{bx}\partial_x(u^2)}_{H^s} \leq K' t^{-\frac{s+3+3n}{2}}.
\end{align}
From \eqref{eq-0-2-4} and \eqref{eq-0-2-5} we obtain \eqref{eq-0-2-2} for $n+1$, completing the induction.
$e^{bx}u \in C^{\infty}((0,\infty), H^{\infty}(\RTL))$ follows the estimate \eqref{eq-0-2-2}.

\qed

\section{Linearized operator}

In this section, we show the properties of the linearized operator of \eqref{ZKeq} around $\tilde{R}_c$.
We define the linearized operator $\mathbb{L}_c$ of \eqref{Seq} around $\tilde{Q}_c$ by
\[\mathbb{L}_c=S''_c(\tilde{Q}_c)=-\Delta + c -2 \tilde{Q}_c\]
and the linearized operator $\mathcal{L}_c$ of \eqref{SKdV} around $Q_c$ by
\[\mathcal{L}_c=-\partial_x^2 + c -2Q_c.\]
Then, the linearized operator of \eqref{ZKeq} around $\tilde{R}_c$ is $\partial_x \mathbb{L}_c$.
From Theorem 3.4 in \cite{C G N T}, $\mathcal{L}_c$ has the only one negative eigenvalue 
\[-\lambda_c=-\frac{5c}{4}\]
and an eigenfunction $(Q_c)^{\frac{3}{2}}$ corresponding to $-\lambda_c$.

\begin{proposition}\label{prop-properties-LO}
Let $c>0$.
\begin{enumerate}
  \setlength{\parskip}{0.1cm} 
  \setlength{\itemsep}{0.05cm} 
\renewcommand{\labelenumi}{\rm (\roman{enumi})}
\item If $0<L\leq \frac{2}{\sqrt{5c}}$, then $\partial_x \mathbb{L}_c$ has no eigenvalues with a positive real part.
\item If $0<L< \frac{2}{\sqrt{5c}}$, then 
\[\mbox{\rm Ker}( \mathbb{L}_c) = \mbox{\rm Span}\{\partial_x \tilde{Q}_c\}.\]
\item If $L=\frac{2}{\sqrt{5c}}$, 
\[\mbox{\rm Ker}(\mathbb{L}_c)=\mbox{\rm Span}\Bigl\{\partial_x \tilde{Q}_c, (\tilde{Q}_c)^{\frac{3}{2}} \cos \frac{y}{L}, (\tilde{Q}_c)^{\frac{3}{2}} \sin \frac{y}{L}\Bigr\}.\]
\item If $L>\frac{2}{\sqrt{5c}}$, then $\partial_x \mathbb{L}_c$ has a positive eigenvalue and the number of eigenvalue of $\partial_x \mathbb{L}_c$ with a positive real part is finite.
\end{enumerate}
Here, $\mbox{\rm Span}\{u_1,\dots,u_n\}$ is the vector space spanned by vectors $u_1,\dots,u_n$.
\end{proposition}

\proof
By the Fourier expansion, we have for $u \in H^1(\RTL)$ 
\begin{align}\label{eq-2-1}
( \mathbb{L}_cu)(x,y)=\sum_{n=-\infty}^{\infty}\Bigl(\mathcal{L}_c+\frac{n^2}{L^2}\Bigr)u_n(x)e^{\frac{iny}{L}},
 \end{align}
where
\[u(x,y)=\sum_{n=-\infty}^{\infty}u_n(x)e^{\frac{iny}{L}}.\]
From the equation \eqref{eq-2-1}, we obtain that $\partial_x\mathbb{L}_c$ has an eigenvalue $\lambda$ if and only if there exists $n \in \Z$ such that $\partial_x(\mathcal{L}_c+n^2/L^2)$ has an eigenvalue $\lambda$.
By Theorem 3.4 in \cite{P W}, the essential spectrum of $\partial_x\mathcal{L}_c$ is the imaginary axis.
Moreover, from Theorem 3.1 in \cite{P W}, the number of eigenvalues of $\partial_x(\mathcal{L}_c+n^2/L^2)$ with a positive real part less than or equal to the number of negative eigenvalues of $\mathcal{L}_c+n^2/L^2$.
In the case  $L\leq \frac{2}{\sqrt{5c}}$, since $n^2/L^2\geq \lambda_c$ for all $n\neq 0$, $\mathcal{L}_c+n^2/L^2$ has no negative eigenvalues and (i) is verified.
The kernel of $\partial_x(\mathcal{L}_c+n^2/L^2)$ is trivial if and only if the kernel of $\mathcal{L}_c+n^2/L^2$ is trivial.
Therefore, for $L> \frac{2}{\sqrt{5c}}$ the kernel of $\partial_x\mathbb{L}_c$ is spanned by $\partial_x \tilde{Q}_c$.
In the case $L= \frac{2}{\sqrt{5c}}$, the kernel of $\partial_x\mathbb{L}_c$ is spanned by $\partial_x \tilde{Q}_c$, $(Q_c)^{\frac{3}{2}}\cos \frac{y}{L}$ and $(Q_c)^{\frac{3}{2}}\sin \frac{y}{L}$.
Thus, (ii) and (iii) are verified.

To prove (iv), we apply Evans' function method in \cite{P W}.
We consider the following equation:
\begin{align}\label{eq-2-2}
\partial_x\Bigl(\mathcal{L}_c+ a \Bigl)u-\lambda u=0.
\end{align}
The equation \eqref{eq-2-2} is equivalent to the first order system
\begin{align}\label{eq-2-3}
\partial_x \vec{u} = A(a,\lambda,x)\vec{u}
\end{align}
where 
\[
\vec{u}=
\begin{pmatrix}
u\\
\partial_x u\\
\partial_x^2 u
\end{pmatrix}, \quad
A(a,\lambda,x)=
\begin{pmatrix}
0 & 1 & 0 \\
0 &  0 & 1\\
-2\partial_x Q_c(x) -\lambda & c+a-2Q_c(x) & 0
\end{pmatrix}.
\]
First, we show that $A(a,\lambda,x)$ satisfies the assumption {\bf H1}, {\bf H2}, {\bf H3} and {\bf H4} in Section 1 of \cite{P W}.
Then, $A(a,\lambda,x)$ is analytic in $\lambda$ and $a$ for each $x$, so {\bf H1} holds true.
Let
\[A_{\infty}(a,\lambda)=
\begin{pmatrix}
0 & 1 &0 \\
0 & 0 & 1\\
-\lambda & c+a &0
\end{pmatrix}.\]
Then, $\lim_{|x| \to \infty}A(a,\lambda,x)=A_{\infty}(a,\lambda)$ and $A(a,\lambda,x)$ satisfies {\bf H2} and {\bf H4}.
We define
\[\mu_1(a,\lambda):=\inf\{\re \mu;  \mu \mbox{ is an eigenvalue of } A_{\infty}(a,\lambda)\},\]
\[\mu_2(a, \lambda):=\inf\{\re \mu; \re \mu> \re \mu_1(a, \lambda), \mu \mbox{ is an eigenvalue of } A_{\infty}(a,\lambda)\}.\]
Let 
\[J=\{(a,\lambda) \in  \C^2; A_{\infty}(a,\lambda) \mbox{ has some purely imaginary eigenlvalue}\}.\]
We define $J_+$ be the connected component of $\C^2 \setminus J$ which contains $\{a\geq 0\}\times \{\lambda>0\}$.
Form the perturbation theory of matrices,  the number of eigenvalues counting multiplicity of $A_{\infty}(a,\lambda)$ having the negative real part is constant for $(a,\lambda) \in J_+$.
Since $A_{\infty}(0,\lambda)$ has the only one simple negative eigenvalue for $\lambda>0$, the number of eigenvalues counting multiplicity of $A_{\infty}(a,\lambda)$ having the negative real part is 1 for $(a,\lambda) \in J_+$.
Therefore, for $(a,\lambda) \in J_+$
\[ \mu_1(a,\lambda)<0<\mu_2(a,\lambda)\]
Moreover, for $a>-c/2$ 
\begin{align*}
\mu_1(a,0)<0\leq \mu_2(a,0).
\end{align*}
By the perturbation theory of matrices, there exists a domain $ \tilde{\Omega}$ in $\C^2$ such that $\{a \geq 0\} \times \{\lambda\geq 0\} \subset \tilde{\Omega}$ for $(a,\lambda) \in \tilde{\Omega}$ and $A_{\infty}(a,\lambda)$ has the unique eigenvalue with the smallest real part $\mu_1(a,\lambda)$, which is simple and 
\begin{align}\label{eq-2-5}
\mu_1(a,\lambda)<\mu_2(a,\lambda)
\end{align}
which implies {\bf H3}.
Therefore, $A(a,\lambda,x)$ satisfies the assumption {\bf H1}, {\bf H2}, {\bf H3} and {\bf H4} in Section 1 of \cite{P W}, so we can define Evans' function $D(a,\lambda)$ for $(a,\lambda) \in \tilde{\Omega}$ by Definition 1.8 in \cite{P W}.
For $(a_0,\lambda_0) \in \tilde{\Omega}$ with $\re \lambda_0 >0$, from Proposition 1.9 in \cite{P W} the kernel of the operator $\partial_x(\mathcal{L}_c+a_0)-\lambda_0$  is non-trivial if and only if $D(a_0,\lambda_0)=0$.
Since $A(a,\lambda,x)$ is analytic in $a$ and $\lambda$ for each fixed $x$, Evans' function $D(a,\lambda)$ is also analytic in $a$ and $\lambda$ for $(a,\lambda) \in \tilde{\Omega}$.

Let
\[\mathscr{P}(\nu)=\nu^3-(c+a)\nu+\lambda\]
denote the characteristic polynomial of $A_{\infty}$ and
\[\tilde{\mathscr{P}}(\nu)=\nu^3+\lambda, \quad \mathscr{Q}(\nu)=-(c+a)\nu.\]
Then, the roots $\nu_0$ of $\tilde{\mathscr{P}}$ are the cube roots of $-\lambda$, and for $|\nu-\nu_0|=o(1)$ as $|\lambda| \to \infty$ we have 
\[\mathscr{Q}(\nu)=-(c+a)\nu_0(1+o(1)), \quad \frac{\partial \tilde{\mathscr{P}}}{\partial \nu}(\nu)=3\nu_0^2(1+o(1)), \quad \biggl| \frac{\mathscr{Q}(\nu_0)}{\frac{\partial \tilde{\mathscr{P}}}{\partial \nu}(\nu_0)}  \biggr|=\frac{|c+a|}{2|\lambda|^{\frac{1}{3}}}.\]
We choose $\rho(\lambda)=\rho_0|c+a|/3|\lambda|^{\frac{1}{3}}$ for any $\rho_0>1$
Then, the assumption of Lemma 1.20 in \cite{P W} are satisfied and the roots of $\mathscr{P}(\nu)=0$ are given by 
\begin{align}\label{eq-2-3-1}
\nu=(-\lambda)^{\frac{1}{3}}+O(|c+a||\lambda|^{-\frac{1}{3}})
\end{align}
as $\lambda \to \infty$.
From \eqref{eq-2-3-1} for any labeling $\nu_1(a,\lambda),\nu_2(a,\lambda),\nu_3(a,\lambda)$ of roots of $A_{\infty}(a,\lambda)$ we have 
\[ \biggl| \frac{\nu_k}{\frac{\partial \mathscr{P}}{\partial \lambda} (\nu_j)}\biggr| = \frac{|\lambda|^{\frac{1}{3}}}{3|\lambda|^{\frac{2}{3}}}(1+o(|c+a|))=O((1+|a|)|\lambda|^{-\frac{1}{3}}),\]
as $|\lambda| \to \infty$ in $\tilde{\Omega}$.
To apply Corollary 1.19 in \cite{P W}, we obtain that the hypotheses of Proposition 1.17 in \cite{P W} hold.
By Corollary 1.18 in \cite{P W}, it follows that $D(a,\lambda) \to 1$ as $|\lambda| \to \infty$ in $\tilde{\Omega}$ for each fixed $a$.
So for $0\leq a \leq\lambda_c$,
 \begin{align}\label{eq-behavior-D}
D(a,\lambda) \to 1 \mbox{ as } \lambda \to \infty.
\end{align}

Since 
\[\partial_x \mathcal{L}_c \partial_x Q_c=0, \quad\mathcal{L}_c \partial_x Q_c =0 \]
and
\[\partial_xQ_c(x) e^{\sqrt{c}x} \to -6c^{\frac{3}{2}} \mbox{ as } x \to \infty , \quad  Q_c(x) e^{-\sqrt{c}x} \to 6c \mbox{ as } x \to -\infty,\]
from (1.35) in \cite{P W} and $D(0,0)=0$ we have
\begin{align}
\frac{\partial D}{\partial a}(0,0)=&\frac{1}{\frac{\partial \mathscr{P}(\mu)}{\partial \mu}|_{(\mu,a,\lambda) =(-\sqrt{c},0,0)}}\int_{-\infty}^{\infty} \frac{Q_c}{6c} \partial_a[-\partial_x(\mathcal{L}_c+a)+\lambda]|_{(a,\lambda)=(0,0)}\frac{\partial_x Q_c}{-6c^{\frac{3}{2}}} dx \notag \\
=& \frac{-1}{72c^{\frac{7}{2}}} \int_{-\infty}^{\infty} |\partial_x Q_c|^2dx<0. \label{eq-2-6}
\end{align}
From Theorem 3.4 in \cite{C G N T} we have that the kernel of $\mathcal{L}_c+a$ on $L^2(\R)$ is trivial for $0<a<\lambda_c$.
If there exists $ 0< a_0< \lambda_c$ satisfying $D(a_0,0)=0$, then there exists a solution $u_0$ of $\partial_x(\mathcal{L}_c+a_0)u=0$ such that for all $\varepsilon >$ there is $C_{\varepsilon }>0 $ satisfying that
\[ |u_0(x)|+|\partial_x u_0(x)|+|\partial_x^2 u_0(x)| \leq C_{\varepsilon} e^{-(\mu_1-\varepsilon)x} \mbox{ as } x  \to \infty\]
and
\[ | u_0(x)| \leq C_{\varepsilon}e^{- \varepsilon x} \mbox{ as } x \to -\infty.\]
Since $((\mathcal{L}_c+a)u)(x) \to 0$ as $x \to \infty$, $u_0$ is a solution $(\mathcal{L}_c+a)u=0$.
By the property of solutions of ordinary differential equations, any solution of $(\mathcal{L}_c+a)u=0$ decays or grows exponentially tend to $-\infty$.
Thus, there are no solutions of $(\mathcal{L}_c+a)u=0$ which grows subexponentially tend to $-\infty$ and decays exponentially tend to $\infty$.
This contradicts that $u_0$ is a solution of $(\mathcal{L}_c+a)u=0$.
Thus, $D(a,0)\neq 0$ for $0<a<\lambda_c$.
Since $D(a,\lambda)$ is real and continuous for real numbers $a$ and $\lambda$ in $J_+$, by \eqref{eq-2-6} $D(a,0)$ is negative for $0<a<\lambda_c$.
From \eqref{eq-behavior-D}, for $a$ there exists $\lambda(a)>0$ such that $D(a,\lambda(a))=0.$
Therefore, $\partial_x(\mathcal{L}_c+a)$ has a positive eigenvalue $\lambda(a)$ for $0<a<\lambda(a)$.
Thus, $\partial_x\mathbb{L}_c$ has a positive eigenvalue for $L>\sqrt{\lambda_c}$.


\qed

To prove the estimate of the propagator $e^{\partial_x(\mathcal{L}+a)t}$ we apply the following  Gearhart--Greiner--Herbst--Pr\"uss theorem, see  \cite{RN}.

\begin{theorem}\label{thm-sm}
Let $\mathcal{A}$ be a generator of a strongly continuous semigroup on a complex Hilbert space $(\mathcal{H},\norm{\cdot}_{\mathcal{H}})$.
Then for each $t>0$, the following spectral mapping theorem is valid
\begin{align*}
\sigma(e^{\mathcal{A}})\setminus \{0\} =\{e^{\lambda};& \mbox{ either $\mu_k:=\lambda+2\pi i k \in \sigma (\mathcal{A}) $ for some $ k \in\Z$} \\
& \mbox{ or the sequence } \{\norm{(\mu_k-\mathcal{A})^{-1}}_{\mathcal{H} \to \mathcal{H}}\}_{k \in \Z} \mbox{ is unbounded}\}
\end{align*}

\end{theorem}

\begin{proposition}\label{prop-est-semi}
Let $a,s\geq 0$ and $s \in \Z$.
Then, for $\varepsilon >0$ there exists $C=C(\varepsilon ,s)>0$ for $u \in H^s(\R)$ and $t>0$, 
\begin{align}\label{est-semi}
\norm{e^{\partial_x(\mathcal{L}_c+a)t}u}_{H^s(\R)}\leq Ce^{(\mu(a)+\varepsilon)t}\norm{u}_{H^s(\R)} 
\end{align}
where $ \mu(a)$ is the maximum of the real pert of elements in $\sigma(\partial_x(\mathcal{L}_c+a))$.
\end{proposition}

\proof

By the compact perturbation theory the essential spectrum of $\partial_x(\mathcal{L}_c+a)$ is the essential spectrum of $\partial_x^3$, so the essential spectrum of $\partial_x(\mathcal{L}_c+a)$ is the imaginary axis.
If we show the sequence $\{\norm{(\lambda+2\pi i k - \partial_x(\mathcal{L}_c+a))^{-1}}_{H^s \to H^s}\}_k$ is bounded for all $\re \lambda >\mu(a)$, we can show the estimate \eqref{est-semi} by applying Theorem \ref{thm-sm} and Lemma 2 and 3 in \cite{S S}(see also the proof of Lemma 3.2 in \cite{YY1}).
If $s\geq 1$, we have that  for $u \in H^s(\R)$
\begin{align*}
&\norm{(\lambda+2\pi i k - \partial_x(\mathcal{L}_c+a))^{-1}u}_{H^s}\\
\lesssim &   \norm{(\lambda+2\pi i k - \partial_x(\mathcal{L}_c+a))^{-1}\partial_x u}_{H^{s-1}}+\norm{(\lambda+2\pi i k - \partial_x(\mathcal{L}_c+a))^{-1}u}_{H^{s-1}}.
\end{align*}
Here, we use the boundedness of $ (-2(Q_c)_{xx}-2(Q_c)_x\partial_x)(\lambda+2\pi i k - \partial_x(\mathcal{L}_c+a))^{-1}$ on $H^{s-1}(\R)$.
Therefore, the boundedness of the sequence $\{\norm{(\lambda+2\pi i k - \partial_x(\mathcal{L}_c+a))^{-1}}_{H^s \to H^s}\}_k$ follows the boundedness of the sequence $\{\norm{(\lambda+2\pi i k - \partial_x(\mathcal{L}_c+a))^{-1}}_{L^2 \to L^2}\}_k$.
Thus, we prove the boundedness of the sequence $\{\norm{(\lambda+2\pi i k - \partial_x(\mathcal{L}_c+a))^{-1}}_{L^2 \to L^2}\}_k$.
For $\beta \in \C$ we have 
\[(i\beta - (i\partial_x)(\mathcal{L}_c+a))^{-1}=(I+A_{\beta}B)^{-1} (i\beta - (i\partial_x) ((i\partial_x)^2+c+a))^{-1},\]
where 
\begin{align*}
A_{\beta}=&2(i\partial_x)\{i\beta-(i\partial_x)((i\partial_x)^2+c+a)\}^{-1}\sqrt{Q_c}\\
B=&\sqrt{Q_c}.
\end{align*}
Since 
\[(I+A_{\beta}B)^{-1}=I-A_{\beta}(I+BA_{\beta})^{-1}B,\]
for $\re \lambda >\mu(a)$ the sequence $\{\norm{(\lambda+2\pi i k - \partial_x(\mathcal{L}_c+a))^{-1}}_{L^2 \to L^2}\}_k$ is bounded if and only if $\{\norm{(I+BA_{\lambda+2\pi ik})^{-1}}_{L^2 \to L^2}\}_k$ is bounded.
For $u \in L^2(\R)$ we have
\begin{align*}
\norm{BA_{\lambda+2\pi ik}u}_{L^2}=&\norm{\sqrt{Q_c}2(i\partial_x)(i\lambda- 2\pi k  - (i\partial_x)((i\partial_x)^2+c+a))^{-1}(\sqrt{Q_c}u)}_{L^2}\\
\lesssim &  \norm{\eta (i\lambda- 2\pi k + \eta (\eta^2+c+a))^{-1}}_{L^1}\norm{u}_{L^2}.
\end{align*}
Let 
\[p(\eta,k)=-\im \lambda-2\pi k + \eta (\eta^2+c+a).\]
From \eqref{eq-2-3-1} in the proof of Proposition \ref{prop-properties-LO}, for $k \in \Z$ there exist roots $\alpha_j(k)$ $(j=1,2,3)$ of $p(\eta,k)=0$ satisfies
\[ \alpha_j(k) = (2\pi k)^{\frac{1}{3}}\omega_3^j+O(|k|^{-\frac{1}{3}})\]
as $|k| \to \infty$, where $\omega_3$ is a primitive root of $\eta^3-1=0$.
Since $|\im \alpha_j|=\sqrt{3}(2\pi k)^{\frac{1}{3}}+O(|k|^{-\frac{1}{3}})$ as $|k| \to \infty$, we have 
\begin{align*}
&\int_{-\infty}^{\infty} \frac{|\eta|}{|i\lambda-2\pi k + \eta(\eta^2+c+a)|} d\eta\\
 \leq
& \frac{1}{\sqrt{2}}\int_{-\infty}^{\infty} \frac{|\eta|}{|\re \lambda| + |\eta-\alpha_1(k)| |\eta-\alpha_2(k)| |\eta-\alpha_3(k)|} d\eta\\
\leq &\frac{\sqrt{2}|k|^{-\frac{1}{3}}}{|\re \lambda|}\sup_{-|k|^{-1}<\xi <|k|^{-1}} |\xi+\alpha_3(k)k^{-\frac{1}{3}}|\\
&+\frac{|k|^{-\frac{1}{3}}}{\sqrt{2}}\int_{(-1, -|k|^{-1}) \cup (|k|^{-1}, 1)} |\xi|^{-1}d\xi \sup_{\xi \in \R} \frac{|\xi+\alpha_3(k)k^{-\frac{1}{3}}|}{|\xi-(\alpha_1(k)-\alpha_3(k))k^{-\frac{1}{3}}| |\xi-(\alpha_2(k)-\alpha_3(k))k^{-\frac{1}{3}}|} \\
&+\frac{|k|^{-\frac{1}{3}}}{\sqrt{2}}\int_{(-\infty, -1) \cup (1, \infty)} \frac{|\xi+\alpha_3(k)k^{-\frac{1}{3}}|}{ |\xi-(\alpha_1(k)-\alpha_3(k))k^{-\frac{1}{3}}| |\xi-(\alpha_2(k)-\alpha_3(k))k^{-\frac{1}{3}}| } d\xi\\
\lesssim &  |k|^{-\frac{1}{3}}\log |k|.
\end{align*}
Hence, we obtain there exists $C>0$ such that
\[\norm{BA_{\lambda+2\pi ik}u}_{L^2}\leq C |k|^{-\frac{1}{3}}(\log |k|)\norm{u}_{L^2}.\]
Since $\norm{BA_{\lambda+2\pi ik}}_{L^2 \to L^2} \to 0$ as $|k| \to \infty$, $\{\norm{(I+BA_{\lambda+2\pi ik})^{-1}}_{L^2 \to L^2}\}_k$ is bounded.
Thus, we obtain the conclusion.

\qed

\section{Orbital instability}
In this section, we prove (ii) of Theorem \ref{orbital-stability} by applying the argument in \cite{R T 1}.
We assume $L>2/\sqrt{5c}$.
Let $\mu_{\max} $ be the largest eigenvalue of $\partial_x\mathbb{L}_c$.
Then, there exists a positive integer $k_0$ such that the largest eigenvalue of $\partial_x(\mathcal{L}_c+k_0^2/L^2)$ is $\mu_{\max}$.
Let $\chi $ be a eigenfunction of $\partial_x(\mathcal{L}_c+k_0^2/L^2)$ corresponding to $\mu_{\max}$.
Since $\mu_{\max}>0$, from the dichotomy for ordinary differential equations $\chi \in H^s(\R)$ for $s>0$.
For $\delta>0$ we define the solution $u^{\delta}$ of \eqref{ZKeq} with initial data $\delta \chi \cos \frac{k_0y}{L}+\tilde{Q}_c$ and we set $v^{\delta}(t,x,y)=u^{\delta}(t,x+ct,y)-\tilde{Q}_c(x)$.
Then, we have $v^{\delta}(0,x,y)=\delta \chi(x)\cos \frac{k_0y}{L}$ and
\[\partial_t v^\delta + \partial_x\mathbb{L}_cv^{\delta}+\partial_x(v^{\delta})^2=0.\]
We define $V_K^s$ as the function space
\[V_K^s=\Bigl\{u \in L^2(\RTL); u(x,y)=\sum_{j=-K}^K u_j(x)e^{\frac{ijk_0y}{L}}, u_j \in H^s(\R)   \Bigr\},\]
and we define a norm of $V_K^s$ as
\[\norm{u}_{V_K^s}=\sup_{|j|\leq K} \norm{u_j}_{H^s(\R)}, \mbox{ for } u=\sum_{j=-K}^Ku_je^{\frac{ijk_0y}{L}} \in V_K^s.\]
To show the smallness of the high frequency part of $v^{\delta}$, we consider an approximate solution 
\[v_{M}^{\delta}=\sum_{l=1}^{M}\delta^{l}w_l, \quad w_l \in V^{s-l+1}_{l}\]
where $w_1$ is the solution of
\[\partial_t w +\partial_x\mathbb{L}_cw=0, \quad w(0,x,y)=\chi(x) \cos \frac{k_0y}{L},\]
and $w_l$ is the solution of
\[ \partial_t w +\partial_x\mathbb{L}_cw + \partial_x\Bigl( \sum_{\substack{l_1,l_2 \geq 1,\\ l_1+l_2=l}}w_{l_1}w_{l_2}\Bigr)=0, \quad w(0,x,y) =0.\]
Then, $v_M^{\delta}$ satisfies 
\[ \partial_t v_M^{\delta} + \partial_x\mathbb{L}_cv_M^{\delta} + \partial_x(v_M^{\delta})^2=F,\]
where
\[F=\delta^M\partial_x\Bigl(\sum_{\substack{1\leq l_1, l_2 \leq M,\\ l_1+l_2>M} }\delta^{l_1+l_2-M} w_{l_1}w_{l_2} \Bigr).\]

From Proposition \ref{prop-est-semi}, we have the following lemma.
\begin{lemma}\label{semi-est}
For $K,s,\varepsilon >0$ there exists $C_{K,s,\varepsilon }>0$ such that for $u \in V_K^s$ 
\[ \norm{e^{t\partial_x\mathbb{L}_c}u}_{V_K^s} \leq C_{K,s,\varepsilon } e^{(\mu_{\max}+\varepsilon )t} \norm{u}_{V_K^s}.\]
\end{lemma}

Let $w^{\delta}=v^{\delta}-v_M^{\delta}$.
Then, we have
\begin{align}\label{3-e-0}
\partial_t w^{\delta}+\partial_x \mathbb{L}_cw^{\delta} +2\partial_x(w^{\delta}v_M^{\delta})+\partial_x(w^{\delta}w^{\delta})+F=0.
\end{align}
Therefore,
\begin{align}
\frac{d\norm{w^{\delta}}_{L^2}^2}{d t}=&  \int_{\RTL}\bigl((w^{\delta})^2\partial_x\tilde{Q}_c - (w^{\delta})^2\partial_x v_M^{\delta} -Fw^{\delta}\bigr)dxdy \notag\\
\leq &  \Bigl(1+\norm{\partial_xv_M^{\delta}}_{L^{\infty}}+\norm{\partial_x \tilde{Q}_c}_{L^{\infty}}\Bigr)\norm{w^{\delta}}_{L^2}^2 +\norm{F}_{L^2}^2 \label{3-e-1}
\end{align}
From Lemma \ref{semi-est} we have that for $\varepsilon _0>0$  there exists $C_{M,s,\varepsilon _0}>0$ such that
\[\norm{w_l(t)}_{H^s}\leq C_{M,s,\varepsilon _0} e^{l(\mu_{\max}+\varepsilon_0 )t}.\]
Therefore,  there exists $C_{M,\varepsilon _0}>0$ such that we have
\begin{align}
\norm{\partial_xv_M^{\delta}(t)}_{L^{\infty}} \leq & C_{M,\varepsilon _0} (\delta e^{(\mu_{\max}+\varepsilon_0)t}+\delta^M e^{M(\mu_{\max}+\varepsilon_0)t})\label{3-e-2}\\
\norm{F}_{L^2}\leq & C_{M,\varepsilon _0}\delta^{M+1}e^{(M+1)(\mu_{\max}+\varepsilon_0)t}.\label{3-e-3}
\end{align}
We set $T_{\delta,\varepsilon }=(\log(\varepsilon )-\log(\delta))/2\mu_{\max}$.
Since $e^{(\mu_{\max}+\varepsilon_0)t}\leq \varepsilon/\delta$ for $0<t\leq T_{\delta,\varepsilon}$, by \eqref{3-e-1}--\eqref{3-e-3} we have
\[\frac{d\norm{w^{\delta}(t)}_{L^2}^2}{d t}\leq \Bigl(1+\norm{\partial_x\tilde{Q}_c}_{L^\infty}+2\varepsilon C_{M,\varepsilon _0}\Bigr)\norm{w^{\delta}(t)}_{L^2}^2 + C_{M,\varepsilon _0}^2 \delta^{2(M+1)}e^{2(M+1)(\mu_{\max}+\varepsilon_0)t}\]
for any $0<\varepsilon <1$ and $0<t\leq T_{\delta,\varepsilon}$.
Thus,
\[\frac{d}{d t}\Bigl(e^{-(1+\norm{\partial_x\tilde{Q}_c}_{L^\infty}+2\varepsilon C_{M,\varepsilon _0})t}\norm{w^{\delta}(t)}_{L^2}^2\Bigr)\leq C_{M,\varepsilon _0}^2 \delta^{2(M+1)}e^{2(M+1)(\mu_{\max}+\varepsilon_0)t-(1+\norm{\partial_x\tilde{Q}_c}_{L^\infty}+2\varepsilon C_{M,\varepsilon _0})t}.\]
If we choose large $M$ and small $\varepsilon (M)$ satisfying
\[2(M+1)(\mu_{\max}+\varepsilon_0)-(1+\norm{\partial_xQ_x}_{L^\infty}+2\varepsilon C_{M,\varepsilon _0})>0,\]
then we obtain 
\[\norm{w^{\delta}(t)}_{L^2}^2\leq C'_{M,\varepsilon _0} \delta^{2(M+1)}e^{2(M+1)(\mu_{\max}+\varepsilon_0)t}\]
for $0<t\leq T_{\delta,\varepsilon}$.
Hence, there exists $C''_{M,\varepsilon _0} >0$ such that 
\[\norm{w^{\delta}(T_{\delta,\varepsilon})}_{L^2}\leq C''_{M,\varepsilon _0}\varepsilon ^{M+1}\]
for small $\varepsilon >0$.
Let $P_0$ be a projection satisfying
\[(P_{0} u)(x,y)=\int_{\T_L} u(x,z)dz, \mbox{ for } (x,y) \in \RTL.\]
From the definition of $v_M^{\delta}$ and the estimate \eqref{3-e-0} we have 
\begin{align*}\label{3-e-4}
\norm{(Id-P_0)v_M^{\delta}(t)}_{L^2} \geq \sqrt{\pi}\norm{\chi}_{L^2}\delta e^{\mu_{\max}t}-C_{\varepsilon _0} (\delta^2 e^{2\mu_{\max}t}+\delta^M e^{M\mu_{\max}t})
\end{align*}

\begin{align*}
\inf_{a \in \R} \norm{u^{\delta}(T_{\delta,\varepsilon},\cdot,\cdot)-\tilde{Q}_c(\cdot +a,\cdot)}_{L^2} \geq & \norm{(Id-P_0)(u^{\delta}(T_{\delta,\varepsilon})-R_c(T_{\delta,\varepsilon}))}_{L^2}\\
=& \norm{(Id-P_0)(v^{\delta}(T_{\delta,\varepsilon}))}_{L^2}\\
\geq & \norm{(Id-P_0)v_M^{\delta}(T_{\delta,\varepsilon})}_{L^2}-\norm{w^{\delta}(T_{\delta,\varepsilon})}_{L^2}\\
\geq & \sqrt{\pi}\norm{\chi}_{L^2} \varepsilon - C_{M,\varepsilon _0}'''\varepsilon ^{2}.
\end{align*}
Thus, if we choose 
\[\varepsilon_1 = \frac{\sqrt{\pi}\norm{\chi}_{L^2}}{2 C_{M,\varepsilon _0}'''},\]
for any $\delta>0$ there exists $T_{\delta,\varepsilon}>0$ such that
\[\inf_{a \in \R} \norm{u^{\delta}(T_{\delta,\varepsilon},\cdot,\cdot)-\tilde{Q}_c(\cdot +a)}_{L^2}\geq \frac{\sqrt{\pi}\norm{\chi}_{L^2} \varepsilon_1}{2}.\]
This completes the proof of (ii) in Theorem \ref{orbital-stability}.

\section{Orbital stability}
In this section, we prove (i) of Theorem \ref{orbital-stability} by applying the arguments in \cite{G S S 1} and \cite{YY2}.
We write the outline of the proof of (i) of Theorem \ref{orbital-stability}.

Theorem 3.3 in \cite{G S S 1} yields the following coercive type lemma for $\mathcal{L}_{c_0}$.
\begin{lemma}\label{lem-coer0}
Let $c_0>0$.
There exists $k_0>0$ such that for $u\in H^1(\R)$ with $(u, Q_{c_0})_{L^2(\R)}=(u,\partial_x Q_{c_0})_{L^2(\R)}=0$, 
\[\tbr{\mathcal{L}_{c_0}u,u}_{H^{-1}(\R),H^1(\R)} \geq k_0 \norm{u}_{H^1(\R)}^2.\]
\end{lemma}

\subsection{Non-critical case $L < \frac{2}{\sqrt{5c_0}}$}
To show the orbital stability of $\tilde{R}_{c_0}$ for $L<\frac{2}{\sqrt{5c_0}}$, we apply the argument in \cite{MIW1} (see also \cite{YY1, YY3}).
Let $L< \frac{2}{\sqrt{5c_0}}$.
By the Fourier expansion \eqref{eq-2-1} we have for $u \in H^1(\RTL)$ 
\[\tbr{\mathbb{L}_{c_0}u,u}_{H^{-1}(\RTL),H^1(\RTL)} = \sum_{n=-\infty}^{\infty} \Tbr{\Bigl( \mathcal{L}_{c_0} +\frac{n^2}{L^2} \Bigr)u_n,u_n}_{H^{-1}(\R),H^1(\R)},\]
where
\[u(x,y)=\sum_{n=-\infty}^{\infty}u_n(x)e^{\frac{iny}{L}}.\]
Since $\lambda_{c_0}<L^{-2}$, $\mathcal{L}_{c_0}+ n^2/L^2$ is positive for $|n| \geq 1$.
From Lemma \ref{lem-coer0} there exists $K_0>0$ such that for $u \in H^1(\RTL)$ with $(u, \tilde{Q}_{c_0})_{L^2(\RTL)}=(u,\partial_x\tilde{Q}_{c_0})_{L^2(\RTL)}=0$, we have
\begin{align}\label{eq-5-1}
\tbr{\mathbb{L}_{c_0}u,u}_{H^{-1}(\RTL),H^1(\RTL)} \geq K_0 \norm{u}_{H^1(\RTL)}^2.
\end{align}
Combing \eqref{eq-5-1} and the proofs of Theorem 3.4 and Theorem 3.5 in \cite{G S S 1}, we obtain the orbital stability of $\tilde{R}_{c_0}$.

\subsection{Critical case $L= \frac{2}{\sqrt{5c_0}}$}
The proof of the orbital stability of $\tilde{R}_{c_0}$ for $L=\frac{2}{\sqrt{5c_0}}$ is similar to the proof of (i) of Theorem 1.4 in \cite{YY2} (see also the proof of (i) of Theorem 1.4 \cite{YY3}).
Let $L=\frac{2}{\sqrt{5c_0}}$.
In this case, from (iii) of Proposition \ref{prop-properties-LO} the linearized operator $\mathbb{L}_{c_0}$ has an extra eigenfunction corresponding to the zero eigenvalue.
Therefore, we have to recover the degeneracy of the kernel of $\mathbb{L}_{c_0}$ from nonlinearity of \eqref{ZKeq}.
We define the action $S_c(u)$ by $E(u)+cM(u)$.


\begin{lemma}\label{lem-p-3-1}
There exist a neighborhood $U$ of $(0,0)$ and a $C^2$ function $ \gamma_c(\vec{a}): U  \to \R $
such that $\gamma_c(0,0)=c$ and for $\vec{a} \in U$ and $|c-c_0|<c_0/2$  
\[ M\bigl(\Theta(\vec{a},\gamma_{c}(\vec{a}))\bigr)=M(\tilde{Q}_{c}),\]
\begin{align}\label{eq-5-2}
\gamma_c(\vec{a})-c =-\frac{cC_{2,c_0}}{3\norm{\tilde{Q}_{c_0}}_{L^2}^2} |\vec{a}|^2+o(|\vec{a}|^2),
\end{align}
where $\Theta(\vec{a},c)(x,y)=cc_0^{-1}\varphi_{c_0}(\vec{a})(\sqrt{cc_0^{-1}}x,y)$.
\end{lemma}
\proof
Let
\[\gamma_c(\vec{a})=c_0\Bigl( \norm{\tilde{Q}_{c}}_{L^2}^2\norm{\varphi_{c_0}(\vec{a})}_{L^2}^{-2} \Bigr)^{\frac{2}{3}}.\]
By the definition of $\Theta$ we have 
\[M\bigl(\Theta(\vec{a},\gamma_c(\vec{a}))\bigr)=M(\tilde{Q}_{c}).\]
Since $\norm{\tilde{Q}_{c}}_{L^2}^2\norm{\tilde{Q}_{c_0}}_{L^2}^{-2}=c^{\frac{3}{2}}c_0^{-\frac{3}{2}}$, we have
\begin{align*}
\gamma_c(\vec{a})=&c-c\frac{\norm{\varphi_{c_0}(\vec{a})}_{L^2}^{\frac{4}{3}}-\norm{\tilde{Q}_{c_0}}_{L^2}^{\frac{4}{3}}}{\norm{\varphi_{c_0}(\vec{a})}_{L^2}^{\frac{4}{3}}}
= c-\frac{cC_{2,c_0}}{3\norm{\tilde{Q}_{c_0}}_{L^2}^2} |\vec{a}|^2 +o(|\vec{a}|^2).
\end{align*}


\qed

Next, we investigate the difference between $\Theta$ and $\tilde{Q}_c$ on the action $S_c$.
\begin{lemma}\label{lem-p-3-2}
For $\vec{a} \in U$ and $|c-c_0|<c_0/2$,
\begin{align}\label{eq-a-term}
S_{c}\bigl(\Theta(\vec{a},\gamma_c(\vec{a}))\bigr)-S_{c}(\tilde{Q}_{c})=&\Bigl(\frac{c}{c_0}\Bigr)^{\frac{5}{2}}\frac{5c_0C_{2,c_0}\norm{\tilde{Q}_{c_0}^{\frac{3}{2}}\cos \frac{y}{L}}_{L^2}^2}{48\norm{\tilde{Q}_{c_0}}_{L^2}^2}|\vec{a}|^4 \notag \\
&+ \Bigl(1-\frac{c}{c_0}\Bigr)\norm{\partial_y \Theta(\vec{a},\gamma_c(\vec{a}))}_{L^2}^2 +o(|\vec{a}|^4) 
\end{align}
as $|\vec{a}| \to 0$.
\end{lemma}
\proof
First, we consider the case $c=c_0$.
From the expansion
\begin{align}\label{eq-5-exp}
 \Theta(\vec{a}, \gamma_{c_0}(\vec{a}))=\varphi_{c_0}(\vec{a})(x,y) + (\gamma_{c_0}(\vec{a})-c_0) \partial_c Q_{c_0} +O\bigl((|\vec{a}|+(\gamma_{c_0}(\vec{a})-c_0))(\gamma_{c_0}(\vec{a})-c_0)\bigr),
 \end{align}
we have
\begin{align*}
S_{c_0}\bigl(\Theta(\vec{a},\gamma_{c_0}(\vec{a})\bigr)-S_{c_0}(\tilde{Q}_{c_0})
=& S_{\check{c} (\vec{a})}(\varphi_{c_0}(\vec{a}))-S_{c_0}(\tilde{Q}_{c_0}) + (c_0-\check{c} (\vec{a}))M(\tilde{Q}_{c_0})\\
&+\frac{1}{2}(\gamma_{c_0}(\vec{a})-c_0)^2(S_{c}''(\tilde{Q}_{c_-})\partial_c\tilde{Q}_{c_0},\partial_c\tilde{Q}_{c_0})_{L^2} + o(|\vec{a}|^4),
\end{align*}
where $\check{c}$ is defined in Proposition \ref{prop-bifurcation}.
Since $\frac{\partial \check{c}}{\partial a_1} (0,0)=0$ and $\frac{\partial^2 \check{c}}{\partial (a_1)^2} (0,0)=\check{c} ''(0)>0$, there exist $\delta_1>0$ and the inverse function $a_1(c)$ of $\check{c} (a_1,0)$ on from $[ c_0, \check{c} (\delta_1,0))$ to $[0, \delta_1)$.
For $c_1,c_2$ with $c_1\neq c_2$ 
\begin{align*}
&\frac{S_{c_1}(\varphi_{c_0}(c_1))-S_{c_2}(\varphi_{c_0}(c_2))}{c_1-c_2}\\
=&\frac{\bigl(S_{c_2}''(\varphi_{c_0}(c_2))(\varphi_{c_0}(c_1)-\varphi_{c_0}(c_2)),\varphi_{c_0}(c_1)-\varphi_{c_0}(c_2)\bigr)_{L^2}}{2(c_1-c_2)}+M(\varphi_{c_0}(c_1))\\
&+\frac{o((\varphi_{c_0}(c_1)-\varphi_{c_0}(c_2))^2)}{c_1-c_2}\\
\to & M(\varphi_{c_0}(c_2)) \mbox{ as } c_1 \to c_2,
\end{align*}
where $\varphi_{c_0}(c)=\varphi_{c_0}(a_1(c),0)$.
Since $S_{c_0}''(\tilde{Q}_{c_0})\partial_{a_1}\varphi_{c_0}(a_1,a_2)|_{(a_1,a_2)=(0,0)}=\mathbb{L}_{c_0}(\tilde{Q}_{c_0}^{\frac{3}{2}}\cos \frac{y}{L})=0$, for $c>c_0$ 
\begin{align*}
&\frac{S_{c}(\varphi_{c_0}(c))-S_{c_0}(\tilde{Q}_{c_0})}{c-c_0}\\
=&\frac{\bigl(S_{c_0}''(\tilde{Q}_{c_0})(\varphi_{c_0}(c)-\tilde{Q}_{c_0}),\varphi_{c_0}(c)-\tilde{Q}_{c_0}\bigr)_{L^2}}{\check{c} ''(0)a_1(c)^2+o(a_1(c)^2)}+M(\varphi_{c_0}(c))+\frac{o((\varphi_{c_0}(c)-\tilde{Q}_{c_0})^2)}{\check{c} ''(0)a_1(c)^2+o(a_1(c)^2)}\\
\to & M(\tilde{Q}_{c_0}) \mbox{ as } c \downarrow  c_0.
\end{align*}
Therefore, $S_c(\varphi_{c_0}(c))$ is $C^1$ and $\partial_c S_c(\varphi_{c_0}(c))=M(\varphi_{c_0}(c))$.
By the same way we obtain that $M(\varphi_{c_0}(c))$ is $C^1$ and
\[ \lim_{c \downarrow c_0} \frac{M(\varphi_{c_0}(c))-M(\tilde{Q}_{c_0})}{c-c_0}=\frac{C_{2,c_0}}{2\check{c} ''(0)}.\]
Thus, we have
\begin{align}
&S_{\check{c} (|\vec{a}|,0)}(\varphi_{c_0}(|\vec{a}|,0))-S_{c_0}(\tilde{Q}_{c_0}) + (c_0-\check{c} (|\vec{a}|,0))M(\tilde{Q}_{c_0})\notag \\
=& \frac{C_{2,c_0}}{4\check{c} ''(0)}(\check{c} (|\vec{a}|,0)-c_0)^2+o((\check{c} (|\vec{a}|,0)-c_0)^2) \notag \\
=&\frac{C_{2,c_0}\check{c} ''(0)}{16}|\vec{a}|^4+o(|\vec{a}|^4).\label{eq-5-3}
\end{align}
From Lemma \ref{lem-p-3-1} and $S_{c_0}''(\tilde{Q}_{c_0})\partial_c\tilde{Q}_{c_0}=-\tilde{Q}_{c_0}$,
\begin{align}\label{eq-5-4}
(\gamma_{c_0}(\vec{a})-c_0)^2\bigl(S_{c_0}''(\tilde{Q}_{c_0})\partial_c\tilde{Q}_{c_0},\partial_c\tilde{Q}_{c_0}\bigr)_{L^2}= -\frac{c_0C_{2,c_0}^2}{12\norm{\tilde{Q}_{c_0}}_{L^2}^2}|\vec{a}|^4+o(|\vec{a}|^4).
\end{align}
Since 
\begin{align*}
&S_{\check{c} (\vec{a})}(\varphi_{c_0}(\vec{a}))-S_{c_0}(\tilde{Q}_{c_0}) + (c_0-\check{c} (\vec{a}))M(\tilde{Q}_{c_0})\\
=&S_{\check{c} (|\vec{a}|,0)}(\varphi_{c_0}(\vec{a}))-S_{c_0}(\tilde{Q}_{c_0}) + (c_0-\check{c} (|\vec{a}|,0))M(\tilde{Q}_{c_0}),
\end{align*}
from \eqref{eq-5-3} and \eqref{eq-5-4} we obtain \eqref{eq-a-term} for $c=c_0$.

Next, we consider the general cases.
Since $M(\Theta(\vec{a},\gamma_c(\vec{a})))=M(\tilde{Q}_c)$, we have
\begin{align*}
&S_c\bigl(\Theta(\vec{a},\gamma_c(\vec{a}))\bigr)-S_c(\tilde{Q}_c)\\
=&\Bigl(\frac{c}{c_0}\Bigr)^{\frac{5}{2}}\Bigl( S_{c_0}\bigl( \Theta(\vec{a},\gamma_{c_0}(\vec{a})) \bigr) -S_{c_0}(\tilde{Q}_{c_0})\Bigr)  + \Bigl(1-\frac{c}{c_0}\Bigr) \norm{\partial_y \Theta(\vec{a},\gamma_c(\vec{a}))}_{L^2}^2.
\end{align*}
Therefore, we obtain \eqref{eq-a-term} for $c>0$.

\qed

We define a distance $\mbox{dist}_{c}$ and neighborhoods $N_{\varepsilon ,c}$ and $N_{\varepsilon ,c}^0$ of $\tilde{Q}_{c_0}$ by 
\[\mbox{dist}_c(u) =\inf_{x \in \R} \norm{u(\cdot,\cdot)-\tilde{Q}_{c}(\cdot-x,\cdot)}_{H^1},\]
\[N_{\varepsilon ,c}=\{ u \in H^1(\RTL); \mbox{dist}_c(u)<\varepsilon \},\]
\[N_{\varepsilon ,c}^l=\{u \in N_{\varepsilon ,c}; M(u)=M(\tilde{Q}_l)\}.\]
In the following lemma, to get a orthogonal condition we decompose functions in $N_{\varepsilon ,c}$.
\begin{lemma}\label{lem-p-3-3}
Let $\varepsilon >0$ sufficiently small.
Then, there exist $K_1>0$, $C^2$ functions $\rho:N_{\varepsilon ,c_0} \to \R$, $ c:N_{\varepsilon ,c_0} \to \R$, $\vec{a}=(a_1,a_2): N_{\varepsilon ,c_0} \to U$ and $\eta: N_{\varepsilon ,c_0} \to H^1(\RTL)$ such that for $u \in N_{\varepsilon ,c_0}$ 
\[u(\cdot+\rho(u),\cdot)=\Theta(\vec{a}(u),c(u))(\cdot,\cdot)+\eta(u)(\cdot,\cdot),\]
\begin{align}\label{eq-5-5} 
|c(u)-c_0|+|\vec{a}(u)|+\norm{\eta(u)}_{H^1} \leq K_1 \mbox{\rm dist}_{c_0}(u),
\end{align}
and $\bigl(\eta(u), \Theta(\vec{a}(u),c(u))\bigr)_{L^2}=\bigl(\eta(u), \partial_x \Theta(\vec{a}(u),c(u))\bigr)_{L^2}=\bigl(\eta(u),\partial_{a_1}\Theta(\vec{a}(u),c(u))\bigr)_{L^2}=\bigl(\eta(u),\partial_{a_2}\Theta(\vec{a}(u),c(u))\bigr)_{L^2}=0$.
\end{lemma}
\proof
We define
\[G(u,c,\rho,a_1,a_2)=
\begin{pmatrix}
(u(\cdot+\rho,\cdot)-\Theta(\vec{a},c),\Theta(\vec{a},c))_{L^2}\\
(u(\cdot+\rho,\cdot)-\Theta(\vec{a},c),\partial_x \Theta(\vec{a},c))_{L^2}\\
(u(\cdot+\rho,\cdot)-\Theta(\vec{a},c),\partial_{a_1}\Theta(\vec{a},c))_{L^2}\\
(u(\cdot+\rho,\cdot)-\Theta(\vec{a},c),\partial_{a_2}\Theta(\vec{a},c))_{L^2}
\end{pmatrix}.\]
Then, $G(\tilde{Q}_{c_0},c_0,0,0,0)=0$.
Since
\[
\frac{\partial G}{\partial (c,\rho,a_1,a_2)}\Bigl|_{\substack{u=\tilde{Q}_{c_0},c=c_0 \\ \rho=a_1=a_2=0}}=
\begin{pmatrix}
-(\partial_c\tilde{Q}_{c_0},\tilde{Q}_{c_0})_{L^2} &0&0&0\\
0& \norm{\partial_x \tilde{Q}_{c_0}}_{L^2}^2& 0 &0\\
0 & 0& -\norm{\tilde{Q}_{c_0}^{\frac{3}{2}}\cos \frac{y}{L}}_{L^2}^2 &0\\
0 & 0& 0& -\norm{\tilde{Q}_{c_0}^{\frac{3}{2}}\sin \frac{y}{L}}_{L^2}^2
\end{pmatrix}
\]
is regular, 
from the implicit function theorem for small $\varepsilon >0$ there exists $C^2$ functions $c,\rho,a_1,a_2:N_{\varepsilon ,c_0} \to \R$ such that for $u \in N_{\varepsilon ,c_0}$ 
\[G(u,c(u),\rho(u),a_1(u),a_2(u))=0.\]
Therefore,
\[ \eta(u)=u(\cdot+\rho(u),\cdot)-\Theta(\vec{a}(u),c(u))\]
satisfies the orthogonal conditions, where $\vec{a}(u)=(a_1(u),a_2(u))$.
The inequality \eqref{eq-5-5} follows the implicit function theorem and the definition of $\eta$.
\qed

In the following lemma, we estimate $\norm{\Theta(\vec{a}(u),c(u))-\Theta(\vec{a}(u),\gamma(\vec{a}(u)))}_{H^1}$ on $N_{\varepsilon ,c_0}^0$.
\begin{lemma}\label{lem-p-3-4}
Let $\varepsilon >0$ sufficiently small.
There exists $C>0$ such that for $|l-c_0|<\varepsilon^{1/2}$ and $u \in N_{\varepsilon ,c_0}^l$,
\[ \norm{\Theta\bigl(\vec{a}(u),\gamma_l(\vec{a}(u))\bigr)-\Theta(\vec{a}(u),c(u))}_{H^1}\leq C \norm{\eta(u)}_{L^2}^2,\]
\begin{align}\label{eq-est-c-1}
| \gamma_l(\vec{a}(u))-c(u)| \lesssim M(\eta(u)).
\end{align}
\end{lemma}
\proof
For $u \in N_{\varepsilon ,c_0}^l$, 
\begin{align*}
M\bigl(\Theta(\vec{a}(u),\gamma_l(\vec{a}(u))\bigr)=M(\tilde{Q}_{l})=&M\bigl(\eta(u)+\Theta(\vec{a}(u),c(u))\bigr)\\
=&M(\eta(u))+M\bigl(\Theta(\vec{a}(u),c(u))\bigr).
\end{align*}
For sufficiently small $\varepsilon >0$, we have
\[ |c(u)-c_0|+|\gamma_l(\vec{a}(u))-c_0| < \frac{c_0}{2}.\]
Therefore,
\begin{align*}
M(\eta(u))=M\bigl(\Theta(\vec{a}(u),\gamma_l(\vec{a}(u)))\bigr)-M\bigl(\Theta(\vec{a}(u),c(u))\bigr)=&(\gamma_l(\vec{a}(u))^{\frac{3}{2}}-c(u)^{\frac{3}{2}})M(\varphi_{c_0}(\vec{a}(u)))\\
\gtrsim& \gamma_l(\vec{a}(u))-c(u) \geq 0.
\end{align*}
Since 
\[ \Theta\bigl(\vec{a}(u),\gamma_l(\vec{a}(u))\bigr)-\Theta(\vec{a}(u),c(u))=\bigl(\gamma_l(\vec{a}(u))-c(u)\bigr)\partial_c\tilde{Q}_{c_0}+o(\gamma_l(\vec{a}(u))-c(u)),\]
\begin{align*}
\norm{\Theta(\vec{a}(u),\gamma_l(\vec{a}(u)))-\Theta(\vec{a}(u),c(u))}_{H^1} \lesssim \gamma_l(\vec{a}(u))-c(u) \lesssim M(\eta(u)).
\end{align*}

\qed

Next we show the coerciveness of $S_{c_0}''(\tilde{Q}_{c_0})$ on a subspace of $H^1(\RTL)$.
\begin{lemma}\label{lem-p-coer}
There exist $k_2>0$ and $\varepsilon _0>0$ such that for $a_1,a_2 \in (-\varepsilon _0,\varepsilon _0)$ and $c \in (c_0-\varepsilon _0,c_0+\varepsilon _0)$, if $w \in H^1(\RTL)$ satisfies 
\[(w,\Theta(\vec{a},c))_{L^2}=(w,\partial_x\Theta(\vec{a},c))_{L^2}=(w,\partial_{a_1}\Theta(\vec{a},c))_{L^2}=(w,\partial_{a_2}\Theta(\vec{a},c))_{L^2}=0,\]
 then
\[\tbr{S_{c_0}''(\Theta(\vec{a},c))w,w}_{H^{-1}(\RTL),H^1(\RTL)} \geq k_2 \norm{w}_{H^1}^2.\]
\end{lemma}
\proof
By the definition of $S_c$, $S_{c_0}''(\tilde{Q}_{c_0})=\mathbb{L}_{c_0}$.
Since $\mathcal{L}_{c_0}+n^2L^{-2}$ is positive for $|n|\geq 2$, from Lemma \ref{lem-coer0} we obtain that there exists $k_2'>0$ such that for $u \in H^1(\RTL)$ satisfying $(u,\tilde{Q}_{c_0})_{L^2}=(u,\partial_x\tilde{Q}_{c_0})_{L^2}=(u,\tilde{Q}_{c_0}^{\frac{3}{2}}\cos \frac{y}{L})_{L^2}=(u,\tilde{Q}_{c_0}^{\frac{3}{2}}\sin \frac{y}{L})_{L^2}=0$, 
\[\tbr{S_{c_0}''(\tilde{Q}_{c_0})u,u}_{H^{-1}(\RTL),H^1(\RTL)} \geq k_2'\norm{u}_{H^1}^2.\]
By a continuity argument we obtain the conclusion.

\qed

Next, we show (i) of Theorem \ref{orbital-stability}
\proof[Proof of (i) of Theorem \ref{orbital-stability}]
Let $\varepsilon  >0$ sufficiently small.
Applying Lemma \ref{lem-p-3-1}--\ref{lem-p-coer}, we obtain that for $u\in N_{\varepsilon ,c_0}^{c_0}$
\begin{align*}
&S_{c_0}(u)-S_{c_0}(\tilde{Q}_{c_0})\\
=& S_{c_0}\bigl(\Theta(\vec{a}(u),c(u))+\eta(u)\bigr)-S_{c_0}(\tilde{Q}_{c_0})\\
=& S_{c_0}\bigl(\Theta(\vec{a}(u),\gamma_{c_0}(\vec{a}(u)))\bigr)-S_{c_0}(\tilde{Q}_{c_0})\\
&+ \tbr{S_{c_0}'\bigl(\Theta(\vec{a}(u),\gamma_{c_0}(\vec{a}(u)))\bigr), \eta(u)+\Theta(\vec{a}(u),c(u))-\Theta\bigl(\vec{a}(u),\gamma_{c_0}(\vec{a}(u))\bigr)}_{H^{-1},H^1}\\
&+\frac{1}{2}\tbr{S_{c_0}''\bigl(\Theta(\vec{a}(u),\gamma_{c_0}(\vec{a}(u)))\bigr)\eta(u),\eta(u)}_{H^{-1},H^1} + o(\norm{\eta(u)}_{H^1}^2)\\
\geq & \frac{5c_0C_{2,c_0}\norm{\tilde{Q}_{c_0}^{\frac{3}{2}}\cos \frac{y}{L}}_{L^2}^2}{48\norm{\tilde{Q}_{c_0}}_{L^2}^2}|\vec{a}(u)|^4+k_2\norm{\eta(u)}_{H^1}^2\\
&+ \tbr{S_{c_0}'\bigl(\Theta(\vec{a}(u),\gamma_{c_0}(\vec{a}(u)))\bigr), \eta(u)}_{H^{-1},H^1}+o(\norm{\eta(u)}_{H^1}^2+|\vec{a}(u)|^4).
\end{align*}
Since $S_{c_0}''(\tilde{Q}_{c_0})\partial_c \tilde{Q}_{c_0}=-\tilde{Q}_{c_0}$ and the expansion \eqref{eq-5-exp}, from Lemma \ref{lem-p-3-1} we have
\begin{align*}
&\tbr{S_{c_0}'\bigl(\Theta(\vec{a}(u),\gamma_{c_0}(\vec{a}(u)))\bigr), \eta(u)}_{H^{-1},H^1}\\
=&\tbr{S_{\check{c} (\vec{a}(u))}'\bigl(\Theta(\vec{a}(u),\gamma_{c_0}(\vec{a}(u)))\bigr),\eta(u)}_{H^{-1},H^1}\\
=&\tbr{\bigl( S_{\check{c} (\vec{a}(u))}''(\varphi_{c_0}(\vec{a}(u))-S_{c_0}''(\tilde{Q}_{c_0}) \bigr)(\gamma_{c_0}(\vec{a}(u))-c_0)\partial_c \tilde{Q}_{c_0}, \eta(u)}_{H^{-1},H^1}\\
&+ (\gamma_{c_0}(\vec{a}(u))-c_0)\bigl(\tilde{Q}_{c_0},\eta(u)\bigr)_{L^2}+o(|\vec{a}(u)|^4+\norm{\eta(u)}_{H^1}^2)\\
=&o(|\vec{a}(u)|^4+\norm{\eta(u)}_{H^1}^2).
\end{align*}
Therefore, there exist $\varepsilon _*, k_*>0$ such that for $u \in N_{\varepsilon _*,c_0}^{c_0}$ 
\begin{align}\label{eq-5-7}
S_{c_0}(u)-S_{c_0}(\tilde{Q}_{c_0})\geq k_*( |\vec{a}(u)|^4+\norm{\eta(u)}_{H^1}^2).
\end{align}

Now we suppose there exist $\varepsilon _0>0$, a sequence $\{u_n\}_n$ of solutions to \eqref{ZKeq} and a sequence $\{t_n\}$ such that $t_n>0$, $u_n(0) \to \tilde{Q}_{c_0}$ as $n \to \infty$ in $H^1$ and $\mbox{\rm dist}_{c_0}(u_n(t_n))>\varepsilon _0$.
Let $v_n=M(\tilde{Q}_{c_0})^{-\frac{1}{2}}M(u_n)^{-\frac{1}{2}}u_n(t_n)$.
Then we have $M(v_n)=M(\tilde{Q}_{c_0})$, $\lim_{n \to \infty}\norm{v_n - u_n(t_n)}_{H^1}=0$ and $\lim_{n \to \infty} S_{c_0}(v_n) = S_{c_0}(\tilde{Q}_{c_0})$.
Thus, by \eqref{eq-5-7} $\lim_{n \to \infty}\vec{a}(v_n)=0$ and $\eta(v_n) \to 0$ as $n \to \infty$ in $H^1$.
Since $\lim_{n \to \infty} \gamma_{c_0}(\vec{a}(v_n))=c_0$, we have $\lim_{n \to \infty}c(v_n) =c_0$.
Hence, $\lim_{n \to \infty} \mbox{\rm dist}_{c_0}(u_n(t_n)) =0$.
This is a contradiction.
We complete the proof of (i) of Theorem \ref{orbital-stability}.

In the following corollary, we estimates the size of the modulation parameters.
\begin{corollary}\label{cor-initial-c}
Let $c_0>0$ and $L=\frac{2}{\sqrt{5c_0}}$.
Then, there exist $\delta_0, C>0$ such that for $0<\delta<\delta_0$ and $u_0 \in H^1(\RTL)$ with $\norm{u_0-\tilde{Q}_{c_0}}_{H^1}<\delta$, the solution $u$ of $\eqref{ZKeq}$ corresponding to the initial data $u_0$ satisfies 
\[ |c(u(t))-c_0|+|\vec{a}(u(t))|^2 \leq C \delta, \quad t \in \R,\]
where $c(u)$ and $\vec{a}(u)$ are defined in Lemma \ref{lem-p-3-3}.
\end{corollary}
\proof
We choose $\varepsilon >0$ which is sufficiently small.
By (i) of Theorem \ref{orbital-stability}, there exists $\delta_1>0$ such that for any solution $u$ with $\norm{u-\tilde{Q}_{c_0}}_{H^1}=\delta<\delta_1$ satisfies $u(t) \in N_{\varepsilon ,c_0}$ for $t \in \R$.
We define $c_m>0$ as 
\[\norm{u_0}_{L^2}=\norm{\tilde{Q}_{c_m}}_{L^2}.\] 
Applying Lemma \ref{lem-p-3-1}--\ref{lem-p-3-4}, we obtain
\begin{align*}
&S_{c_m}(u)-S_{c_m}(\tilde{Q}_{c_m})\\
=&\frac{1}{2}\tbr{S_{c_0}''\bigl(\Theta(\vec{a}(u),\gamma_{c_m}(\vec{a}(u)))\bigr)\eta(u),\eta(u)}_{H^{-1},H^1}  \\
&+\Bigl( \frac{c_m}{c_0} \Bigr)^{\frac{5}{2}}\frac{5c_0C_{2,c_0}\norm{\tilde{Q}_{c_0}^{\frac{3}{2}}\cos \frac{y}{L}}_{L^2}^2}{48\norm{\tilde{Q}_{c_0}}_{L^2}^2}|\vec{a}(u)|^4 +\Bigl(1-\frac{c_m}{c_0}  \Bigr)\norm{\partial_y \Theta\bigl(\vec{a}(u),\gamma_{c_m}(\vec{a}(u))\bigr)}_{L^2}^2 \\ 
&+o(|\vec{a}(u)|^4+\norm{\eta(u_0)}_{H^1}^2)
\end{align*}
as $\delta \to 0$.
Since $|c_0-c_m|\lesssim \delta$ and
\begin{align*}
\partial_y \Theta\bigl(\vec{a}(u),\gamma_{c_m}(\vec{a}(u))\bigr)(x,y)&=- \frac{a_1(u)\gamma_{c_m}(\vec{a}(u))}{c_0L}\tilde{Q}_{c_0}^{\frac{3}{2}}\Bigl(\sqrt{\frac{\gamma_{c_m}(\vec{a}(u))}{c_0}}x,y \Bigr)\sin \frac{y}{L} \\
&+  \frac{a_2(u)\gamma_{c_m}(\vec{a}(u))}{c_0L}\tilde{Q}_{c_0}^{\frac{3}{2}}\Bigl(\sqrt{\frac{\gamma_{c_m}(\vec{a}(u))}{c_0}}x,y \Bigr)\cos \frac{y}{L} +O(|\vec{a}(u)|^2),
\end{align*}
there exist $k_3,k_4>0$ such that $k_3$ and $k_4$ are not depend on $c_m$, and 
\begin{align}\label{eq-or-1}
S_{c_m}(u)-S_{c_m}(\tilde{Q}_{c_m}) \geq k_3 \norm{\eta(u)}_{H^1}^2 + k_3|\vec{a}(u)|^2(|\vec{a}(u)|^2-\delta k_4 )+ o(|\vec{a}(u)|^4+\norm{\eta(u)}_{H^1}^2).
\end{align}
Using the conservation laws and \eqref{eq-5-5} , we obtain 
\begin{align}\label{eq-or-2}
S_{c_m}(u)-S_{c_m}(\tilde{Q}_{c_m})&=S_{c_m}(u_0)-S_{c_m}(\tilde{Q}_{c_m}) \notag \\
\lesssim & \norm{\eta(u_0)}_{H^1}^2+|\vec{a}(u_0)|^2 \lesssim \delta^2.
\end{align} 
From \eqref{eq-or-1} and \eqref{eq-or-2}, we have that there exist $\delta_*,k_5>0$ such that if $0<\delta<\delta_*$, then 
\[\norm{\eta(u)}_{H^1}^2 + |\vec{a}(u)|^2(|\vec{a}(u)|^2-\delta k_4 )-k_5\delta^2 \leq 0.\]
Therefore, there exists $C(k_4,k_5)>0$ such that
\[ |\vec{a}(u)|^2+\norm{\eta(u)}_{H^1}\leq C(k_4,k_5) \delta.\]
Applying \eqref{eq-est-c-1}, we have
\[|c_0-c(u)|\lesssim \norm{\eta(u)}_{L^2}^2+|\gamma_{c_m}(\vec{a}(u))-c_m| + |c_0-c_m| \lesssim \delta.\]

\qed

\section{Liouville property}
In this section, we prove the Liouville property of \eqref{ZKeq}.
First, we show the following equation of the integration of $Q_{c}$.
\begin{lemma}\label{lem-Q-pro}
Let $p,c>0$.
Then, we have
\begin{align}\label{Q-pro}
 \int_{\R} Q_{c}^{p+1} dx = \frac{3pc}{2p+1 }\int_{\R} Q_ c^p dx.
 \end{align}
\end{lemma}
\proof
Since
\begin{align}\label{Q-eq-1}
 -\partial_x^2 Q_c +cQ_c-Q_c^2=0,
\end{align}
we have
\begin{align*}
\int_{\R} Q_c^{p+1}dx=& - \int_{\R} Q_c^{p-1} \partial_x^2 Q_c dx +c \int_{\R} Q_c^pdx\\
=& (p-1) \int_{\R}Q_c^{p-2}(\partial_x Q_c)^2 dx +c \int_{\R}Q_c^pdx.
\end{align*}
 Multiplying \eqref{Q-eq-1} by $\partial_x Q_c$ and integrating this, we obtain 
 \[ -(\partial_x Q_c)^2 +c Q_c^2 - \frac{2}{3}Q_c^3 =0.\]
 Thus,
 \begin{align}\label{Q-eq-2}
\int_{\R} Q_c^{p+1}dx=(p-1) \int_{\R}Q_c^{p-2}\Bigl(c Q_c^2 - \frac{2}{3}Q_c^3\Bigr) dx  +c \int_{\R}Q_c^pdx
\end{align}
 which implies \eqref{Q-pro}.

\qed

Let 
\[\phi_c(x)=-\frac{\partial_x Q_c(x)}{Q_c(x)}=\sqrt{c}\tanh \frac{\sqrt{c}x}{2}.\]
Then, $\phi_c(x) \to \pm \sqrt{c}$ as $x \to \pm \infty$ and 
\[\partial_x  \phi_c(x)=\frac{c}{2} \cosh^{-2} \frac{\sqrt{c}x}{2}=\frac{1}{3}Q_c.\]
We introduce the following coerciveness type lemma in \cite{M M 3}.
\begin{lemma}\label{lem-coer}
For $u \in H^1(\R)$ 
\begin{align*}
-\int_{\R} \partial_x u \mathcal{L}_c(u\phi_c)dx=& \frac{3}{2}\int_{\R}\Bigl(\partial_x\Bigl( \frac{u}{Q_c} \Bigr)\Bigr)^2 Q_c^2 \partial_x \phi_cdx\\
\geq & \frac{5c}{8}\Bigl( \int_{\R} 3|u|^2 \partial_x\phi_cdx -\norm{Q_{c}}_{L^3(\R)}^{-3}\Bigl(\int_{\R}u Q_c^2 dx \Bigr)^2\Bigr)
\end{align*}
\end{lemma}
\proof


Let $v=\frac{u}{Q_c}$.
Since 
\begin{align*}
\mathcal{L}_c(u\phi_c)=\mathcal{L}_c(v\partial_xQ_c)=-2\partial_xv \partial_x^2Q_c-\partial_x^2 v \partial_xQ_c,
\end{align*}
we have
\begin{align*}
-\int_{\R}\partial_x u \mathcal{L}_c(u\phi_c)dx=\int_{\R}\partial_x(Q_cv)\mathcal{L}_c(v\partial_xQ_c)dx
=\frac{1}{2}\int_{\R}(\partial_xv)^2Q_c^3dx
\end{align*}
Let $ w=vQ_c^{\frac{3}{2}}$.
Using
\[\partial_x^2 Q_c=cQ_c -Q_c^2, \quad (\partial_x Q_c)^2 = c Q_c^2-\frac{2}{3} Q_c^3,\]
we obtain that
\begin{align*}
\frac{1}{2}\int_{\R}(\partial_xv)^2Q_c^3dx
=&\frac{1}{2}\int_{\R}w\Bigl( -\partial_x^2w + \frac{3}{2} \partial_x^2Q_cQ_c^{-1}w +\frac{3}{4}(\partial_x Q_c)^2 Q_c^{-2}w \Bigr)dx\\
=&\frac{1}{2}\int_{\R}w\Bigl(\mathcal{L}_c+\frac{5c}{4}\Bigr)wdx
\end{align*}
From the properties of $\mathcal{L}_c$ the operator $\mathcal{L}_c+\frac{5c}{4}$ is non-negative and the kernel of $\mathcal{L}_c+\frac{5c}{4}$ is spanned by $ Q_c^{\frac{3}{2}}$.
Moreover, the second eigenvalue of $\mathcal{L}_c+\frac{5c}{4}$ is $\frac{5c}{4}$.
Therefore, we have
\begin{align*}
\int_{\R}w\Bigl(\mathcal{L}_c+\frac{5c}{4}\Bigr)wdx \geq & \frac{5c}{4}\Bigl( \norm{w}_ {L^2}^2-\norm{Q_c}_{L^3(\R)}^{-3}\Bigl(\int_{\R}w Q_c^{\frac{2}{3}} dx \Bigr)^2\Bigr)\\
=& \frac{5c}{4}\Bigl( \int_{\R} 3|u|^2 \partial_x\phi_cdx -\norm{Q_c}_{L^3(\R)}^{-3}\Bigl(\int_{\R}u Q_c^2 dx \Bigr)^2\Bigr)
\end{align*}

\qed

\subsection{Monotonicity properties}
In this subsection, we show the monotonicity properties of \eqref{ZKeq}.
By Proposition \ref{prop-WPW}, the equation \eqref{ZKeq} has the Kato type local smoothing effect.
Therefore, the proof of the monotonicity properties is similar to one in \cite{M M 3,C M P S}.
Thus, we omit the detail of proofs in this subsection, see Section 3 in \cite{C M P S}.

We define $\psi_R \in C^{\infty}(\R,\R)$ by
\begin{align}\label{def-psi}
\psi_R(x)=\frac{2}{\pi}\arctan(e^{x/R}), \quad x \in \R.
\end{align}
Then, we have $\lim_{x \to \infty} \psi_R(x)=1$, $\lim_{x \to -\infty} \psi_R(x)=0$,
\[\partial_x \psi_R(x)=\frac{1}{\pi R \cosh(x/R)} \mbox{ and } |\partial_x^3\psi_R(x)| \leq \frac{1}{R^2} \partial_x \psi_R(x).\]
Let $\varepsilon ,\beta ,c_0>0$ and $u$ be a solution to \eqref{ZKeq} satisfying that there exists $\rho \in C(\R,\R)$ such that 
\begin{align}\label{eq-6-1}
\norm{u(t,\cdot,\cdot)- \tilde{Q}_{c_0}(\cdot- \rho(t),\cdot)}_{H^1}<\varepsilon_0 , \quad t \in \R
\end{align}
and 
\begin{align}\label{eq-6-1-1}
|\dot{\rho}(t)-c_0|\leq c_0/2, \quad t \in \R.
\end{align}
For $x_0,t_0, t \in \R$ we define
\[ \tilde{x}=\tilde{x}(x_0,t_0,t)=x-\rho(t_0)+\frac{\beta(t_0-t)}{2}-x_0,\]
\[\tilde{x}_- = \tilde{x}(-x_0,t,t_0,)\]
\[I_{x_0,t_0}(u(t))=\int_{\RTL} |u(t,x,y)|^2 \psi_R(\tilde{x}(t))dxdy,\]
and 
\[ I_{x_0,t_0}^-(u(t))=\int_{\RTL} |u(t,x,y)|^2 \psi_R(\tilde{x}_-(t))dxdy.\]

In the following lemma, we show the property of the parameter $\rho$ (see Lemma 3.2 in \cite{C M P S}).
\begin{lemma}\label{lem-rho-pro}
Assume that $u \in C(\R,H^1(\RTL))$ is a solution to $\eqref{ZKeq}$ satisfying $\eqref{eq-6-1}$, $\eqref{eq-6-1-1}$ and 
that there exist $\tilde{\rho} \in C(\R,\R)$ and $ C, \delta_0>0$ such that
\begin{align}\label{ass-decay-1}
\int_{\T_L} |u(t,x+\tilde{\rho}(t),y)|^2dy \leq C e^{-\delta_0|x|}, \quad (t,x) \in \R^2.
\end{align}
If $0<\varepsilon _0< \frac{1}{2}\norm{\tilde{Q}_{c_0}}_{L^2(|x|\leq 1)}$,
 then $u$ satisfies 
\begin{align}\label{ass-decay}
\int_{\T_L} |u(t,x+\rho(t),y)|^2dy \lesssim e^{-\delta_0|x|}, \quad (t,x) \in \R^2,
\end{align}
where
\[\norm{u}_{L^2(|x|\leq R)}^2=\int_{|x|\leq R} |u|^2dxdy.\]
\end{lemma}


The following two lemmas show the $L^2$-monotonicity property of \eqref{ZKeq}.
\begin{lemma}\label{lem-L2-mono}
Let $0< \beta < c_0/2 $.
Assume that $u \in C(\R,H^1(\RTL))$ is a solution to $\eqref{ZKeq}$ satisfying $\eqref{eq-6-1}$ and $\eqref{eq-6-1-1}$.
Then, for $x_0>0,t_0 \in \R $, $R\geq 2/\sqrt{\beta}$ and $t \leq t_0$
\begin{align}\label{L2-mono-1}
 I_{x_0,t_0}(u(t_0))-I_{x_0,t_0}(u(t)) \lesssim e^{-x_0/R},
\end{align}
if $\varepsilon _0>0$ in $\eqref{eq-6-1}$ is chosen small enough.
Moreover, if $u$ satisfies the decay assumption $\eqref{ass-decay}$, then
\begin{align}\label{L2-mono-2}
\int_{\RTL}|u(t_0,x,y)|^2&\psi_R(\tilde{x}(t_0))dxdy\notag \\
&+\int_{-\infty}^{t_0} \int_{\RTL}(|\nabla u|^2 +|u|^2)(t,x,y) \partial_x\psi_R(\tilde{x}(t))dxdydt \lesssim e^{-x_0/R}.
\end{align}
\end{lemma}

\begin{lemma}\label{lem-L2-mono-}
Let $0< \beta < c_0/2 $.
Assume that $u \in C(\R,H^1(\RTL))$ is a solution to $\eqref{ZKeq}$ satisfying $\eqref{eq-6-1}$ and $\eqref{eq-6-1-1}$.
Then, for $x_0>0,t_0 \in \R $, $R\geq 2/\sqrt{\beta}$ and $t \geq t_0$
\begin{align}\label{L2-mono-1-}
 I_{x_0,t_0}^-(u(t))-I_{x_0,t_0}^-(u(t_0)) \lesssim e^{-x_0/R},
\end{align}
if $\varepsilon _0>0$ in $\eqref{eq-6-1}$ is chosen small enough.
\end{lemma}
The proof of Lemma \ref{lem-L2-mono} follows the proof of Lemma  3.3 in \cite{C M P S}.
The proof of Lemma \ref{lem-L2-mono-} is similar to the proof of Lemma 4.9 in \cite{C M P S}.


We define a functional $J$ by 
\[J_{x_0,t_0}(u(t))=\int_{\RTL} \Bigl(|\nabla u|^2 -\frac{2}{3}u^3\Bigr)(t,x,y) \psi_R(\tilde{x})dxdy.\]
In the following lemma, we show the monotonicity property for $J$.
\begin{lemma}\label{lem-E-mono}
Let $ 0<\beta<c_0/2$.
Assume that $u \in C(\R, H^1(\RTL))$ is a solution to $\eqref{ZKeq}$ satisfying $\eqref{eq-6-1}$ and $\eqref{eq-6-1-1}$.
Then, for $x_0>0$, $t_0 \in \R$, $R\geq 2/\sqrt{\beta}$ and $t \leq t_0$
\begin{align}\label{E-mono-1}
J_{x_0,t_0}(u(t_0))-J_{x_0,t_0}(u(t)) \lesssim e^{-x_0/R}.
\end{align}
Moreover,  if $u$ satisfies the decay assumption $\eqref{ass-decay}$, then 
\begin{align}\label{E-mono-2}
&\int_{\RTL} |\nabla u|^2 (t_0) \psi_R(\tilde{x}(t_0))dxdy \\
&+ \int_{-\infty}^{t_0}\int_{\RTL} (|\nabla^2 u(t)|^2 + |\nabla u(t)|^2 +u(t)^4)(\partial_x \psi_R)(\tilde{x}(t)) dxdydt \lesssim e^{-x_0/R}.
\end{align}
\end{lemma}
The proof of Lemma \ref{lem-E-mono} is similar to the proof of Lemma 3.4 in \cite{C M P S}.

The following proposition shows the boundedness of higher Sobolev norm of solutions satisfying the decay assumption \eqref{ass-decay}.
\begin{proposition}\label{cor-high-bdd}
Let $ 0<\beta<c_0/2$ and $k \in \Z_+$.
Assume that $u \in C(\R, H^1(\RTL))$ is a solution to $\eqref{ZKeq}$ satisfying $\eqref{eq-6-1}$, $\eqref{eq-6-1-1}$ and the decay assumption $\eqref{ass-decay}$.
If $\varepsilon _0>0$ in $\eqref{eq-6-1}$ is sufficiently small, there exist $\tilde{\delta}, C=C(k)>0$ such that 
\begin{align}\label{high-bdd}
 \sup_{t \in \R} \int_{\RTL} (\partial^{\alpha} u)^2(x+\rho(t),y) e^{\tilde{\delta}|x|} dxdy \leq C,
 \end{align}
for $\alpha \in (\N_0)^2 $ satisfying $|\alpha|\leq k$.
\end{proposition}
The proof of this proposition is same as the proof of Corollary 3.9 in \cite{C M P S}.


\subsection{Critical case $L= \frac{2}{\sqrt{5c_0}}$}
In this section, we show the Liouville property for  $L = \frac{2}{\sqrt{5c_0}}$.

\begin{lemma}\label{lem-c-mod}
There exist $\varepsilon _0,K_0>0$ such that for any $0<\varepsilon <\varepsilon _0$ the following is true.
For any solution $u \in C(\R,H^1(\RTL))$ of $\eqref{ZKeq}$ satisfying 
\begin{align*}
\inf_{b \in \R} \norm{u(t,\cdot,\cdot)-Q_{c_0}(\cdot-b,\cdot)}_{H^1}\leq \varepsilon 
\end{align*}
there exist $ \vec{a}=(a_1,a_2) \in C^1(\R,\R^2)$ and $\rho, c \in C^1(\R,\R)$  uniquely such that 
\begin{align}\label{eq-6-c-1}
\eta(t,x,y)=u(t,x+\rho(t),y)-\Theta(\vec{a}(t),c(t))
\end{align}
satisfies for all $t\in \R$ 
\begin{align}
&|c(t)-c_0|+|a_1(t)|+|a_2(t)|+\norm{\eta(t)}_{H^1}\leq K_0\varepsilon, \label{eq-6-c-2} \\
\int_{\RTL} \eta(t)\partial_x\Theta(\vec{a}(t),c(t))&dxdy=\int_{\RTL}\eta(t)\Theta(\vec{a}(t),c(t))dxdy\notag\\
=\int_{\RTL}\eta(t)\Theta(\vec{a}(t),c(t))^{\frac{3}{2}}&\cos \frac{y}{L}dxdy
=\int_{\RTL}\eta(t)\Theta(\vec{a}(t),c(t))^{\frac{3}{2}}\sin \frac{y}{L}dxdy=0 \label{eq-6-c-3}
\end{align} 
and
\begin{align}
|\dot{\vec{a}}(t)| \leq & K_0\norm{\eta(t)}_{L^2}, \label{eq-6-c-2-1}\\
|\dot{c}(t)|\leq & \varepsilon K_0\norm{\eta(t)}_{L^2}.\label{eq-6-c-2-2}\\
|\dot{\rho}(t)-\hat{c}(t)| \leq & K_0 (\norm{\eta(t)}_{L^2}+|c-c_0||\vec{a}|)\label{eq-6-c-2-3},
\end{align}
where $\hat{c}(t)=c_0^{-1}c(t)\check{c} (\vec{a}(t))$.
\end{lemma}
\proof
From Lemma \ref{lem-p-3-3}, there exist $C^1$ mappings $\rho(t)=\rho(u(t)),c(t)=c(u(t)), \vec{a}(t)=\vec{a}(u(t)), \eta(t)=\eta(u(t))$ satisfying \eqref{eq-6-c-1}--\eqref{eq-6-c-3}.
By the calculation we have
\begin{align}\label{eq-6-c-5}
 \eta_t
=& \partial_x(-\Delta \eta -2 \Theta \eta -\eta^2 - \Delta \Theta - \Theta^2) + \dot{\rho}\partial_x(\eta+\Theta) -\dot{\vec{a}}\cdot \partial_{\vec{a}}\Theta  -\dot{c} \partial_c\Theta\notag \\
=& \partial_x(\mathbb{L}_{c_0} \eta -\eta^2) + \partial_x S_{c_0}'(\Theta) - 2\partial_x((\Theta-\tilde{Q}_{c_0})\eta) + (\dot{\rho}-c_0) \partial_x(\eta+\Theta) -\dot{\vec{a}}\cdot \partial_{\vec{a}}\Theta  -\dot{c} \partial_c\Theta,
\end{align}
where $S_{c_0}'(\Theta)= - \Delta \Theta +c_0\Theta - \Theta^2$ and $\dot{\vec{a}}\cdot \partial_{\vec{a}}\Theta=\dot{a}_1 \partial_{a_1}\Theta + \dot{a}_2 \partial_{a_2} \Theta$.
From \eqref{eq-6-c-3} and $\Theta (x,y)=\Theta (-x,y)$ we obtain that
\begin{align*}
0=& \frac{d}{d t} \int_{\RTL}\eta\Theta dxdy \\
=&- \int_{\RTL} \Theta \dot{\vec{a}}\cdot \partial_{\vec{a}}\Theta dxdy  -\dot{c}\int_{\RTL}\tilde{Q}_{c_0} \partial_c\tilde{Q}_{c_0} dxdy \\
&+O\bigl(\norm{ \eta}_{L^2} (\norm{\eta}_{L^2} +|\vec{a}|+|c-c_0|+|\dot{\vec{a}}| + |\dot{c}|) \bigr).
\end{align*}
By the expansion
\[ \Theta(\vec{a},c)=\tilde{Q}_{c_0} + O(|\vec{a}|+|c-c_0|), \quad \dot{\vec{a}}\cdot \partial_{\vec{a}}\Theta(\vec{a},c)= \dot{a}_1 Q_{c_0}^{\frac{3}{2}} \cos \frac{y}{L} + \dot{a}_2 Q_{c_0}^{\frac{3}{2}} \sin \frac{y}{L} + O((|\vec{a}|+|c-c_0|)|\dot{\vec{a}}|),\]
we have 
\[\int_{\RTL} \Theta \dot{\vec{a}}\cdot \partial_{\vec{a}}\Theta dxdy=O((|\vec{a}|+|c-c_0|)|\dot{\vec{a}}|).\]
Since $\int_{\RTL}\tilde{Q}_{c_0} \partial_c\tilde{Q}_{c_0} dxdy \neq 0$, 
\begin{align}\label{eq-6-c-6}
|\dot{c}|=O\bigl(\norm{ \eta}_{L^2} (\norm{ \eta}_{L^2} +|\vec{a}|+|c-c_0|+|\dot{\vec{a}}| )+(|\vec{a}|+|c-c_0|)|\dot{\vec{a}}| \bigr).
\end{align}
From \eqref{eq-6-c-3}, \eqref{eq-6-c-6} and $\Theta (x,y)=\Theta (-x,y)$, we obtain that 
\begin{align*}
0=& \frac{d}{d t} \int_{\RTL}\eta\Theta^{\frac{3}{2}}\cos \frac{y}{L} dxdy \\
=&- \dot{a}_1 \int_{\RTL}\Bigl( \Theta^{\frac{3}{2}}\cos \frac{y}{L}\Bigr)  \partial_{a_1}\Theta dxdy +O\bigl(\norm{ \eta}_{L^2} + (|\vec{a}|+|c-c_0|)|\dot{\vec{a}}|\bigr).
\end{align*}
Since
\[ \int_{\RTL}\Bigl( \Theta^{\frac{3}{2}}\cos \frac{y}{L}\Bigr)\partial_{a_1}\Theta dxdy= \int_{\RTL}\Bigl( \tilde{Q}_{c_0}^{\frac{3}{2}}\cos \frac{y}{L}\Bigr)^2dxdy +O (|\vec{a}|+|c-c_0|),\]
we obtain 
\begin{align}\label{eq-6-c-7}
|\dot{a}_1|=O\bigl(\norm{ \eta}_{L^2} + (|\vec{a}|+|c-c_0|)|\dot{a}_2|\bigr).
\end{align}
By the same way, from \eqref{eq-6-c-6} and \eqref{eq-6-c-7} we get
\begin{align}\label{eq-6-c-8}
|\dot{a}_2|=O\bigl(\norm{ \eta}_{L^2}\bigr).
\end{align}
The estimates \eqref{eq-6-c-2-1} and \eqref{eq-6-c-2-2} follow \eqref{eq-6-c-6}--\eqref{eq-6-c-8}.
By the similar computation to \eqref{eq-6-c-5}, 
\begin{align}\label{eq-6-c-9}
\eta_t=\partial_x(\mathbb{L}_{\hat{c}} \eta -\eta^2) + \partial_x S_{\hat{c}}'(\Theta) - 2\partial_x((\Theta-\tilde{Q}_{\hat{c}})\eta) + (\dot{\rho}-\hat{c}) \partial_x(\eta+\Theta) -\dot{\vec{a}}\cdot \partial_{\vec{a}}\Theta  -\dot{c} \partial_c\Theta.
\end{align}
By the definition of $\Theta$ and $\hat{c}$ we have
\begin{align}\label{eq-6-c-10}
S_{\hat{c}}'(\Theta)
=& \frac{c^2}{c_0^2}(-\Delta \varphi_{c_0}+ \check{c}  \varphi_{c_0} -(\varphi_{c_0})^2) + \frac{c(c-c_0)}{c_0^2} \partial_y^2 \varphi_{c_0}\notag \\
=&\frac{c-c_0}{c_0} \partial_y^2\Theta.
\end{align}
Since 
\begin{align}\label{eq-6-p}
S_{\hat{c}}'(\Theta)=\frac{c-c_0}{c_0} \partial_y^2\Theta = O(|c_0-c||\vec{a}|),
\end{align}
from \eqref{eq-6-c-3}, \eqref{eq-6-c-6}--\eqref{eq-6-c-8}, we obtain that
\begin{align*}
0=& \frac{d}{d t} \int_{\RTL}\eta(\partial_x \Theta) dxdy \\
=&(\dot{\rho}-\hat{c}) \int_{\RTL}(\partial_x \tilde{Q}_{c_0})^2 dxdy+O\bigl(\norm{\eta}_{L^2} +|c_0-c||\vec{a}|\bigr).
\end{align*}
Thus, the estimate \eqref{eq-6-c-2-3} holds.

\qed

Next we prove the following Liouville type theorem.
\begin{theorem}\label{thm-c-Liouville}
Let $c_0>0$ and $L=\frac{2}{\sqrt{5c_0}}$.
There exists $\varepsilon  _0>0$ satisfies the following.
For any solution $u \in C(\R,H^1(\RTL))$ to $\eqref{ZKeq}$ satisfying $\eqref{eq-6-1}$ and $\eqref{ass-decay}$, there exist $c_+>0, \vec{a}_+=(a_{1,+},a_{2,+})$ and $\rho_0$ such that
\[u(t,x,y)=\Theta(\vec{a}_+,c_+)(x-\hat{c}_+t+\rho_0,y),\]
\[|c_+-c_0||\vec{a}_+|=0,\]
where 
\[\hat{c}_+=
\begin{cases}
c_+, & \vec{a}_+=(0,0),\\
\check{c} (\vec{a}_+) , & c_+=c_0.
\end{cases}\]
\end{theorem}

\proof

Let $u \in C(\R,H^1(\RTL))$ be solution to \eqref{ZKeq} satisfying \eqref{eq-6-1} and \eqref{ass-decay}.
From Lemma \ref{lem-rho-pro} and \ref{lem-c-mod}, $\rho$ in Lemma \ref{lem-c-mod} satisfies \eqref{eq-6-1} and \eqref{ass-decay}.
Let $\eta(t),c(t),\vec{a}(t),\hat{c}(t)$ be in Lemma \ref{lem-c-mod}.
We define
\[ v=\mathbb{L}_{\hat{c}}\eta-\eta^2.\]
Then, $v$ has the following almost orthogonal condition.
\begin{align}\label{con-ort-1}
\int_{\RTL}v\partial_x \tilde{Q}_{\hat{c}}dxdy=\int_{\RTL}(\mathbb{L}_{\hat{c}}\eta-\eta^2)\partial_x \tilde{Q}_{\hat{c}}dxdy = O(\norm{\eta}_{L^2}^2),
\end{align}
\begin{align}\label{con-ort-2}
\int_{\RTL} v \partial_c\tilde{Q}_{\hat{c}}dxdy
=-\int_{\RTL}\eta\tilde{Q}_{\hat{c}}dxdy+O(\norm{\eta}_{L^2}^2)
=O(\norm{\eta}_{L^2}(|c-c_0|+|\vec{a}|+\norm{\eta}_{L^2})),
\end{align}
\begin{align}\label{con-ort-3}
\int_{\RTL} v \tilde{Q}_{\hat{c}}^{\frac{3}{2}} \cos \frac{y}{L}dxdy
=& \int_{\RTL}\eta \mathbb{L}_{\hat{c}}\biggl(\tilde{Q}_{\hat{c}}^{\frac{3}{2}} \cos \frac{y}{L} \biggr)dxdy+O(\norm{\eta}_{L^2}^2) \notag \\
=&O(\norm{\eta}_{L^2}(|c-c_0|+|\vec{a}|+\norm{\eta}_{L^2})),
\end{align}
\begin{align}\label{con-ort-4}
\int_{\RTL} v \tilde{Q}_{\hat{c}}^{\frac{3}{2}} \sin \frac{y}{L}dxdy=O(\norm{\eta}_{L^2}(|c-c_0|+|\vec{a}|+\norm{\eta}_{L^2})).
\end{align}
From the orthogonal conditions \eqref{eq-6-c-3} and Lemma \ref{lem-p-coer}, we have 
\begin{align*}
(v,\eta)_{L^2(\RTL)}=(\mathbb{L}_{\hat{c}}\eta,\eta)_{L^2(\RTL)}-\norm{\eta}_{L^3(\RTL)}^3 \geq k_2 \norm{\eta}_{H^1}^2 +O(\norm{\eta}_{H^1}^3).
\end{align*}
Therefore, if $\varepsilon _0>0$ is sufficiently small, then for $t \in \R$
\begin{align}\label{con-coer}
\norm{\eta}_{H^1}\lesssim \norm{v}_{L^2}.
\end{align}
By \eqref{eq-6-c-9}, we have
\begin{align}\label{eq-c-li}
v_t=&\mathbb{L}_{\hat{c}}\eta_t+(\partial_t \mathbb{L}_{\hat{c}})\eta-2\eta\eta_t
= \mathbb{L}_{\hat{c}}\partial_x v  + \mathbb{L}_{\hat{c}} \partial_xS_{\hat{c}}'(\Theta)  +R(\eta,\vec{a},c),
\end{align}
where
\begin{align*}
R(\eta,\vec{a},c)=&-2\eta\partial_x(v+S_{\hat{c}}'(\Theta)) + (\dot{\rho}-\hat{c}) \mathbb{L}_{\hat{c}}\partial_x(\eta+\Theta) -2  (\dot{\rho}-\hat{c})\eta \partial_x(\eta+\Theta)\\
&+(\mathbb{L}_{\hat{c}}-2\eta)( 2\partial_x((\tilde{Q}_{\hat{c}}-\Theta)\eta) - \dot{\vec{a}}\cdot \partial_{\vec{a}}\Theta-\dot{c}\partial_c\Theta) +(\partial_t \mathbb{L}_{\hat{c}})\eta.
\end{align*}
Therefore, 
\begin{align}\label{eq-c-li-1}
&-\frac{1}{2} \frac{d}{d t} \int_{\RTL}(v+S_{\hat{c}}'(\Theta))^2 \phi_{\hat{c}}dxdy \notag \\
=& - \int_{\RTL} (\mathbb{L}_{\hat{c}}\partial_x v)v \phi_{\hat{c}}dxdy  - \int_{\RTL} (\mathbb{L}_{\hat{c}} \partial_xS_{\hat{c}}'(\Theta))v \phi_{\hat{c}}dxdy - \int_{\RTL} (\mathbb{L}_{\hat{c}}\partial_x v)S_{\hat{c}}'(\Theta) \phi_{\hat{c}}dxdy\notag \\
&-\int_{\RTL} (\mathbb{L}_{\hat{c}} \partial_xS_{\hat{c}}'(\Theta))S_{\hat{c}}'(\Theta) \phi_{\hat{c}}dxdy -\int_{\RTL}R(\eta,\vec{a},c)(v+S_{\hat{c}}'(\Theta))\phi_{\hat{c}}dxdy \notag \\
&- \int_{\RTL}(\partial_t S_{\hat{c}}'(\Theta))(v+S_{\hat{c}}'(\Theta))\phi_{\hat{c}}dxdy-\frac{1}{2} \int_{\RTL}(v+S_{\hat{c}}'(\Theta))^2 \dot{\hat{c}}\partial_c\phi_{\hat{c}}dxdy.
\end{align}
We estimate each term in \eqref{eq-c-li-1} separately.

(I) The estimate of $- \int_{\RTL} (\mathbb{L}_{\hat{c}}\partial_x v)v \phi_{\hat{c}}dxdy $.
From the Fourier expansion $v(t,x,y)=v_0 (t,x) + \sum_{n =1}^{\infty} (v_{n,1}(t,x) \cos \frac{y}{L}+v_{n,2}(t,x) \sin \frac{y}{L})$ and Lemma \ref{lem-coer}  we have
\begin{align}\label{eq-c-li-2}
&- \int_{\RTL} (\mathbb{L}_{\hat{c}}\partial_x v)v \phi_{\hat{c}}dxdy\notag \\
=& -2\pi L\int_{\RTL}\Bigl(\mathcal{L}_{\hat{c}}\partial_x v_0\Bigr) v_0 \phi_{\hat{c}} dx - \pi L \sum_{n \in \Z_+, j=1,2}\int_{\R}   \Bigl(\Bigl(\mathcal{L}_{\hat{c}}+\frac{n^2}{L^2} \Bigr)\partial_x v_{n,j}\Bigr) v_{n,j} \phi_{\hat{c}} dx \notag \\
\geq & \frac{5\pi L\hat{c}}{4}\Bigl( \int_{\R} |v_0|^2 Q_{\hat{c}} dx- \norm{Q_{\hat{c}}}_{L^3(\R)}^{-3} \Bigl( \int_{\R}v_0 Q_{\hat{c}}^2dx \Bigr)^2 \Bigr) \notag \\
&+ \pi L\Bigl( \Bigl(\frac{5\hat{c}}{8}+ \frac{1}{6L^2} \Bigr)\int_{\R} (|v_{1,1}|^2+|v_{1,2}|^2) Q_{\hat{c}} dx - \frac{5\hat{c}}{8}\norm{Q_{\hat{c}}}_{L^3(\R)}^{-3} \Bigl( \int_{\R}v_{1,1} Q_{\hat{c}}^2dx \Bigr)^2  \notag \\
&- \frac{5\hat{c}}{8}\norm{Q_{\hat{c}}}_{L^3(\R)}^{-3} \Bigl( \int_{\R}v_{1,2} Q_{\hat{c}}^2dx \Bigr)^2\Bigr) \notag \\
&+ \pi L \sum_{n=2}^{\infty}  \frac{n^2}{6L^2} \int_{\R} (|v_{n,1}|^2+|v_{n,2}|^2) Q_{\hat{c}} dx.
\end{align}
From the almost orthogonal condition \eqref{con-ort-2} 
\begin{align}\label{eq-c-li-3}
\Bigl| \int_{\R} v_0 Q_{c_0}^2 dx\Bigr|^2 \leq & \Bigl| \int_{\R} v_0\Bigl( Q_{c_0}^2 - \norm{Q_{c_0}^{-\frac{1}{2}} \partial_c Q_{c_0}}_{L^2(\R)}^{-2} \int_{\R} Q_{c_0}\partial_c Q_{c_0} dx' \partial_c Q_{c_0} \Bigr)dx\Bigr|^2 \notag \\
&+ O((|\vec{a}|+|c-c_0|+\norm{\eta}_{L^2})\norm{\eta}_{L^2}\norm{v}_{L^2}) \notag \\
\leq & \norm{v_0 Q_{c_0}^{\frac{1}{2}}}_{L^2(\R)}^2\norm{Q_{c_0}^{\frac{3}{2}} - \norm{Q_{c_0}^{-\frac{1}{2}} \partial_c Q_{c_0}}_{L^2(\R)}^{-2} \int_{\R} Q_{c_0}\partial_c Q_{c_0} dx' Q_{c_0}^{-\frac{1}{2}}\partial_c Q_{c_0} }_{L^2(\R)}^2\notag \\
&+ O((|\vec{a}|+|c-c_0|+\norm{\eta}_{L^2})\norm{\eta}_{L^2}\norm{v}_{L^2})
\end{align}
Since
\begin{align*}
 &\norm{Q_{c_0}^{\frac{3}{2}} - \norm{Q_{c_0}^{-\frac{1}{2}} \partial_c Q_{c_0}}_{L^2(\R)}^{-2} \int_{\R} Q_{c_0}\partial_c Q_{c_0} dx' Q_{c_0}^{-\frac{1}{2}}\partial_c Q_{c_0} }_{L^2(\R)}^2
 <  \norm{Q_{c_0}}_{L^3(\R)}^3,
\end{align*}
from \eqref{eq-c-li-3} we obtain 
\begin{align}\label{eq-c-li-4}
&\int_{\R} |v_0|^2 Q_{c_0} dx - \norm{Q_{c_0}}_{L^3(\R)}^{-3} \Bigl( \int_{\R}v_0 Q_{c_0}^2dx \Bigr)^2 \notag \\
\gtrsim & \int_{\R} |v_0|^2 Q_{c_0} dx+ O((|\vec{a}|+|c-c_0|+\norm{\eta}_{L^2})\norm{\eta}_{L^2}\norm{v}_{L^2}).
\end{align}
On the other hand, from the almost orthogonal conditions \eqref{con-ort-3} and \eqref{con-ort-4} 
\begin{align}\label{eq-c-li-5}
\Bigl| \int_{\R}v_{1,j} Q_{c_0} ^2 dx\Bigr|^2 =& \Bigl| \int_{\R} v_{1,j} \Bigl( Q_{c_0}^{\frac{3}{2}} -\norm{Q_{c_0}}_{L^2(\R)}^{-2}  \int_{\R} Q_{c_0}^{\frac{5}{2}} dx' Q_{c_0}\Bigr)Q_{c_0}^{\frac{1}{2}} dx \Bigr|^2 \notag \\
&+ O((|\vec{a}|+|c-c_0|+\norm{\eta}_{L^2})\norm{\eta}_{L^2}\norm{v}_{L^2})\notag \\
\leq & \norm{Q_{c_0}^{\frac{3}{2}}-\norm{Q_{c_0}}_{L^2(\R)}^{-2}  \int_{\R} Q_{c_0}^{\frac{5}{2}} dx' Q_{c_0}}_{L^2(\R)}^2\int_{\R} |v_{1,j}|^2 Q_{c_0} dx  \notag \\
&+ O((|\vec{a}|+|c-c_0|+\norm{\eta}_{L^2})\norm{\eta}_{L^2}\norm{v}_{L^2}).
\end{align}
By Lemma \ref{lem-Q-pro}, 
\begin{align}\label{eq-c-li-6}
\int_{\R} Q_{c_0}^{\frac{7}{2}}dx= \frac{5c_0}{4} \int_{\R} Q_{c_0}^{\frac{5}{2}}dx \quad \int_{\R} Q_{c_0}^{\frac{5}{2}}dx = \frac{9c_0}{8} \int_{\R} Q_{c_0}^{\frac{3}{2}}dx.
\end{align}
From H\"older's inequality, \eqref{eq-c-li-5} and \eqref{eq-c-li-6}, we obtain
\begin{align}\label{eq-c-li-7}
&\int_{\R} |v_{1,j}|^2 Q_{c_0} dx - \norm{Q_{c_0}}_{L^3(\R)}^{-3} \Bigl( \int_{\R}v_{1,j} Q_{c_0}^2dx \Bigr)^2 \notag \\
\geq & \norm{Q_{c_0}^{\frac{5}{4}}}_{L^2(\R)}^{-2}\norm{Q_{c_0}^{\frac{7}{4}}}_{L^2(\R)}^{-1} \norm{Q_{c_0}^{\frac{3}{4}}}_{L^2(\R)}^{-1} \Bigl( \int_{\R} Q_{c_0}^{\frac{5}{2}}dx \Bigr)^2 \int_{\R}|v_{1,j}|^2 Q_{c_0} dx \notag \\
&+ O((|\vec{a}|+|c-c_0|+\norm{\eta}_{L^2})\norm{\eta}_{L^2}^2)\notag \\
=& \sqrt{\frac{9}{10}} \int_{\R}|v_{1,j}|^2 Q_{c_0} dx + O((|\vec{a}|+|c-c_0|+\norm{\eta}_{L^2})\norm{\eta}_{L^2}\norm{v}_{L^2}).
\end{align}
By \eqref{eq-c-li-3}, \eqref{eq-c-li-4} and \eqref{eq-c-li-7}, we obtain that there exists $k_3>0$ such that
\begin{align}\label{eq-c-li-term1}
&- \int_{\RTL} (\mathbb{L}_{\hat{c}}\partial_x v)v \phi_{\hat{c}}dxdy\notag \\
\geq & k_3 \int_{\R} |v_0|^2 Q_{c_0} dx + \pi L\Bigl(\frac{5c_0}{8}\sqrt{\frac{9}{10}} + \frac{1}{6L^2} \Bigr)\int_{\R} (|v_{1,1}|^2+|v_{1,2}|^2) Q_{c_0} dx  \notag \\
&+ \pi L \sum_{n = 2}^{\infty}  \frac{n^2}{6L^2} \int_{\R}( |v_{n,1}|^2+|v_{n,2}|^2) Q_{c_0} dx \notag \\
&+ O((|\vec{a}|+|c-c_0|+\norm{\eta}_{L^2})\norm{\eta}_{H^1}\norm{v}_{L^2}) .
\end{align}

(II) The estimate of $-\int_{\RTL} (\mathbb{L}_{\hat{c}} \partial_xS_{\hat{c}}'(\Theta))S_{\hat{c}}'(\Theta) \phi_{\hat{c}}dxdy$.
Since
\begin{align}\label{eq-exp-theta}
S_{\hat{c}}'(\Theta)=\frac{c-c_0}{c_0} \Bigl( a_1 \tilde{Q}_{c_0}^{\frac{3}{2}} \cos \frac{y}{L}+a_2 \tilde{Q}_{c_0}^{\frac{3}{2}} \sin \frac{y}{L} \Bigr) + O((|c-c_0|+|\vec{a}|)|c-c_0||\vec{a}|),
\end{align}
we have
\begin{align}\label{eq-c-li-7-1}
&-\int_{\RTL} (\mathbb{L}_{\hat{c}} \partial_xS_{\hat{c}}'(\Theta))S_{\hat{c}}'(\Theta) \phi_{\hat{c}}dxdy \notag \\
=& -\frac{\pi L(c-c_0)^2 |\vec{a}|^2 }{c_0^2}\int_{\R} \Bigl((\mathcal{L}_{c_0} +\frac{1}{L^2}) \partial_x Q_{c_0}^{\frac{3}{2}} \Bigr) Q_{c_0}^{\frac{3}{2}} \phi_{c_0}dx+O((|c-c_0|+|\vec{a}|)|c-c_0|^2|\vec{a}|^2).
\end{align}
From \eqref{Q-eq-1} and \eqref{Q-eq-2} 
\begin{align}
\partial_x Q_{c_0}^{\frac{3}{2}} =& \frac{3}{2} Q_{c_0}^{\frac{1}{2}}\partial_x Q_{c_0} \label{Q-eq-t1},\\
\partial_x^2 Q_{c_0}^{\frac{3}{2}}
=& \frac{9c_0}{4}Q_{c_0}^{\frac{3}{2}} -2 Q_{c_0}^{\frac{5}{2}}, \label{Q-eq-t2}\\
\partial_x^3 Q_{c_0}^{\frac{3}{2}}=& \frac{27c_0}{8} Q_{c_0}^{\frac{1}{2}}\partial_xQ_{c_0}-5Q_{c_0}^{\frac{3}{2}}\partial_xQ_{c_0}, \label{Q-eq-t3}\\
\partial_x^4Q_{c_0}^{\frac{3}{2}}=& \frac{3}{2}\partial_x^3(Q_{c_0}^{\frac{1}{2}}\partial_xQ_{c_0}) 
=\frac{81c_0^2}{16}Q_{c_0}^{\frac{3}{2}} -17 c_0 Q_{c_0}^{\frac{5}{2}} + 10Q_{c_0}^{\frac{7}{2}}. \label{Q-eq-t4}
\end{align}
From \eqref{Q-eq-t1}--\eqref{Q-eq-t4} we have
\begin{align*}
&\Bigl(\Bigl( \mathcal{L}_{c_0} + \frac{1}{L^2} \Bigr)\partial_x Q_{c_0}^{\frac{3}{2}}\Bigr)Q_{c_0}^{\frac{3}{2}} \phi_{c_0} 
=-2c_0 Q_{c_0}^4+\frac{4}{3}Q_{c_0}^5.
\end{align*}
Applying Lemma \ref{lem-Q-pro}, we obtain 
\begin{align}\label{eq-c-li-term2}
&-\int_{\RTL} (\mathbb{L}_{\hat{c}} \partial_xS_{\hat{c}}'(\Theta))S_{\hat{c}}'(\Theta) \phi_{\hat{c}}dxdy \notag \\
=& \frac{\pi L(c-c_0)^2 |\vec{a}|^2 }{6c_0^2}\int_{\R} Q_{c_0}^5dx+O((|c-c_0|+|\vec{a}|)|c-c_0|^2|\vec{a}|^2).
\end{align}

(III) The estimate of $- \int_{\RTL} (\mathbb{L}_{\hat{c}} \partial_xS_{\hat{c}}'(\Theta))v \phi_{\hat{c}}dxdy- \int_{\RTL} (\mathbb{L}_{\hat{c}}\partial_x v)S_{\hat{c}}'(\Theta) \phi_{\hat{c}}dxdy$.
By the Fourier expansion and \eqref{eq-exp-theta}
\begin{align}\label{eq-c-li-8}
&- \int_{\RTL} (\mathbb{L}_{\hat{c}} \partial_xS_{\hat{c}}'(\Theta))v \phi_{\hat{c}}dxdy - \int_{\RTL} (\mathbb{L}_{\hat{c}}\partial_x v)S_{\hat{c}}'(\Theta) \phi_{\hat{c}}dxdy\notag \\
= &- \frac{\pi L(c-c_0)}{c_0} \int_{\R} \Bigl((\mathcal{L}_{c_0}  + \frac{1}{L^2} ) \partial_xQ_{c_0}^{\frac{3}{2}}\Bigr)\Bigl( a_1  v_{1,1} + a_2 v_{1,2} \Bigr)\phi_{c_0} dx \notag \\
&+\frac{\pi L(c-c_0)}{c_0} \int_{\R} \Bigl(\partial_x(\mathcal{L}_{c_0}  + \frac{1}{L^2} )( Q_{c_0}^{\frac{3}{2}}\phi_{c_0})\Bigr)\Bigl( a_1  v_{1,1} + a_2 v_{1,2} \Bigr) dx \notag \\
&+O((|c-c_0|+|\vec{a}|)|c-c_0||\vec{a}|\norm{v}_{L^2}).
\end{align}
From \eqref{Q-eq-t1}--\eqref{Q-eq-t4} we have
\begin{align*}
&\Bigl(\Bigl( \mathcal{L}_{c_0} + \frac{1}{L^2}\Bigr) \partial_x Q_{c_0}^{\frac{3}{2}}\Bigr)\phi_{c_0} 
=- 2c_0 Q_{c_0}^{\frac{5}{2}} + \frac{4}{3} Q_{c_0}^{\frac{7}{2}}.
\end{align*}
By the similar computation we have
\begin{align*}
&\partial_x(\mathcal{L}_{c_0}  + \frac{1}{L^2} )( Q_{c_0}^{\frac{3}{2}}\phi_{c_0})
=-\frac{10c_0}{3}Q_{c_0}^{\frac{5}{2}} + \frac{8}{3} Q_{c_0}^{\frac{7}{2}}.
\end{align*}
Therefore, applying Lemma \ref{lem-Q-pro}, we obtain that 
\begin{align}\label{eq-c-li-term3}
&\Bigl|- \int_{\RTL} (\mathbb{L}_{\hat{c}} \partial_xS_{\hat{c}}'(\Theta))v \phi_{\hat{c}}dxdy - \int_{\RTL} (\mathbb{L}_{\hat{c}}\partial_x v)S_{\hat{c}}'(\Theta) \phi_{\hat{c}}dxdy\Bigr|\notag \\
= & \Bigl| \frac{\pi L(c-c_0)}{c_0} \int_{\R} \Bigl(-\frac{4c_0}{3}Q_{c_0}^{\frac{5}{2}} + \frac{4}{3}Q_{c_0}^{\frac{7}{2}}\Bigr)\Bigl( a_1  v_{1,1} + a_2 v_{1,2} \Bigr)dx\Bigr| \notag \\
&+O((|c-c_0|+|\vec{a}|)|c-c_0||\vec{a}|\norm{v}_{L^2})\notag \\
\leq &\frac{ 3 c_0 \pi L}{4} \int_{\R} (|v_{1,1}|^2+|v_{1,2}|^2) Q_{c_0}dx +  \frac{20|c-c_0|^2|\vec{a}|^2 \pi L}{297c_0^2} \int_{\RTL} Q_{c_0}^5dx \notag \\
&+O((|c-c_0|+|\vec{a}|)|c-c_0||\vec{a}|\norm{v}_{L^2}).
\end{align}

(IV) The estimate of $-\int_{\RTL}R(\eta,\vec{a},c)(v+S_{\hat{c}}'(\Theta))\phi_{\hat{c}}dxdy$.
Since 
\[(\partial_t \mathbb{L}_{\hat{c}})\eta= \dot{\hat{c}} \eta -2 \dot{\hat{c}} (\partial_c \tilde{Q}_{c_0}) \eta \]
and
\[\mathbb{L}_{\hat{c}}\partial_x\eta=2(\partial_x \tilde{Q}_{\hat{c}})\eta +  v_x + \partial_x \eta^2,\]
we have
\begin{align}\label{eq-R-term}
&R(\eta,\vec{a},c)\notag \\
=&-2\eta\partial_x(v+S_{\hat{c}}'(\Theta)) + (\dot{\rho}-\hat{c}) (2(\partial_x \tilde{Q}_{\hat{c}})\eta +  v_x + \partial_x \eta^2)+ (\dot{\rho}-\hat{c})\mathbb{L}_{\hat{c}}\partial_x \Theta -2  (\dot{\rho}-\hat{c})\eta \partial_x(\eta+\Theta) \notag\\
& -2(\partial_x^3 (\tilde{Q}_{\hat{c}}-\Theta))\eta -6 (\partial_x^2 (\tilde{Q}_{\hat{c}}-\Theta))\eta_x + 2(\partial_x (\tilde{Q}_{\hat{c}}-\Theta))v +2(\partial_x (\tilde{Q}_{\hat{c}}-\Theta))\eta^2\notag \\
& -4 (\partial_x (\tilde{Q}_{\hat{c}}-\Theta))\eta_{xx} +4 (\tilde{Q}_{\hat{c}}-\Theta)(\partial_x \tilde{Q}_{\hat{c}})\eta + 2(\tilde{Q}_{\hat{c}}-\Theta)v_x +2(\tilde{Q}_{\hat{c}}-\Theta)\partial_x\eta^2\notag \\
&-4\eta \partial_x((\tilde{Q}_{\hat{c}}-\Theta)\eta)-(\mathbb{L}_{\hat{c}}-2\eta)(  \dot{\vec{a}}\partial_{\vec{a}}\Theta+\dot{c}\partial_c\Theta) +\dot{\hat{c}} \eta -2 \dot{\hat{c}} (\partial_c \tilde{Q}_{c_0}) \eta.
\end{align}
From Lemma \ref{lem-c-mod}, integration by parts, $\mathbb{L}_{\hat{c}}\partial_x \Theta=O(|c-c_0|+|\vec{a}|)$ and $\mathbb{L}_{\hat{c}}\partial_{\vec{a}}\Theta=O(|c-c_0|+|\vec{a}|)$, we obtain
\begin{align}\label{eq-c-li-term4}
-\int_{\RTL}R(\eta,\vec{a},c)(v+S_{\hat{c}}'(\Theta))\phi_{\hat{c}}dxdy
= O((|c-c_0|+|\vec{a}|+\norm{\eta}_{H^1})(|c-c_0|^2|\vec{a}|^2+\norm{v}_{H^1}^2)).
\end{align}

(V) The estimate of $- \int_{\RTL}(\partial_t S_{\hat{c}}'(\Theta))(v+S_{\hat{c}}'(\Theta))\phi_{\hat{c}}dxdy-\frac{1}{2} \int_{\RTL}(v+S_{\hat{c}}'(\Theta))^2 \dot{\hat{c}}\partial_c\phi_{\hat{c}}dxdy$.
Since 
\[\partial_t S_{\hat{c}}'(\Theta))=\frac{\dot{c}}{c_0}\partial_y^2\Theta + \frac{c-c_0}{c_0}\dot{\vec{a}}\cdot \partial_{\vec{a}}\partial_y^2\Theta +\frac{c-c_0}{c_0} \dot{c}\partial_c\partial_y^2\Theta=O((|c-c_0|+|\vec{a}|)\norm{\eta}_{L^2}),\]
from Lemma \ref{lem-c-mod} we have
\begin{align}\label{eq-c-li-term5}
&- \int_{\RTL}(\partial_t S_{\hat{c}}'(\Theta))(v+S_{\hat{c}}'(\Theta))\phi_{\hat{c}}dxdy-\frac{1}{2} \int_{\RTL}(v+S_{\hat{c}}'(\Theta))^2 \dot{\hat{c}}\partial_c\phi_{\hat{c}}dxdy \notag \\
=&O((|c-c_0|+|\vec{a}|)(|c-c_0|^2|\vec{a}|^2+\norm{v}_{L^2}^2)).
\end{align}


Therefore, from (I)--(V) we deduce gathering \eqref{eq-c-li-term1}--\eqref{eq-c-li-term5} that there exists $k_4>0$ such that 
\begin{align}\label{eq-c-li-term6}
 &-\frac{1}{2} \frac{d}{d t} \int_{\RTL}(v+S_{\hat{c}}'(\Theta))^2 \phi_{\hat{c}}dxdy \notag \\
\geq & k_4 \Bigl( \int_{\RTL}v^2 \tilde{Q}_{c_0}dxdy + |c-c_0|^2|\vec{a}|^2 \Bigr)+ O\bigl((|c-c_0|+|\vec{a}|+\norm{\eta}_{H^1})(|c-c_0|^2|\vec{a}|^2+\norm{v}_{H^1}^2)\bigr) .
\end{align}
On the other hand, by  \eqref{eq-c-li}, we have
\begin{align}\label{eq-c-li-9}
 &-\frac{1}{2} \frac{d}{d t} \int_{\RTL}v^2xdxdy \notag \\
 =&\frac{1}{2}\int_{\RTL}(3|\partial_x v|^2 +|\partial_y v|^2+ \hat{c} v^2 ) dxdy  - \int_{\RTL} v^2 \partial_x(x\tilde{Q}_{\hat{c}})dxdy \notag \\
&- \int_{\RTL} (\mathbb{L}_{\hat{c}} \partial_x S_{\hat{c}}'(\Theta))vxdxdy - \int_{\RTL} R(\eta,\vec{a},c) v x dxdy.
 \end{align} 
 From Proposition \ref{cor-high-bdd}, 
\begin{align}\label{eq-c-li-9-0}
\Bigl| \int_{\RTL} \eta v_x v xdxdy \Bigr| \leq \norm{x ^2\eta }_{H^1}^{\frac{1}{2}} \norm{\eta}_{H^1}^{\frac{1}{2}}\norm{v}_{H^1}^2=O(\norm{\eta}_{H^1}^{\frac{1}{2}} \norm{v}_{H^1}^2).
\end{align}
By the  similar calculation to \eqref{eq-c-li-9-0}, we have
\begin{align}\label{eq-c-li-10}
\Bigl| \int_{\RTL} R(\eta,\vec{a},c) v x dxdy\Bigr| =O\bigl( (|c-c_0|+|\vec{a}|+\norm{\eta}_{H^1}^{\frac{1}{2}}) \norm{v}_{H^1}^2\bigr).
\end{align}
By the H\"older inequality and Proposition \ref{cor-high-bdd}, we have 
\begin{align}\label{eq-c-li-11}
\Bigl| \int_{\RTL} v^2 \partial_x(x\tilde{Q}_{\hat{c}})dxdy \Bigr| \leq \Bigl(1 + \hat{c}\norm{x^2 (\partial_x \tilde{Q}_{\hat{c}})^2\tilde{Q}_{\hat{c}}^{-1}}_{L^{\infty}}\Bigr) \int_{\RTL} v^2 \tilde{Q}_{\hat{c}} dxdy + \frac{\hat{c}}{8}\norm{v}_{L^2}^2.
\end{align}
By the H\"older inequality, we obtain there exists $C>0$ such that 
\begin{align}\label{eq-c-li-11-1}
\Bigl| \int_{\RTL} (\mathbb{L}_{\hat{c}} \partial_x S_{\hat{c}}'(\Theta))vxdxdy \Bigr| \leq \frac{\hat{c}}{8}\norm{v}_{L^2}^2 + C|c-c_0||\vec{a}|.
\end{align}
We deduce gathering \eqref{eq-c-li-9}--\eqref{eq-c-li-11-1} that
\begin{align}\label{eq-c-li-term7}
 &-\frac{1}{2} \frac{d}{d t} \int_{\RTL}v^2xdxdy \notag \\
 \geq &\frac{1}{4}\int_{\RTL}(|\nabla v|^2 + \hat{c} v^2 ) dxdy  +\Bigl(1 + \hat{c}\norm{x^2 (\partial_x \tilde{Q}_{\hat{c}})^2\tilde{Q}_{\hat{c}}^{-1}}_{L^{\infty}}\Bigr) \int_{\RTL} v^2 \tilde{Q}_{\hat{c}} dxdy \notag \\
&- C|c-c_0||\vec{a}|+O( (|c-c_0|+|\vec{a}|+\norm{\eta}_{H^1}^{\frac{1}{2}}) \norm{v}_{H^1}^2).
 \end{align} 
From \eqref{eq-c-li-term6} and \eqref{eq-c-li-term7}, we obtain 
\begin{align}\label{vir-est}
&-\frac{1}{2} \frac{d}{dt} \int_{\RTL}\Bigl((v+S_{\hat{c}}'(\Theta))^2\phi_{\hat{c}} + \varepsilon_+  xv^2\Bigr)dxdy \notag \\
 \geq & \frac{\varepsilon _+ }{4} \int_{\RTL}(|\nabla v|^2+ c_0 v^2)dxdy + \frac{k_4}{2} |c-c_0|^2|\vec{a}|^2 \notag \\
&+ O( (|c-c_0|+|\vec{a}|+\norm{\eta}_{H^1}^{\frac{1}{2}}) (\norm{v}_{H^1}^2+|c-c_0|^2|\vec{a}|^2)),
\end{align}
where 
\[ \varepsilon _+ = \frac{k_4}{2}\Bigl(1+C+ c_0\norm{x^2 (\partial_x \tilde{Q}_{c_0})^2\tilde{Q}_{c_0}^{-1}}_{L^{\infty}}\Bigr)^{-1}>0.\]
Integrating \eqref{vir-est} between $t_1$ and $t_2$, we have for sufficiently small $\varepsilon _0>0$ 
\begin{align}\label{eq-c-li-12}
&\int_{\RTL} \Bigl(\Bigl(v(t_1)+S_{\hat{c}(t_1)} ' \bigl(\Theta(\vec{a}(t_1),c(t_1))\bigr)  \Bigr)^2\phi_{\hat{c}(t_1)}- \Bigl(v(t_2)+S_{\hat{c}(t_2)} ' \bigl(\Theta(\vec{a}(t_2),c(t_2))\bigr)  \Bigr)^2 \phi_{\hat{c}(t_2)}\notag \\
 &+ xv(t_1)^2 -x v(t_2)^2\Bigr)dxdy \notag \\
 \geq & \int_{t_1}^{t_2} \Bigl(  \frac{\varepsilon _+}{4} \int_{\RTL} (|\nabla v(t)|^2 + c_0 v(t)^2)dxdy + k_4 |c(t)-c_0|^2 |\vec{a}(t)|^2 \Bigr)dt.
 \end{align}
From Proposition \ref{cor-high-bdd} 
\begin{align*}
&\int_{-\infty}^{\infty} \Bigl(  \frac{\varepsilon _+}{4} \int_{\RTL} (|\nabla v(t)|^2 + c_0 v(t)^2)dxdy + k_4 |c(t)-c_0|^2 |\vec{a}(t)|^2 \Bigr)dt \\
\lesssim& \sup_{t \in \R} \Bigl|\int_{\RTL}\Bigl((v+S_{\hat{c}}'(\Theta))^2\phi_{\hat{c}} + \varepsilon_+  xv^2\Bigr)dxdy \Bigr| < \infty.
\end{align*}
Therefore, there exist sequences $\{t_{1,n}\}_n$ and $\{t_{2,n}\}_n$ such that 
\[\lim_{n \to \infty} t_{1,n}= -\infty, \quad \lim_{n \to \infty} t_{2,n} = \infty\]
 and
\begin{align}\label{eq-c-li-13}
& \lim_{n \to \infty}\Bigl|  \frac{\varepsilon _+}{4} \int_{\RTL} (|\nabla v(t_{1,n})|^2 + c_0 v(t_{1,n})^2)dxdy + k_4 |c(t_{1,n})-c_0|^2 |\vec{a}(t_{1,n})|^2 \Bigr|\notag \\
& = \lim_{n \to \infty}\Bigl|  \frac{\varepsilon _+}{4} \int_{\RTL} (|\nabla v(t_{2,n})|^2 + c_0 v(t_{2,n})^2)dxdy + k_4 |c(t_{2,n})-c_0|^2 |\vec{a}(t_{2,n})|^2 \Bigr|=0.
\end{align}
Combining \eqref{eq-c-li-12} and \eqref{eq-c-li-13}, we obtain that 
\begin{align*}
 \int_{-\infty}^{\infty} \Bigl(  \frac{\varepsilon _+}{4} \int_{\RTL} (|\nabla v(t)|^2 + c_0 v(t)^2)dxdy + k_4 |c(t)-c_0|^2 |\vec{a}(t)|^2 \Bigr)dt=0
 \end{align*}
which implies $v\equiv 0$ and $|c-c_0||\vec{a}|\equiv 0$.
By \eqref{con-coer} and $v \equiv 0$, we have $\eta\equiv 0$.
Therefore, we obtain the conclusion.

\qed

\subsection{Non-critical case $L < \frac{2}{\sqrt{5c_0}}$}
In this subsection, we show the Liouville property for  $L < \frac{2}{\sqrt{5c_0}}$.
Since the proof of the Liouville property for  $L < \frac{2}{\sqrt{5c_0}}$ is similar to the proof of the Liouville property for  $L = \frac{2}{\sqrt{5c_0}}$, we omit the detail of the proof.
\begin{lemma}\label{lem-nc-mod}
Let $c_0>0$.
There exist $\varepsilon _0,K>0$ such that for any $0<\varepsilon <\varepsilon _0$ the following is true.
For any solution $u \in C(\R,H^1(\RTL))$ of $\eqref{ZKeq}$ satisfying 
\begin{align*}
\inf_{b \in \R} \norm{u(t,\cdot,\cdot)-Q_{c_0}(\cdot-b,\cdot)}_{H^1}\leq \varepsilon 
\end{align*}
there exist  $\rho_1, c \in C^1(\R,\R)$  uniquely such that 
\begin{align*}
\eta(t,x,y)=u(t,x+\rho(t),y)-Q_{c(t)}(x)
\end{align*}
satisfies for all $t\in \R$ 
\begin{align*}
|c(t)-c_0|+&\norm{\eta(t)}_{H^1}\leq K_0\varepsilon,\\
\int_{\RTL} \eta(t)\partial_xQ_{c(t)}dxdy=&\int_{\RTL}\eta(t)Q_{c(t)}dxdy=0
\end{align*} 
and
\begin{align*}
|\dot{c}(t)|^{\frac{1}{2}}+|\dot{\rho}(t)-c(t)|\leq K_0\norm{\eta(t)}_{L^2}.
\end{align*}
\end{lemma}
The following is Liouville property in the non-critical case.
\begin{theorem}\label{thm-nc-Liouville}
Let $c_0>0$ and $L<\frac{2}{\sqrt{5c_0}}$.
There exists $\varepsilon  _0>0$ satisfies the following.
For any solution $u \in C(\R,H^1(\RTL))$ to $\eqref{ZKeq}$ satisfying $\eqref{eq-6-1}$ and $\eqref{ass-decay}$, there exist $c_+>0$ and $\rho_0\in \R $ such that
\[u(t,x,y)=Q_{c_+}(x-c_+t+\rho_0,y).\]
\end{theorem}
\begin{remark}
The proof of Theorem \ref{thm-nc-Liouville} is easier then the proof of Theorem \ref{thm-c-Liouville}. 
In the case $L < \frac{2}{\sqrt{5c_0}}$, $\mathbb{L}_{c_0}$ has the following coercive type estimate.
There exists $k_{5}>0$ such that for $\eta \in H^1(\RTL)$ with $(\eta, \partial_x\tilde{Q}_{c_0})_{L^2}=(\eta,\tilde{Q}_{c_0})_{L^2}=0$, 
\[ \tbr{\mathbb{L}_{c_0}\eta,\eta}_{H^{-1},H^1} \geq k_5 \norm{\eta}_{H^1}^2.\]
Therefore, from Lemma \ref{lem-coer} we can show a coercive type estimate for the virial identity 
\[-\frac{1}{2} \frac{d}{dt} \int_{\RTL}v^2(\phi_{c} + \varepsilon_+  x)dxdy
 \geq  \frac{\varepsilon_+  }{4} \int_{\RTL}(|\nabla v|^2+ c_0 v^2)dxdy 
+ o(\norm{v}_{H^1}^2) \mbox{ as } \varepsilon _0 \to 0,\]
for sufficiently small $\varepsilon_+ >0$.
\end{remark}

\section{Asymptotic stability}
In this section, we prove Theorem \ref{asymptotic-stability} by applying the monotonicity property and the Liouville property in Section 6.
We follow the argument Martel and Merle \cite{M M 1, M M 2, M M 3} for the generalized KdV equation and C\^ote et al. \cite{C M P S} for the Zakharov--Kuznestov equation on $\R^2$.

\subsection{Critical case $L = \frac{2}{\sqrt{5c_0}}$}

In this subsection. we consider the critical case $L = \frac{2}{\sqrt{5c_0}}$.
The following proposition shows the compactness of the orbit of solutions in $H^1(x>-A)$.
\begin{proposition}\label{prop-c-cpt}
Let $c_0>0$ and $L = \frac{2}{\sqrt{5c_0}}$.
There exists $ 0<\varepsilon _*<\varepsilon _0$ such that if $0<\varepsilon \leq \varepsilon _*$ and $u \in C(\R,H^1(\RTL))$ is a solution to $\eqref{ZKeq}$ satisfying $\sup_{t \in \R} {\rm dist}_{c_0} (u(t))< \varepsilon$ then the following holds true.
For any sequence $\{t_n\}_n$ with $\lim_{n \to \infty }t_n=\infty$, there exists a subsequence $\{t_{n_k}\}_k$ and $\tilde{u}_0 \in H^1(\RTL)$ such that 
\[ \norm{u(t_{n_k},\cdot+\rho(t_{n_k}),\cdot) - \tilde{u}_0}_{  H^1(x>-A)} \to 0 \mbox{ as } k \to \infty\]
for any $A>0$, where $\rho$ is the function associated to $u$ given by Lemma \ref{lem-c-mod}.
Moreover, the solution $\tilde{u}$ of $\eqref{ZKeq}$ with $\tilde{u}(0)=\tilde{u}_0$ satisfies
\begin{align}\label{eq-7-cond-1}
 \norm{\tilde{u}(t,\cdot + \tilde{\rho}(t),\cdot)-\tilde{Q}_{c_0}}_{H^1} \lesssim \varepsilon _0, \quad t \in \R
\end{align}
and
\begin{align}\label{eq-7-cond-2}
 \int_{\T_L}|\tilde{u}(t,x+\tilde{\rho}(t),y)|^2 dy \lesssim e^{-\delta_1|x|}, \quad (t,x) \in \R^2
\end{align}
for some $\delta_1>0$, where $\tilde{\rho}$ is the function associated to $\tilde{u}$ given by Lemma \ref{lem-c-mod} and $\tilde{\rho}(0)=0$.
\end{proposition}
Since the proof of Proposition \ref{prop-c-cpt} is similar to the proof of Proposition 4.1 in \cite{C M P S}, we omit the proof.

Next, we show (ii) of Theorem \ref{asymptotic-stability}.

\proof[Proof of {\rm (ii)} of Theorem \ref{asymptotic-stability}]

Let $\beta>0$ and $u$ be a solution to \eqref{ZKeq} with ${\rm dist}_{c_0}(u(0))<\varepsilon $.
By Theorem \ref{orbital-stability} if $\varepsilon $ is sufficiently small, then ${\rm dist}_{c_0}(u(t))<\varepsilon _*$.
Let $\rho$, $c$ and $\vec{a}$ be functions associated to $u$ given by Lemma \ref{lem-c-mod}.
From Proposition \ref{prop-c-cpt}, for any sequence $\{t_n\}_n$ with $\lim_{n \to \infty}t_n =\infty$, there exist a subsequence $\{t_{n_k}\}_k$, $\tilde{c}_0>0$, $\tilde{\vec{a}}_0 \in \R^2$ and $\tilde{u}_0 \in H^1(\RTL)$ such that for $A >0$
\begin{align*}
u(t_{n_k},\cdot+\rho(t_{n_k}),\cdot) \underset{ k \to \infty}{ \rightarrow } \tilde{u}_0 \mbox{ in } H^1(x > -A), \quad 
c(t_{n_k}) \underset{ k \to \infty}{ \rightarrow } \tilde{c}_0 \mbox{ and } \vec{a}(t_{n_k}) \underset{ k \to \infty}{ \rightarrow } \tilde{\vec{a}}_0.
\end{align*}
Moreover, the solution $\tilde{u}$ of \eqref{ZKeq} with $\tilde{u}(0)=\tilde{u}_0$ satisfies \eqref{eq-7-cond-1} and \eqref{eq-7-cond-2}.
Let $\tilde{\rho}$, $\tilde{c}$ and $\tilde{\vec{a}}$ be functions associated to $\tilde{u}$ given by Lemma \ref{lem-c-mod}.
By uniqueness of the decomposition in Lemma \ref{lem-c-mod}, we have $\tilde{\rho}(0) =0$, $\tilde{c}(0)=\tilde{c}_0$ and $\tilde{\vec{a}}(0)=\tilde{\vec{a}}_0$.
Applying Theorem \ref{thm-c-Liouville}, we obtain that there exist $\rho_0\in \R$, $c_1 >0$ and $\vec{a}_1 \in \R^2$ such that $|c_1-c_0||\vec{a}_1|=0$ and 
\[\tilde{u}(t,x,y)=\Theta(\vec{a}_1,c_1)(x-\hat{c}_1t-\rho_0,y),\]
where
\[ \hat{c}_1=
\begin{cases}
c_1, & \vec{a}_1=(0,0),\\
\check{c}(\vec{a}_1), & c_1=c_0.
\end{cases}\]
By uniqueness of the decomposition in Lemma \ref{lem-c-mod}, $\rho_0=0$, $c_1=\tilde{c}_0$, $\vec{a}_1 =\tilde{\vec{a}}_0$ and $|\tilde{c}_0-c_0||\tilde{\vec{a}}_0|=0$.
Since for any sequence $\{t_n\}_n$ with $\lim_{n \to \infty}t_n =\infty$ there exists a subsequence $\{t_{n_k}\}_k$ such that
\[ u(t_{n_k},\cdot+\rho(t_{n_k}),\cdot) - \Theta(\vec{a}(t_{n_k}),c(t_{n_k}))  \underset{ k \to \infty}{ \rightarrow } 0  \mbox{ in } H^1(x > -A) \mbox{ and } |c(t_{n_k})-c_0||\vec{a}(t_{n_k})| \underset{ k \to \infty}{ \rightarrow } 0 ,\]
we obtain
\begin{align}\label{eq-7-c-1}
u(t,\cdot+\rho(t),\cdot) - \Theta(\vec{a}(t),c(t))  \underset{ t \to \infty}{ \rightarrow } 0  \mbox{ in } H^1(x > -A)
\end{align}
and
\begin{align}\label{eq-7-c-1-0}
|c(t)-c_0||\vec{a}(t)| \underset{ t \to \infty}{ \rightarrow } 0.
\end{align}
Moreover, \eqref{eq-7-c-1} implies that for $R>0$ and $x_0 \in \R$
\begin{align}\label{eq-7-c-1-1}
\lim_{t \to \infty} \int_{\RTL}(|\nabla \eta|^2+|\eta|^2)(t,x-\rho(t),y) \psi_R(x-\rho(t)+x_0)dxdy =0,
\end{align}
where $\eta(t,x,y)=u(t, x+\rho(t),y)-\Theta(\vec{a}(t),c(t))$.
By \eqref{eq-7-c-1}, for any $\alpha>0$ and $R > 2/\sqrt{\beta}$ there exist $x_1 \in \R$ and $T_1>0$ such that for $x_0 > x_1$ and $t >T_1$
\begin{align}\label{eq-7-c-2}
\Bigl| \int_{\RTL} |u(t,x,y)|^2\psi_R(x-\rho(t)+x_0)dxdy-\int_{\RTL}|\Theta(\vec{a}(t),c(t))|^2dxdy\Bigr| < \alpha,
\end{align}
where $\psi_R$ is defined by \eqref{def-psi}.
From Lemma \ref{lem-L2-mono-} there exists $x_2 \in \R$ such that for $x_0 \geq x_2$ and $t\geq t'$
\begin{align}\label{eq-7-c-3}
\int_{\RTL}|u(t,x,y)|^2&\psi_R(x-\rho(t)+x_0)dxdy  \notag \\
&- \int_{ \RTL}|u(t',x,y)|^2\psi_R(x-\rho(t')+x_0)dxdy \leq \alpha.
\end{align}
By \eqref{eq-7-c-2} and \eqref{eq-7-c-3} we have that for $t\geq t' >T_0$
\[\int_{\RTL}|\Theta(\vec{a}(t),c(t))|^2dxdy \leq \int_{\RTL} |\Theta(\vec{a}(t'),c(t'))|^2 dxdy +3 \alpha.\]
Since for any $\alpha>0$
\[\limsup_{t \to \infty} \int_{\RTL}|\Theta(\vec{a}(t),c(t))|^2dxdy \leq \liminf_{t' \to \infty} \int_{\RTL} |\Theta(\vec{a}(t'),c(t'))|^2 dxdy +3 \alpha,\]
$\int_{\RTL}|\Theta(\vec{a}(t),c(t))|^2dxdy$ has the limit as $t \to \infty$.
By the definition of $\Theta$, we have
\begin{align}\label{eq-7-c-4} \lim_{ t \to \infty} \int_{\RTL}|\Theta(\vec{a}(t),c(t))|^2dxdy = \lim_{t \to \infty} \Bigl(\frac{c(t)}{c_0}\Bigr)^{\frac{3}{2}}\norm{\varphi_{c_0}(\vec{a}(t))}_{L^2}^2.
\end{align}
Since $\norm{\varphi_{c_0}(\vec{a}(t))}_{L^2}^2$ is strictly increasing with respect to $|\vec{a}|$,
from \eqref{eq-7-c-1-0} and \eqref{eq-7-c-4} the $\omega$-limit set of $(|\vec{a}(t)|,c(t))$
is at most two points.
By the continuity of $\vec{a}(t)$ and $c(t)$, the $\omega$-limit set of $(|\vec{a}(t)|,c(t))$ is the one point set which implies there exist $a_+\geq 0$ and $c_+>0$ such that
\begin{align}\label{eq-7-c-5}
\lim_{t \to \infty} |\vec{a}(t)| = a_+, \quad \lim_{t \to \infty} c(t)=c_+.
\end{align}
Therefore, by Corollary \ref{cor-initial-c} we have 
\[|c_+-c_0|+|a_+|^2 \lesssim \norm{u(0)-\tilde{Q}_{c_0}}_{H^1}.\]

Next, we improve the convergence of \eqref{eq-7-c-6}.
By Lemma \ref{lem-L2-mono}, for all $t_1\leq t_2$, $x_0>0$ and $R>2/\sqrt{\beta}$ 
\begin{align}\label{eq-7-c-7}
\int_{\RTL}|u(t_2,x,y)|^2\psi_R(\tilde{x}(t_1,t_2)) dxdy - \int_{\RTL} |u(t_1,x,y)|^2 \psi_R(\tilde{x}(t_1,t_1))dxdy\leq C e^{-x_0/R},
\end{align}
where $\tilde{x}(t,\tau)=x-\rho(t)-\frac{\beta}{2}(\tau-t)+x_0$.
By \eqref{eq-6-c-3} if $R > \frac{1}{4c_0} $ and $\varepsilon _0$ is sufficiently small, then we have
\begin{align}\label{eq-7-c-8}
&\Bigl| \int_{\RTL} \eta(t,x,y)\Theta(\vec{a}(t),c(t))\psi_R(x+x_0)dxdy\Bigr| \notag \\
=& \Bigl| \int_{\RTL} \eta(t,x,y)\Theta(\vec{a}(t),c(t))(1-\psi_R(x+x_0))dxdy\Bigr| \notag \\
 =& \norm{\eta(t)}_{L^2} \norm{\Theta(\vec{a}(t),c(t))(1-\psi_R(x+x_0))}_{L^2}\lesssim e^{-x_0/R}.
\end{align}
Since
\begin{align*}
 \bigl(\eta(t,x-\rho(t),y)\bigr)^2=&\bigl(u(t,x,y)\bigr)^2 -2 \eta(t,x-\rho(t),y)\Theta(\vec{a}(t),c(t))(x-\rho(t),y)\notag \\
&-\bigl( \Theta(\vec{a}(t),c(t))(x-\rho(t),y)\bigr)^2,
\end{align*}
from \eqref{eq-7-c-7} and \eqref{eq-7-c-8} we have that there exists $C_0>0$ such that
\begin{align*}
&\int_{\RTL} \bigl(\eta(t_2,x-\rho(t_2),y)\bigr)^2 \psi_R(\tilde{x}(t_1,t_2))dxdy  \notag \\
&-  \int_{\RTL}\bigl(\eta(t_1,x-\rho(t_1),y)\bigr)^2 \psi_R(\tilde{x}(t_1,t_1))dxdy \notag \\
\leq & C_0(e^{-x_0/R} + |c(t_1)-c(t_2)|+||\vec{a}(t_1)|^2-|\vec{a}(t_2)|^2|).
\end{align*}
For $t>0$ large enough, we define $0<t'<t$ such that $\rho(t')+\frac{\beta}{2}(t-t') -x_0 =\beta t$.
Then, we have $t' \to \infty$ as $ t \to \infty$.
\begin{align*}
&\int_{\RTL} \bigl( \eta(t,x-\rho(t),y)\bigr)^2 \psi_R(x-\beta t) dxdy  \\
\leq & \int_{\RTL} \bigl( \eta(t',x-\rho(t'),y)\bigr)^2 \psi_R(x-\rho(t')+x_0) dxdy \notag \\
&+ C_0(e^{-x_0/R} + |c(t_1)-c(t_2)|+||\vec{a}(t_1)|^2-|\vec{a}(t_2)|^2|).
\end{align*}
From \eqref{eq-7-c-1-1} and \eqref{eq-7-c-5}, we obtain for any $x_0>0$
\[ \limsup_{ t \to \infty} \int_{\RTL} \bigl( \eta(t,x-\rho(t),y)\bigr)^2 \psi_R(x-\beta t) dxdy \leq C_0e^{-x_0/R}.\]
Therefore, 
\begin{align}\label{eq-7-c-9}
\lim_{t\to \infty} \int_{\RTL} \bigl( \eta(t,x-\rho(t),y)\bigr)^2 \psi_R(x-\beta t) dxdy =0.
\end{align}
From Lemma \ref{lem-E-mono} we have 
for all $t_1\leq t_2$, $x_0>0$ and $R>2/\sqrt{\beta}$ 
\begin{align}\label{eq-7-c-10}
J_{x_0,t_1}(u(t_2)) - J_{x_0,t_1}(u(t_1))\leq C e^{-x_0/R}.
\end{align}
Moreover, we have
\begin{align}\label{eq-7-c-11}
&\Bigl |\int_{\RTL} \bigl(u(t_2,x,y)\bigr)^3 \psi_R(\tilde{x}(t_1,t_2))dxdy -\int_{\RTL}\bigl(u(t_1,x,y)\bigr)^3 \psi_R(\tilde{x}(t_1,t_1))dxdy \Bigr| \notag \\
\lesssim &  \Bigl( \int_{\RTL} \bigl(\eta(t_1,x-\rho(t_1),y) \bigr)^2 \psi_R(\tilde{x}(t_1,t_2))dxdy\Bigr)^{\frac{1}{2}} \notag \\
&+  \Bigl( \int_{\RTL} \bigl(\eta(t_2,x-\rho(t_2),y)\bigr)^2 \psi_R(\tilde{x}(t_1,t_2))dxdy\Bigr)^{\frac{1}{2}}\notag \\
&+  (e^{-x_0/R} + |c(t_1)-c(t_2)|+||\vec{a}(t_1)|-|\vec{a}(t_2)||).
\end{align}
By \eqref{eq-7-c-10} and \eqref{eq-7-c-11}
\begin{align}\label{eq-7-c-12}
&\int_{\RTL} \bigl|\nabla \eta(t_2,x-\rho(t_2),y)\bigr|^2 \psi_R(\tilde{x}(t_1,t_2))dxdy  \notag \\
&-  \int_{\RTL}\bigl|\nabla \eta(t_1,x-\rho(t_1),y)\bigr|^2 \psi_R(\tilde{x}(t_1,t_1))dxdy \notag \\
\lesssim &  \Bigl( \int_{\RTL} \bigl(\eta(t_1,x-\rho(t_1),y) \bigr)^2 \psi_R(\tilde{x}(t_1,t_2))dxdy\Bigr)^{\frac{1}{2}} \notag \\
&+  \Bigl( \int_{\RTL} \bigl(\eta(t_2,x-\rho(t_2),y)\bigr)^2 \psi_R(\tilde{x}(t_1,t_2))dxdy\Bigr)^{\frac{1}{2}}\notag \\
&+  (e^{-x_0/R} + |c(t_1)-c(t_2)|+||\vec{a}(t_1)|-|\vec{a}(t_2)||).
\end{align}
From \eqref{eq-7-c-12} with $t_1=t' $ and $t_2=t$ we obtain that there exists $C>0$ such that 
\begin{align*}
&\int_{\RTL} \bigl|\nabla \eta(t,x-\rho(t),y)\bigr|^2 \psi_R(x-\beta t)dxdy \notag \\
\leq & \int_{\RTL} \bigl|\nabla \eta(t',x-\rho(t'),y)\bigr|^2 \psi_R(x-\rho(t')+x_0)dxdy \notag \\
&+C \Bigl( \int_{\RTL} \bigl(\eta(t,x-\rho(t),y) \bigr)^2 \psi_R(x-\beta t)dxdy\Bigr)^{\frac{1}{2}} \notag \\
&+ C \Bigl( \int_{\RTL} \bigl(\eta(t',x-\rho(t'),y)\bigr)^2 \psi_R(x-\rho(t')+x_0)dxdy\Bigr)^{\frac{1}{2}}\notag \\
&+ C (e^{-x_0/R} + |c(t')-c(t)|+||\vec{a}(t')|-|\vec{a}(t)||).
\end{align*}
Therefore, it follows form \eqref{eq-7-c-1-1}, \eqref{eq-7-c-5} and \eqref{eq-7-c-9} that for $x_0>0$
\[\limsup_{t \to \infty} \int_{\RTL} \bigl|\nabla \eta(t,x-\rho(t),y)\bigr|^2 \psi_R(x-\beta t)dxdy \leq C e^{-x_0/R}\]
which implies
\begin{align}\label{eq-7-c-13}
\lim_{t \to \infty} \int_{\RTL} \bigl|\nabla \eta(t,x-\rho(t),y)\bigr|^2 \psi_R(x-\beta t)dxdy=0.
\end{align}
Then, we define $\rho_2(t)$ by 
\[\rho_2(t)=
\begin{cases}
\Phi^{-1}\Bigl ( \frac{\vec{a}(t)}{|\vec{a}(t)|}\Bigr), & \mbox{if } |\vec{a}(t)|\neq 0 \mbox{ and } a_+\neq 0,\\
0, & \mbox{if otherwise},
\end{cases}
\]
where $\Phi (\theta)=(\cos \theta,-\sin \theta)$ for $\theta \in \R/2\pi \Z$.
Using
\[ \Theta(\vec{a}(t),c(t))(x,y)=\Theta((|\vec{a}(t)|,0),c(t))(x,y-\rho_2(t)),\]
from \eqref{eq-7-c-9} and \eqref{eq-7-c-13} we obtain 
\begin{align}\label{eq-7-c-6}
u(t,\cdot+\rho(t),y+\rho_2(t))-\Theta((a_+,0),c_+) \underset{ t \to \infty}{ \rightarrow } 0  \mbox{ in } H^1(x > \beta t).
\end{align}
From  \eqref{eq-6-c-2-1}--\eqref{eq-6-c-2-3}, \eqref{eq-7-c-1} and \eqref{eq-7-c-5}, we have
\[\lim_{t \to \infty} \dot{c}(t) = \lim_{t \to \infty} |\dot{\vec{a}}(t)|=\lim_{t \to \infty}|\dot{\rho}(t)-\hat{c}_+|=0,\]
where $\hat{c}_+$ is defined by \eqref{def-hatc}.
If $a_+=0$, then $\dot{\rho}_2(t)=0$ for $t>0$.
On the other hand, if $a_+ \neq 0$, then $|\dot{\rho}_2(t)| \lesssim |\dot{\vec{a}}(t)| \to 0$ as $t \to \infty$.

\qed

\subsection{Non-Critical case $L < \frac{2}{\sqrt{5c_0}}$}
In this subsection, we show (i) of Theorem \ref{asymptotic-stability}.
The proof of (i) of Theorem \ref{asymptotic-stability} is similar to the proof of (ii) of Theorem \ref{asymptotic-stability}.
We omit the detail of the proof.

\begin{proposition}\label{prop-nc-cpt}
Let $c_0>0$ and $L < \frac{2}{\sqrt{5c_0}}$.
There exists $ 0<\varepsilon _*<\varepsilon _0$ such that if $0<\varepsilon \leq \varepsilon _*$ and $ u \in C(\R,H^1(\RTL))$ is a solution to $\eqref{ZKeq}$ satisfying $\sup_{t \in \R} {\rm dist}_{c_0} (u(t))< \varepsilon$ then the following holds true.
For any sequence $\{t_n\}_n$ with $\lim_{n \to \infty }t_n=\infty$, there exists a subsequence $\{t_{n_k}\}_k$ and $\tilde{u}_0 \in H^1(\RTL)$ such that 
\[ u(t_{n_k},\cdot+\rho(t_{n_k}),\cdot) \to \tilde{u}_0 \mbox{ in } H^1(x>-A) \mbox{ as } k \to \infty\]
for any $A>0$, where $\rho$ is the function associated to $u$ given by Lemma \ref{lem-nc-mod}.
Moreover, the solution $\tilde{u}$ of $\eqref{ZKeq}$ with $\tilde{u}(0)=\tilde{u}_0$ satisfies
\begin{align}\label{eq-7-cond-1-n}
 \norm{\tilde{u}(t,\cdot + \tilde{\rho}(t),\cdot)-\tilde{Q}_{c_0}}_{H^1} \lesssim \varepsilon _0, \quad t \in \R
\end{align}
and
\begin{align}\label{eq-7-cond-2-n}
 \int_{\T_L}|\tilde{u}(t,x+\tilde{\rho}(t),y)|^2 dy \lesssim e^{-\delta_1|x|}, \quad (t,x) \in \R^2
\end{align}
for some $\delta_1>0$, where $\tilde{\rho}$ is the function associated to $\tilde{u}$ given by Lemma \ref{lem-nc-mod} and $\tilde{\rho}(0)=0$.
\end{proposition}
By applying Theorem \ref{thm-nc-Liouville} and Proposition \ref{prop-nc-cpt} and the similar proof to the proof of (ii) of Theorem \ref{asymptotic-stability}, we obtain (i) of Theorem.

\section*{Acknowledgments}
The author would like to express his great appreciation to Professor Yoshio Tsutsumi for a lot to helpful advices and encouragements.
The author would like to thank Professor Tetsu Mizumachi for his helpful advices.

Yohei Yamazaki

Department of Mathematics

Kyoto University

Kyoto 606-8502

Japan

E-mail address: y-youhei@math.kyoto-u.ac.jp

\begin{thebibliography}{99}




\bibitem{A P S}
	\newblock J. C. Alexander, R. L. Pego and R. L. Sachs,
 	\newblock \emph{On the transverse instability of solitary waves in the Kadomtsev-Petviashvili equation}, 
 	\newblock Phys. Lett. A  (1997), no. 3-4, 187--192. 
 	
\bibitem{TBB}
	\newblock T. B. Benjamin, 
 	\newblock \emph{The stability of solitary waves}, 
 	\newblock Proc. Roy. Soc. (London) Ser. A  \textbf{328} (1972), 153--183. 

\bibitem{JTB}
	\newblock T.J. Bridges, 
 	\newblock \emph{Universal geometric conditions for the transverse instability of solitary waves}, 
 	\newblock Phys. Rev. Lett.  \textbf{84} no. 12 (2000) 2614--2617.
 	
\bibitem{C G N T}
	\newblock S.-M. Chang, S. Gustafson, K. Nakanishi and T.-P. Tsai,
 	\newblock \emph{Spectra of linearized operators for NLS solitary waves}, 
 	\newblock SIAM J. Math. Anal.  \textbf{39} (2007/08), no. 4, 1070--1111.  

\bibitem{C M P S}
	\newblock R. C\^ote, C. Mu\~noz, D. Pilod and G. Simpson, 
 	\newblock \emph{Asymptotic Stability of High-dimensional Zakharov-Kuznetsov Solitons}, 
 	\newblock  Arch. Ration. Mech. Anal. \textbf{220} (2016), no. 2, 639--710.

\bibitem{C P}
	\newblock A. Comech and D. E. Pelinovsky, 
 	\newblock \emph{Purely nonlinear instability of standing waves with minimal energy}, 
 	\newblock   Comm. Pure Appl. Math. \textbf{56} (2003), no. 11, 1565--1607.

\bibitem{D P C}
	\newblock B. Deconinck, D. E. Pelinovsky and J. D. Carter,
 	\newblock \emph{Transverse instabilities of deep-water solitary waves}, 
 	\newblock Proc. R. Soc. Lond. Ser. A Math. Phys. Eng. Sci. \textbf{462} (2006), no.2071, 2671--2694.
 	



\bibitem{AdB}
	\newblock A. de Bouard, 
 	\newblock \emph{Stability and instability of some nonlinear dispersive solitary waves in higher dimension}, 
 	\newblock  Proc. Roy. Soc. Edinburgh Sect. A \textbf{126} (1996), no. 1, 89--112.
 	
\bibitem{AVF}
	\newblock A. V. Faminskii, 
 	\newblock \emph{The Cauchy problem for the Zakharov-Kuznetsov equation}, 
 	\newblock  translation in Differential Equations A \textbf{31} (1995), no. 6, 1002--1012 

\bibitem{G O}
	\newblock V. Georgiev and M. Ohta, 
	\newblock \emph{Nonlinear instability of linearly unstable standing waves for nonlinear Schr\"odinger equations},
	\newblock Journal of the Mathematical Society of Japan \textbf{64} (2012), no. 2, 533--548.
	



\bibitem{EG}
	\newblock E. Grenier, 
 	\newblock \emph{On the nonlinear instability of Euler and Prandtl equations}, 
 	\newblock Comm. Pure Appl. Math. \textbf{53} (2000) 1067--1091.


\bibitem{G S S 1}
     \newblock M. Grillakis, J. Shatah and W. Strauss,
     \newblock \emph{Stability theory of solitary waves in the presence of symmetry. I},
     \newblock J. Funct. Anal. A \textbf{74} (1987), no. 1, 160--197.

\bibitem{G S S 2}
	\newblock M. Grillakis, J. Shatah and W. Strauss, 
	\newblock \emph{Stability theory of solitary waves in the presence of symmetry II}, 
	\newblock J. Funct. Anal. \textbf{94} (1990), no. 2, 308--348.

\bibitem{AG1}
	\newblock A. Gr\"unrock,
	\newblock \emph{A remark on the modified Zakharov-Kuznetsov equation in three space dimensions}, 
	\newblock  Math. Res. Lett. \textbf{21} (2014), no. 1, 127--131.
	
\bibitem{MJ}
	\newblock K. Johnson, 
	\newblock \emph{The transverse instability of periodic waves in Zakharov-Kuznetsov type equations}, 
	\newblock Stud. Appl. Math. \textbf{124} (2010), no. 4, 323--345.
	
\bibitem{G H}
	\newblock A. Gr\"unrock and S. Herr, 
	\newblock \emph{The Fourier restriction norm method for the Zakharov-Kuznetsov equation}, 
	\newblock Discrete Contin. Dyn. Syst. \textbf{34} (2014), no. 5, 2061--2068. 

\bibitem{K P} 
     \newblock B. B. Kadomtsev and V. I. Petviashvili, 
     \newblock \emph{On the stability of solitary waves in weakly dispersive media},
     \newblock Sov. Phys. Dokl. \textbf{15} (1970), 539--541.


\bibitem{TK1}
     \newblock T. Kato, 
     \newblock \emph{On the Cauchy problem for the (generalized) Korteweg-de Vries equation. Studies in applied mathematics},
     \newblock  Adv. Math. Suppl. Stud.,  \textbf{8} 93--128, Academic Press, New York, 1983.

\bibitem{K K P}
     \newblock E. Kirr, P. G. Kevrekidis and D. E. Pelinovsky, 
     \newblock \emph{Symmetry-breaking bifurcation in the nonlinear Schr\"odinger equation with symmetric potentials},
     \newblock Comm. Math. Phys. \textbf{308} (2011), no. 3, 795--844.

\bibitem{L L S}
     \newblock D. Lannes, F. Linares and J.-C. Saut, 
     \newblock \emph{The Cauchy problem for the Euler-Poisson system and derivation of the Zakharov-Kuznetsov equation},
     \newblock Studies in phase space analysis with applications to PDEs, 181--213, Progr. Nonlinear Differential Equations Appl., 84, Birkhauser/Springer, New York, 2013. 

\bibitem{L P 1}
     \newblock  F. Linares and A. Pastor, 
     \newblock \emph{Well-posedness for the two-dimensional modified Zakharov-Kuznetsov equation},
     \newblock   SIAM J. Math. Anal. \textbf{41} (2009), no. 4, 1323--1339.

\bibitem{L P 2}
     \newblock  F. Linares and A. Pastor, 
     \newblock \emph{Local and global well-posedness for the 2D generalized Zakharov-Kuznetsov equation},
     \newblock  J. Funct. Anal. \textbf{206} (2011), no. 4, 1060--1085.

\bibitem{L S}
     \newblock F. Linares and J.-C. Saut, 
     \newblock \emph{The Cauchy problem for the 3D Zakharov-Kuznetsov equation},
     \newblock Discrete Contin. Dyn. Syst. \textbf{24} (2009), no. 2, 547--565. 
     
\bibitem{MM}
     \newblock M. Maeda, 
     \newblock \emph{Stability of bound states of Hamiltonian PDEs in the degenerate cases},
     \newblock J. Funct. Anal.  \textbf{263 } (2012), no. 2, 511--528.
%


\bibitem{M M 1}
	\newblock M. Martel and F. Merle, 
 	\newblock \emph{Asymptotic stability of solitons for subcritical generalized KdV equations}, 
 	\newblock Arch. Ration. Mech. Anal. \textbf{157}  (2001), no. 3, 219--254.

\bibitem{M M 2}
	\newblock M. Martel and F. Merle, 
 	\newblock \emph{Asymptotic stability of solitons of the subcritical gKdV equations revisited}, 
 	\newblock  Nonlinearity \textbf{18} (2005), no. 1, 55--80.
 	
 \bibitem{M M 3}
	\newblock M. Martel and F. Merle, 
 	\newblock \emph{Asymptotic stability of solitons of the gKdV equations with general nonlinearity}, 
 	\newblock  Math. Ann. \textbf{341} (2008), no. 2, 391--427.

\bibitem{TM1}
	\newblock T. Mizumachi, 
 	\newblock \emph{Large time asymptotics of solutions around solitary waves to the generalized Korteweg-de Vries equations}, 
 	\newblock SIAM J. Math. Anal. \textbf{32} (2001), no. 5, 1050--1080
     
\bibitem{TM2}
	\newblock T. Mizumachi, 
 	\newblock \emph{Stability of line solitons for the KP-II equation in $\R^2$}, 
 	\newblock Mem. Amer. Math. Soc. \textbf{238} (2015), no. 1125, vii+95 pp. 
 	
 \bibitem{M T}
	\newblock T. Mizumachi and N. Tzvetkov, 
 	\newblock \emph{Stability of the line soliton of the KP-II equation under periodic transverse perturbations}, 
 	\newblock Math. Ann. \textbf{352} (2012), no. 3, 659--690. 
 	
\bibitem{M P}
	\newblock  L. Molinet and D. Pilod,
	\newblock  \emph{Bilinear Strichartz estimates for the Zakharov-Kuznetsov equation and applications}, 
	\newblock  Ann. Inst. H. Poincare Anal. Non Lineaire \textbf{32} (2015), no. 2, 347--371.

\bibitem{RN}
	\newblock  R. Nagel (ed.),
	\newblock  \emph{One Parameters Semigroups of Positive Operators}, 
	\newblock  Lecture Notes in Math., {\bf 1184}, Springer-Verlag, Berlin, 1984.

\bibitem{MO}
     \newblock M. Ohta, 
     \newblock \emph{Instability of bound states for abstract nonlinear Schr\"odinger equations},
     \newblock J. Funct. Anal.  \textbf{261} (2011), no. 1, 90-110.
 
 \bibitem{P W} 
\newblock R. Pego and M. I. Weinstein,
\newblock \emph{Eigenvalues, and instabilities of solitary waves},
\newblock Philos. Trans. Roy. Soc. London Ser. A \textbf{340} (1992), no. 1656, 47--94. 
 	
\bibitem{P W 2} 
\newblock R. Pego and M. I. Weinstein,
\newblock \emph{Asymptotic stability of solitary waves},
\newblock Comm. Math. Phys. \textbf{164} (1994), no. 2, 305--349. 

\bibitem{R V}
	\newblock F. Ribaud and S. Vento, 
	\newblock \emph{Well-posedness results for the three-dimensional Zakharov-Kuznetsov equation}, 
	\newblock  SIAM J. Math. Anal. \textbf{44} (2012), no. 4, 2289--2304.
	
\bibitem{R T 0}
	\newblock F. Rousset and N. Tzvetkov, 
	\newblock \emph{Transverse nonlinear instability of solitary waves for some Hamiltonian PDE's}, 
	\newblock J. Math. Pures. Appl. \textbf{90} (2008) 550--590.

\bibitem{R T 1}
	\newblock F. Rousset and N. Tzvetkov, 
	\newblock \emph{Transverse nonlinear instability for two-dimensional dispersive models}, 
	\newblock Ann. I. Poincar\'e-AN \textbf{26} (2009) 477--496.

\bibitem{R T 2} 
\newblock F. Rousset and N. Tzvetkov,
\newblock \emph{A simple criterion of transverse linear instability for solitary waves},
\newblock Math. Res. Lett., \textbf{17} (2010), no. 1, 157--169.

\bibitem{R T 3} 
\newblock F. Rousset and N. Tzvetkov,
\newblock \emph{Stability and instability of the KdV solitary wave under the KP-I flow},
\newblock Comm. Math. Phys., \textbf{313} (2012), no. 1, 155--173.


\bibitem{S S} 
\newblock J. Shatah and W. Strauss,
\newblock \emph{Spectral condition for instability},
\newblock Contemp. Math., \textbf{255} (2000), 189--198.


\bibitem{V A}
     \newblock J. Villarroel and M. J. Ablowitz, 
     \newblock \emph{On the initial value problem for the KPII equation with data that do not decay along a line},
     \newblock Nonlinearity, \textbf{17} (2004), no. 5, 1843--1866. 

     \bibitem{MIW1}
     \newblock M. I. Weinstein, 
     \newblock \emph{Modulational stability of ground states of nonlinear Schrodinger equations},
     \newblock SIAM J. Math. Anal., \textbf{16} (1985), 472--491.

     \bibitem{MIW2}
     \newblock M. I. Weinstein, 
     \newblock \emph{Lyapunov stability of ground states of nonlinear dispersive evolution equations},
     \newblock  Comm. Pure Appl. Math., \textbf{39} (1986), no. 1, 51--67.     
     
\bibitem{YY1}
     \newblock Y. Yamazaki, 
     \newblock \emph{Transverse instability for a system of nonlinear Schr\"odinger equations},
     \newblock  Discrete Contin. Dyn. Syst. Ser. B \textbf{19} (2014), no.2, 565--588.
     
\bibitem{YY2}
     \newblock Y. Yamazaki, 
     \newblock \emph{Stability of line standing waves near the bifurcation point for nonlinear Schr\"odinger equations},
     \newblock  Kodai Math. J. \textbf{38} (2015), no. 1, 65--96.
     
     \bibitem{YY3}
     \newblock Y. Yamazaki, 
     \newblock \emph{Transverse instability for nonlinear Schr\"odinger equation with a linear potential},
     \newblock Adv. in Differential Equations \textbf{21} (2016), no. 5-6, 429--462.



\bibitem{Z K}
     \newblock V. E. Zakharov and E. A. Kuznetsov, 
     \newblock \emph{On three dimensional solitons},
     \newblock  Sov. Phys.-JETP,  \textbf{39} (1974), no.2, 285--286.


\end{thebibliography}
\end{document}